\documentclass[12pt, a4paper, oneside, openany]{book}
\usepackage[utf8]{inputenc}
\usepackage[T1]{fontenc}
\usepackage{geometry}
\usepackage[hidelinks]{hyperref}
\usepackage{titlesec}
\usepackage{parskip}
\usepackage{amsmath, amssymb, amsfonts, mathtools}
\usepackage{amsthm}
\usepackage{graphicx}
\usepackage{booktabs}
\usepackage{enumitem}
\usepackage{xcolor}
\usepackage{tcolorbox}
\tcbuselibrary{breakable}

\usepackage{titlesec}
\titleformat{\paragraph}[block]
{\normalfont\normalsize\bfseries}
{\theparagraph}{1em}{}
\titlespacing*{\paragraph}{0pt}{3.25ex plus 1ex minus .2ex}{0.5em}

\usepackage{tikz}
\usepackage{pgfplots}
\pgfplotsset{compat=1.18}
\newcommand{\R}{\mathbb{R}}
\newcommand{\Pbb}{\mathbb{P}}
\newcommand{\E}{\mathbb{E}}
\newcommand{\dd}{\,\mathrm{d}}

\newcommand{\Xspace}{\R^n}
\newcommand{\Uspace}{\R^m}
\newcommand{\Law}{\mathcal{L}}

\newtheorem{theorem}{Theorem}[section]

\theoremstyle{remark}
\newtheorem{remark}[theorem]{Remark}

\tcbset{
  goalbox/.style={
    colback=blue!5!white,
    colframe=blue!60!black,
    breakable
  },
  algorithmbox/.style={
    colback=gray!3,
    colframe=black,
    breakable
  },
  notebox/.style={
    colback=gray!3,
    colframe=gray!65!black,
    breakable
  }
}

\titleformat{\chapter}[display]
  {\normalfont\huge\bfseries}
  {\chaptertitlename\ \thechapter}
  {18pt}
  {\Huge}
\titleformat{\section}{\Large\bfseries}{\thesection}{0.75em}{}[\titlerule]
\titleformat{\subsection}{\large\bfseries}{\thesubsection}{0.6em}{}
\titleformat{\subsubsection}{\normalsize\bfseries}{\thesubsubsection}{0.5em}{}
\title{\textbf{Notes on Generative Modeling for \\ Feedback Control and Planning}}
\author{Karthik Elamvazhuthi}

\begin{document}

\begin{titlepage}
	\centering
	\vspace*{\fill}
	
	\begin{minipage}[c]{0.48\textwidth}
		\raggedright
		{\Huge\bfseries
			Notes on Generative Modeling for Feedback Control
			and Planning\par}
		
		\vspace{1.5cm}
		
		{\Large Karthik Elamvazhuthi\par}
		
		\vspace{0.5cm}
		
		{\large \today\par}
	\end{minipage}
	\hfill
	\begin{minipage}[c]{0.47\textwidth}
		\centering
		\includegraphics[width=\linewidth]{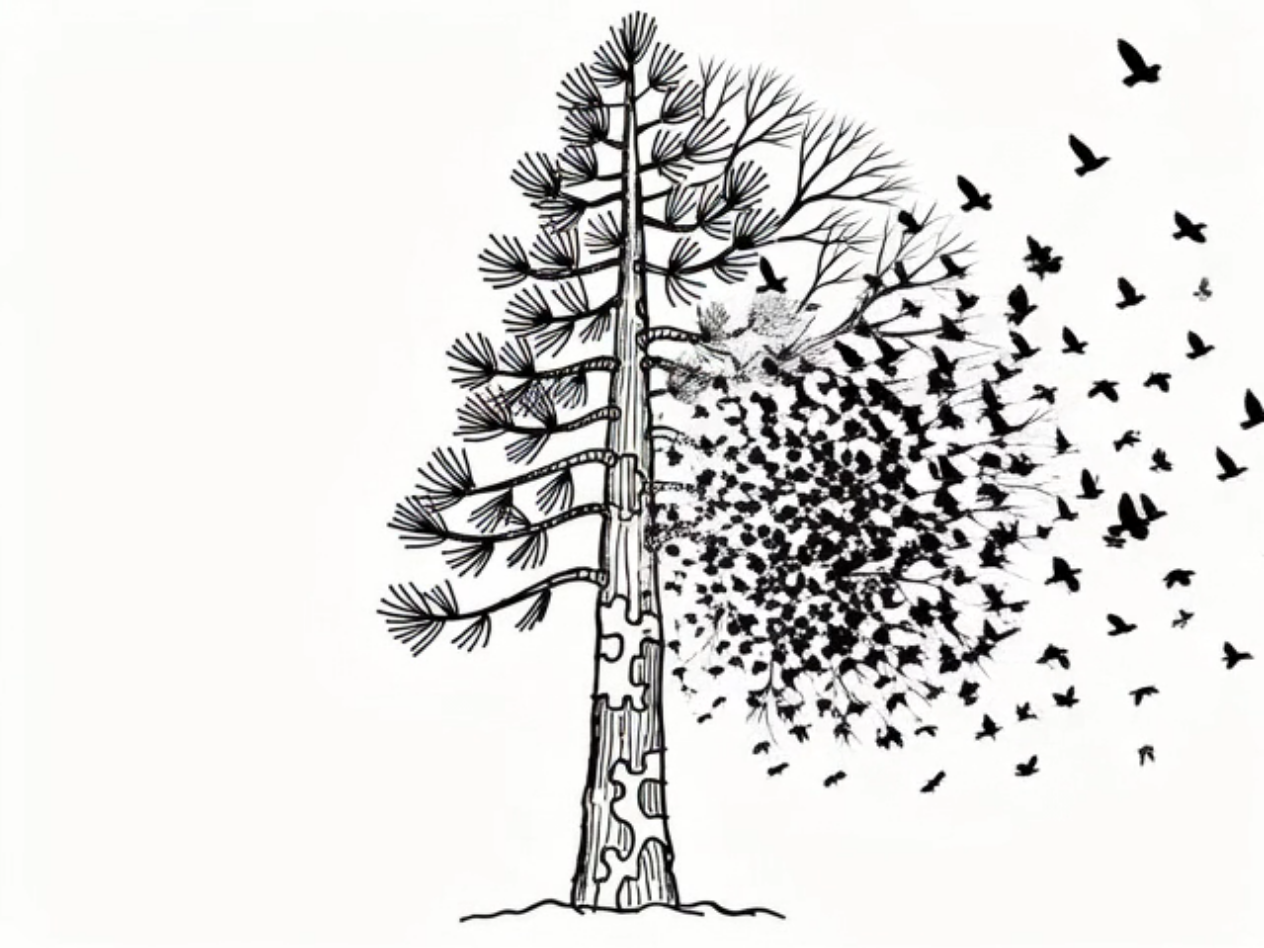}
	\end{minipage}
	
	\vspace*{\fill}
\end{titlepage}

\tableofcontents
\newpage

\chapter{Introduction}
\label{chap:foundations}

\section{From Generative Modeling to Control}

The goal of this series of notes is to frame control as a generative modeling problem and to incorporate modern diffusion-style deep learning methods into planning and control.

Let us first go over generative modeling. Let $X$ be a random variable with density $\rho_X(x)$, and let $Y$ be another random variable with density $\rho_Y(y)$. The generative modeling problem is to find $G_\theta$ such that
\[
    Y \approx G_\theta(X).
\]

One can also use \textit{measure pushforward notation} to pose this problem. If $\mu_X$ denotes the law of $X$, the map $G_\theta:\R^n\to\R^n$ induces the pushforward measure $(G_\theta)_{\#}\mu_X$, defined by
\[
    (G_\theta)_{\#}\mu_X(A)
    =
    \mu_X\bigl(G_\theta^{-1}(A)\bigr)
\]
for sets $A$. When the measures admit densities, we use the same notation informally for the corresponding densities. The problem then we are interested is the following.

\begin{tcolorbox}[goalbox,title=Goal (Transport)]
Find $G_\theta$ such that
\[
    (G_\theta)_{\#}\,\mu_X = \mu_Y,
\]
\end{tcolorbox}

\paragraph{Control as generative modeling}

In the control setting the question transforms as follows. Let
\[
    \frac{\dd X}{\dd t}=f(X,U), \qquad X_0\sim\mu_0.
\]
Can we find a (possibly time-dependent) feedback law $U_\theta(t,x)$ such that the closed-loop dynamics
\[
    \frac{\dd X}{\dd t}=f\bigl(X,U_\theta(t,X)\bigr), \qquad X_0\sim\mu_0
\]
satisfy the terminal objective $X(T) \approx Y$?

We can again pose this problem in the language of probability measures or densities. Let $\Phi_U: \R^n \to \R^n$ be the flow map induced by the controlled dynamics. Given the ODE
\[
    \dot x = f\bigl(x, u(t,x)\bigr), \qquad x(0) = x_0,
\]
define $\Phi_U(x_0) := x(T)$. The problem then becomes:

\begin{tcolorbox}[goalbox,title=Goal (Control as Transport)]
Find $u(t,x)$ such that
\[
    (\Phi_u)_{\#}\,\mu_0 = \mu_\star.
\]
\end{tcolorbox}

At first, it might seem like we have only made the problem harder. What once looked like a simple state-to-state control problem has now been transformed into a probabilistic task of transforming one density into another. A benefit of this viewpoint is that one can view control as a dynamically constrained sampling problem. Since many modern machine learning techniques such as diffusion models and flow matching are techniques for sampling, one can integrate these ideas with control. One especially interesting instance of this integration is using time-reversals of diffusive processes to steer controls systems to a target state.

Another especially appealing aspect of  this perspective is that it simplifies certain issues when it comes to analysis. Once the closed-loop vector field is fixed, the evolution equation for the probability density is linear in the density, even when the underlying state dynamics are nonlinear. Thus, lot of the analysis, which we go over only formally in these notes, reduces to linearity considerations.

We will explore several ways to adapt existing generative modeling methods to the control setting with the following road map:

\begin{tcolorbox}[goalbox,title=Roadmap]
	\begin{enumerate}
		\item \textbf{Measure transport and the continuity equation:} the framework for going from evolution of trajectories to  evolution of probability laws.
		\item \textbf{Flow matching:} feedback synthesis from families of feasible open-loop trajectories, and the noising-denoising principle of control.
		\item \textbf{Reachable-set sampling:} revisit  randomization strategies to sample from reachable sets, and introduce sample efficient alternatives.
		\item \textbf{Controlled normalizing flows:} likelihood-based approach to transport control systems, with application to reachability.
		\item \textbf{Fokker--Planck equations:} dynamics of laws of diffusion processes under white noise forcing and designing control laws to stabilize them to targets.
		\item \textbf{Denoising diffusion:} time reversal of diffusion processes for feedback control.
		\item \textbf{Sub-Laplacian transport:} Connecting analytical properties of second-order differential operators to controllability properties and transport.
	\end{enumerate}
\end{tcolorbox}


We will be restricting our attention in these notes to the continuous-time, continuous-state setting, which provides the cleanest path for adapting modern generative modeling methods from machine learning to control. Discrete-time and discrete-state formulations are certainly possible, but will not be our focus here. There is, perhaps, a mild anxiety among control theorists that the machine-learning bus was missed by not sufficiently embracing the discrete. The prominence of Markov decision process formulations in reinforcement learning can reinforce this impression. These notes take a somewhat unapologetic position in the opposite direction. Rather than discretizing control theory, we will try to ``go full continuum''. While the generative modeling community has embraced measures, partial differential equations, and diffusion operators, it is time for the forces of control and robotics to follow suit. Despite this, no reprieve from the curse of dimensionality is promised. What can be expected, however, for readers willing to get their hands dirty with the algorithms themselves, is an appreciation of the curse from the ``measure-theoretic viewpoint".

The tone of these notes has intentionally been kept informal, with an emphasis on algorithms, the formal introduction of the relevant equations, and the algebraic manipulations needed to work with them. Readers seeking a more rigorous treatment will find precise statements, assumptions, and proofs in the cited references (See the ``Further reading'' subsections). A future version of these notes will incorporate more formal theorems and proofs, closer in style to a traditional control-theory textbook.

\section{Disclaimer}

Generative modeling has become very popular in robotics. Students in this application area might worry if they don't learn the trade they might not be employable. However, the author feels obligated to mention that these set of notes significantly depart from how generative modeling is used in these works. For instance, the dominant framework is to consider generative modeling on entire spaces of {\it trajectories of state or action-state pairs} \cite{janner2022planning,chi2025diffusion}, instead
of the state. Therefore, planning amounts to sampling from this space, satisfying constraints or optimizing for rewards arising from the problem context. This leads to higher dimensionality when applied
to control synthesis, since the set of trajectories lives in a
higher dimensional space in comparison to the set of states. In this sense, this viewpoint is closer in spirit to
path-integral control \cite{kappen2005path} and the broader control-as-inference
literature \cite{levine2018reinforcement}.

In contrast, the approach proposed in these set of notes is {\it Lagrangian} in flavor, as described in the previous section. Control is to be viewed as generative modeling on the {\it state-space} with dynamical constraints on their evolution.

\section{Notions and Notations}

Throughout the notes we use notation some familiar, and some maybe unfamiliar, to practitioners of control. 

\paragraph{Control Systems}

The state is $x\in\Xspace$, the control is $u\in\Uspace$, and
a general controlled system is written
\begin{equation}
	\dot x=f(x,u).
\end{equation}
We will usually work with control affine system, 
\begin{equation}
		f(x,u)
		=
		f_0(x)+G(x)u
		=
		f_0(x)+\sum_{i=1}^m g_i(x)u_i,
	\label{eq:global-control-affine}
\end{equation}
where $f_0$ is the drift, $g_i$ are the control vector fields, and
$G(x)=[g_1(x)\ \cdots\ g_m(x)]$ is the input matrix.
The system is called {\it driftless} when $f_0=0$.

When needed, $x_t^u$ denotes the trajectory generated by
a particular control $u$. The reachable set at time $T$ is
\begin{equation}
	\mathcal R_T(x_0)
	:=
	\left\{
	x_T^u:
	x_0^u=x_0,\;
	u \text{ is admissible}
	\right\}.
\end{equation}

We say that the
system is {\it controllable} on $\Xspace$ if, for every pair
$x_0,x_1\in\Xspace$, there exist a finite time $T>0$ and an
admissible control $u$ such that
\begin{equation}
	x^u(0)=x_0,
	\qquad
	x^u(T)=x_1.
\end{equation}

\paragraph{Probability}

Lowercase letters such as $x$ denote individual states,
while uppercase letters such as $X$ denote random states.
The notation $ X\sim\mu,		\mu=\Law(X),$
means that $X$ has probability law $\mu$: $\mathbb{P}(X \in A) = \int_{A} \mu$ .
For a time-dependent random state, we write
$\mu_t=\Law(X_t)$. When this law admits a density with respect
to Lebesgue measure, we write
\begin{equation}
	\mu_t(\dd x)=\rho_t(x)\dd x,
	\qquad
	\rho_t\geq0,
	\qquad
	\int_{\Xspace}\rho_t(x)\dd x=1.
\end{equation}

A probability law might not have a density. For example,
$\delta_{x_\star}$ denotes the {\it Dirac measure} concentrated
at $x_\star$. We use $\mu_\star$ for a
target law and $\rho_\star$ for its density when one exists. The support $\operatorname{supp}\mu$ consists of those points
whose every neighborhood has positive $\mu$-probability.

\paragraph{Differential operators.}

For a scalar function $\varphi$,
$\nabla\varphi=(\partial_{x_1}\varphi,\dots,
\partial_{x_n}\varphi)^{\mathsf T}$ is its gradient.
For a vector field $F$, its divergence is
$\nabla\cdot F=\sum_{j=1}^n\partial_{x_j}F_j$.
The Laplacian is
$\Delta\varphi=\nabla\cdot\nabla\varphi$.

For the control directions $g_i$, we use the directional
operators and their formal adjoints
\begin{equation}
	\mathcal Y_i\varphi:=g_i\cdot\nabla\varphi,
	\qquad
	\mathcal Y_i^*q:=-\nabla\cdot(g_iq).
	\label{eq:global-directional-operators}
\end{equation}
The word {\it adjoint} refers to the integration-by-parts
identity
\begin{equation}
	\int_\Omega q\,\mathcal Y_i\varphi\,\dd x
	=
	\int_\Omega\varphi\,\mathcal Y_i^*q\,\dd x,
\end{equation}
for compactly supported functions $\phi$.
We similarly write $\mathcal Y_0\varphi=f_0\cdot\nabla\varphi$
for differentiation along the drift.

The horizontal gradient collects derivatives along the
admissible control directions:
\begin{equation}
	\nabla_H\varphi
	:=
	\begin{pmatrix}
		\mathcal Y_1\varphi\\
		\vdots\\
		\mathcal Y_m\varphi
	\end{pmatrix}
	=
	G^{\mathsf T}\nabla\varphi.
\end{equation}
Its associated sub-Laplacian is
$\Delta_H=-\sum_{i=1}^m\mathcal Y_i^*\mathcal Y_i$.
For coordinate vector fields $g_i=e_i$, these reduce to
the ordinary gradient and Laplacian.

These conventions are reused throughout the notes.

\chapter{Measure Transport and the Continuity Equation}
\label{chap:measure-transport}

Chapter~\ref{chap:foundations} posed feedback control as the problem of transporting an initial probability law to a desired terminal law. We now introduce the continuity equation, which will be used in subsequent chapters.

\section{From Trajectories to Densities}
\label{sec:trajectories-to-densities}

Consider the ordinary differential equation
\begin{equation}
    \dot{x}=f(x),
    \qquad
    x(0)=x_0,
\end{equation}
where $f:\R^n\to\R^n$ is a smooth vector field. Under standard regularity assumptions on $f$, there exists a unique solution $x(t)$ to the ODE.

Now suppose that we do not know $x_0$ exactly. Instead, the initial condition is distributed according to a probability density $\rho_0(x)$. Given a set $A\subseteq\R^n$, the probability that the initial condition lies in $A$ is
\[
    \Pbb(X_0\in A)
    =
    \int_A \rho_0(y)\dd y.
\]

We are solving the same ODE for different initial conditions, with different probabilities. What, then, is the evolution of the probability
\[
    \Pbb(X_t\in A)
    =
    \int_A \rho_t(y)\dd y,
\]
where $\rho:[0,T]\times\R^n\to[0,\infty)$ is the evolving density?

\begin{figure}[htbp]
    \centering
    \IfFileExists{ode_to_density.pdf}{%
        \includegraphics[width=\textwidth]{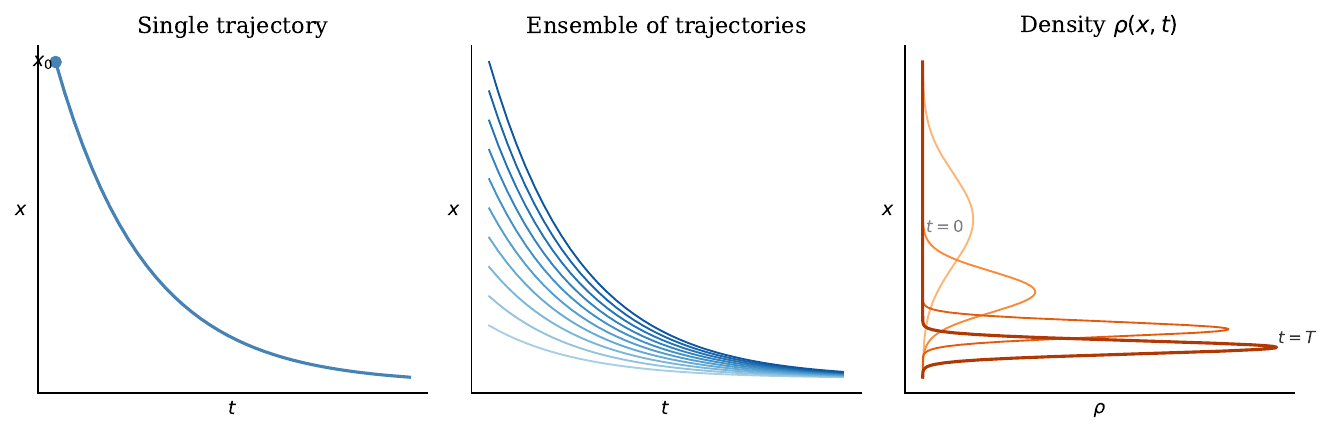}%
    }{%
        \fbox{\parbox[c][2.2in][c]{0.92\textwidth}{\centering
        Figure placeholder: add \texttt{ode\_to\_density.pdf} beside the main \LaTeX{} file.}}%
    }
    \caption{From a single trajectory to an evolving density.
    \textit{Left:} A single solution $x(t)$ of $\dot{x}=f(x)$ from a known initial condition $x_0$.
    \textit{Center:} An ensemble of trajectories from initial conditions drawn from $\rho_0$.
    \textit{Right:} The probability density $\rho_t$ at successive times; as trajectories compress, the density narrows and its peak grows.}
    \label{fig:ode-to-density}
\end{figure}

To derive the equation for $\rho_t$, we first define the pushforward operation. Any measurable map $T:\R^n\to\R^n$ induces a transformation of probability measures. If $\rho_1$ and $\rho_2$ are densities and
\[
    T_{\#}\rho_1=\rho_2,
\]
then, by definition,
\[
    \int_{T^{-1}(A)}\rho_1(x)\dd x
    =
    \int_A\rho_2(x)\dd x
\]
for every set $A\subseteq\R^n$. In probabilistic notation, if $X$ has density $\rho_1$ and $Y=T(X)$, then $Y$ has density $\rho_2$.

An equivalent and particularly convenient characterization uses a smooth test function $\varphi:\R^n\to\R$:
\begin{equation}
    \int_{\R^n}\varphi(T(x))\rho_1(x)\dd x
    =
    \int_{\R^n}\varphi(x)\rho_2(x)\dd x.
    \label{eq:pushforward-test-function}
\end{equation}
This is the pushforward, or change-of-variables, identity that we will use to derive the continuity equation.

Let $\Phi_t:\R^n\to\R^n$ be the flow map of the differential equation, defined by
\[
    \frac{\dd}{\dd t}\Phi_t(x)
    =
    f(\Phi_t(x)),
    \qquad
    \Phi_0(x)=x.
\]
The evolving density is the pushforward of the initial density:
\[
    (\Phi_t)_{\#}\rho_0=\rho_t.
\]
The corresponding operator that transports densities is commonly called the Perron--Frobenius operator. 

Let $\varphi\in C_c^\infty(\R^n)$ be a smooth compactly supported test function. On the one hand,
\[
    \frac{\dd}{\dd t}
    \int_{\R^n}\varphi(x)\rho_t(x)\dd x
    =
    \int_{\R^n}\varphi(x)\partial_t\rho_t(x)\dd x.
\]
On the other hand, using \eqref{eq:pushforward-test-function},
\begin{align*}
    \frac{\dd}{\dd t}
    \int_{\R^n}\varphi(\Phi_t(x))\rho_0(x)\dd x
    &=
    \int_{\R^n}
    \nabla\varphi(\Phi_t(x))\cdot f(\Phi_t(x))
    \rho_0(x)\dd x\\
    &=
    \int_{\R^n}
    \nabla\varphi(x)\cdot f(x)\rho_t(x)\dd x.
\end{align*}
The first equality follows by differentiating under the integral and applying the chain rule; the second follows from the pushforward identity.

Therefore,
\[
    \int_{\R^n}\varphi(x)\partial_t\rho_t(x)\dd x
    =
    \int_{\R^n}\nabla\varphi(x)\cdot f(x)\rho_t(x)\dd x.
\]
Integration by parts gives
\[
    \int_{\R^n}\varphi(x)\partial_t\rho_t(x)\dd x
    =
    -\int_{\R^n}
    \varphi(x)\nabla\cdot\bigl(f(x)\rho_t(x)\bigr)\dd x.
\]
Since this holds for every test function $\varphi\in C_c^\infty(\R^n)$, we obtain the continuity equation:
\begin{equation}
    \boxed{
    \partial_t\rho_t
    +
    \nabla\cdot(f\rho_t)
    =0.
    }
    \label{eq:continuity-equation}
\end{equation}

One can think of $\rho_t$ as a fluid density. Particles at position $x$ move with velocity $f(x)$ and carry probability mass with them. An equivalent derivation begins directly from conservation of mass on an arbitrary control volume.

Consider, for example, the one-dimensional vector field $f(x)=-x$. Then
\[
    \partial_t\rho
    +
    \partial_x(-x\rho)
    =0,
\]
or equivalently
\[
    \partial_t\rho
    =
    \partial_x(x\rho).
\]
The flow is $\Phi_t(x)=e^{-t}x$. If $\rho_0$ is uniform on $[-1,1]$, then
\[
    \rho_t(x)
    =
    \frac{e^t}{2}\,
    \mathbf{1}_{[-e^{-t},e^{-t}]}(x).
\]
Thus the support contracts toward zero while the density height increases so that total probability remains one, exactly as the ensemble picture in Figure~\ref{fig:ode-to-density} suggests.

A useful property of the continuity equation is that, even if the ODE $\dot{x}=f(x)$ is nonlinear, the PDE is linear in the density $\rho_t$ once the vector field $f$ is fixed. 

\section{Convexifying Nonlinear Control Problems}

{\it This section can be skipped in the first reading, if the goal of the reader is to get to the generative modeling part of the notes. It is mainly an advertisement of the measure formulation from an analytical point of view.}

An immediate application of the continuity equation to control problems is reinterpretation of optimal control problems from the density evolution framework. 

Consider the optimal control problem
\begin{equation}
	\begin{aligned}
		\inf_{x,u}\quad
		&\int_0^T
		\bigl(L(x(t))+M(u(t))\bigr)\dd t,\\
		\text{subject to}\quad
		&\dot x(t)
		=
		f_0(x(t))
		+
		\sum_{i=1}^m u_i(t)g_i(x(t))
		=:f(x(t),u(t)),\\
		&x(0)=x_0.
	\end{aligned}
	\label{eq:classical-trajectory-optimization}
\end{equation}
For general nonlinear vector fields, this problem is
non-convex even when $L$ and $M$ are convex. The difficulty
comes from the nonlinear dynamics, which constrain the
admissible state and control trajectories.

This is the classical formulation.  We choose a control to
shape the evolution of a single trajectory. Next, we
translate this viewpoint into a problem over probability
distributions.

Suppose the initial state is distributed according to a
density $\rho_0$, and we apply a feedback control $u(t,x)$.
As discussed earlier, the evolving density $\rho_t$
satisfies the continuity equation
\begin{equation}
	\begin{aligned}
		\partial_t\rho_t
		+
		\nabla\cdot
		\bigl(f(x,u(t,x))\rho_t\bigr)
		&=0,\\
		\rho_{t=0}&=\rho_0.
	\end{aligned}
	\label{eq:density-control-continuity}
\end{equation}
We can now treat this continuity equation as our state
equation instead of the ODE.

The shift toward generative modeling is that, instead
of controlling one state $x(t)$, we control how the whole
distribution $\rho_t$ moves through the state space.
A fixed initial state corresponds to the initial law
$\delta_{x_0}$. Working with densities instead describes
 a whole ensemble of initial states. 
 
 One change in the structure of the problem is immediately clear. Even though the original control system is nonlinear, the new control system {\it bilinear}. That is, for a fixed control $u(t,x)$, the solution of the continuity equation is linear in $\rho$.

When we seek feedback controls $u(t,x)$ instead of
open-loop controls $u(t)$, a natural extension of the
objective is
\begin{equation}
	\inf_{\rho,u}
	\int_0^T\int_{\R^n}
	\bigl(L(x)+M(u(t,x))\bigr)\rho_t(x)
	\dd x\,\dd t,
	\label{eq:density-control-objective}
\end{equation}
subject to \eqref{eq:density-control-continuity}.
This is the {\it expected cost}: we average the state
and control costs against the distribution of states.
Equivalently, when $X_t$ follows the controlled dynamics,
the objective is
\[
\mathbb E \left[ 
\int_0^T 
\bigl(L(X_t)+M(u(t,X_t))\bigr)\dd t
\right].
\]

At this point, the decision variables are the evolving
density $\rho_t$ and the feedback control $u(t,x)$. This
formulation is generally not convex in $(\rho,u)$.

\paragraph{Convexification using flux variables}

There is, however, a useful change of variables that transform the problem into a convex one.
Define the {\it control fluxes}
\begin{equation}
	m_i(t,x)=\rho_t(x)u_i(t,x),
	\qquad i=1,\dots,m,
	\label{eq:control-flux}
\end{equation}
and write $m=(m_1,\dots,m_m)^{\mathsf T}$.
Where $\rho_t>0$, the feedback can be recovered as
$u=m/\rho_t$.

The optimization problem becomes
\begin{equation}
	\begin{aligned}
		\inf_{\rho,m}\quad
		&\int_0^T\int_{\R^n}
		\left[
		L(x)\rho_t(x)
		+
		\rho_t(x)
		M\!\left(\frac{m(t,x)}{\rho_t(x)}\right)
		\right]\dd x\,\dd t,\\
		\text{subject to}\quad
		&\partial_t\rho_t
		+
		\nabla\cdot
		\left(
		f_0(x)\rho_t
		+
		\sum_{i=1}^m g_i(x)m_i
		\right)
		=0,\\
		&\rho_{t=0}=\rho_0,
		\qquad
		\rho_t\geq0,
		\qquad
		\int_{\R^n}\rho_t(x)\dd x=1.
	\end{aligned}
	\label{eq:convex-density-control}
\end{equation}

The function
\[
(\rho,m)\longmapsto
\rho\,M\!\left(\frac{m}{\rho}\right),
\qquad \rho>0,
\]
is called the {\it perspective} of $M$ \cite{combettes2018perspective}. It is jointly
convex in $(\rho,m)$ whenever $M$ is convex. For example,
if $M(u)=\frac12|u|^2$, the control cost becomes
\begin{equation}
	\rho\,M\!\left(\frac{m}{\rho}\right)
	=
	\frac{|m|^2}{2\rho},
	\label{eq:quadratic-flux-cost}
\end{equation}
with value zero at $(\rho,m)=(0,0)$ and $+\infty$
when $\rho=0$ but $m\neq0$.

Moreover, the continuity equation is linear in the
new decision variables. Although $f_0$ and $g_i$ may
depend nonlinearly on $x$, they are fixed coefficients
in this equation. The state cost is also linear in
$\rho$, so convexity of $L$ is not required for this
density formulation.

Thus, for convex $M$, we obtain a convex optimization
problem over densities and control fluxes. A prescribed
terminal density can be imposed through the additional
linear constraint $
\rho_T=\rho_\star$, with the special case $\rho_{\star} = \delta_{x_{\star}}$ when the system has to be steered to a fixed target state $x_{\star}$

This change of variables is central to the dynamical
optimal transport formulation of Benamou and Brenier
\cite{benamou2000computational}, which considers the
fully actuated system $\dot x=u$, quadratic control
cost, and prescribed initial and terminal densities, generalized to control-affine systems in \cite{elamvazhuthi2023dynamical,elamvazhuthi2024benamou}. One can even incorporate control constraints without breaking the convexity structure. For example, let
$u_i(t,x)\in[-\bar u_i,\bar u_i]$. Introducing the same flux
variables $m_i=\rho u_i$, these bounds become
\[
-\bar u_i\rho
\leq m_i
\leq \bar u_i\rho.
\]
Hence, one advantage of the continuity-equation
perspective is that the problem becomes convex in
the density and flux variables. Of course, there is no free lunch. The price of this continuity equation-burger is that
the optimization is now infinite dimensional!

\paragraph{Reachability as a linear problem}
If one did not have a cost function, analytically, the problem simplifies even
further. Suppose we only wish to transport a prescribed initial
density $\rho_0$ to a target density $\rho_1$. 
Thus, the problem becomes a linear feasibility problem:
\begin{align}
	\partial_t\rho
	+\nabla\cdot\left(
	f_0\rho+\sum_{i=1}^m g_i m_i
	\right)
	&=0,\\
	\rho(0,\cdot)&=\rho_0,
	\qquad
	\rho(T,\cdot)=\rho_1,\\
	\rho&\geq0,\\
	-\bar u_i\rho
	\leq m_i
	&\leq\bar u_i\rho,
	\qquad i=1,\dots,m.
\end{align}

\paragraph{Numerical example}

We illustrate the density-based control formulation using the
double-gyre example considered in
\cite{elamvazhuthi2016optimal}. The controlled dynamics on
$\Omega=[0,2]\times[0,1]$ are
\begin{align}
	\dot x
	&=
	-\pi A\sin(\pi F(x,t))\cos(\pi y)+u_1,\\
	\dot y
	&=
	\pi A\cos(\pi F(x,t))\sin(\pi y)
	\partial_x F(x,t)+u_2,
\end{align}
where
\[
F(x,t)
=
\beta\sin(\omega t)x^2
+
\bigl(1-2\beta\sin(\omega t)\bigr)x.
\]
We choose $A=\beta=0.25$ and $\omega=2\pi$, so the
uncontrolled flow is periodic with period one. The drift
generates two interacting gyres, and the controls perturb
this motion to achieve the desired transport.

We consider two cases over the horizon $T=2$. In the first, the
initial distribution is supported uniformly on the left gyre and transported toward a point
near $(1.5,0.5)$, represented numerically by mass concentrated
in a single grid cell. 
In the second case, the initial and target distributions are approximately
uniform on numerically estimated regions of the left
and right gyres. 

Writing the drift as $b(t,x,y)$ and introducing the control
flux $m=\rho u$, we solve
\begin{equation}
	\inf_{\rho,m}
	\frac12\int_0^T\int_\Omega
	\frac{|m(t,x,y)|^2}{\rho_t(x,y)}
	\,\dd x\,\dd y\,\dd t
\end{equation}
subject to
\[
\partial_t\rho+\nabla\cdot(b\rho+m)=0.
\]
The convex problem is discretized on a
$64\times32$ spatial grid with $80$ time intervals.
Particle trajectories are then reconstructed from the
computed total flux and density. The particle trajectories are not simulations of the corresponding Lagrangian trajectories.

\begin{figure}[htbp]
	\centering
	\includegraphics[width=\textwidth]
	{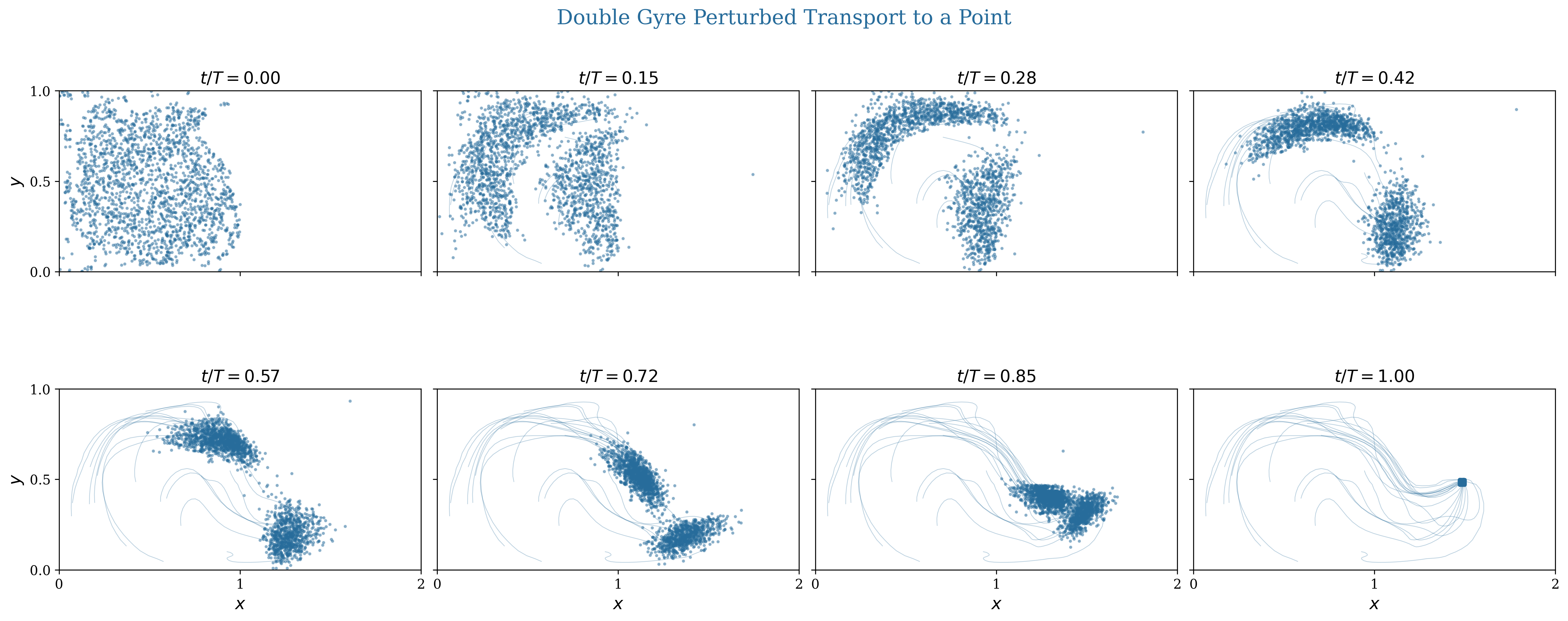}
	\caption{Controlled transport from the left gyre-core
		distribution toward a point target near $(1.5,0.5)$.
		The terminal density is concentrated in 
		a Dirac target.}
	\label{fig:double-gyre-point}
\end{figure}

\begin{figure}[htbp]
	\centering
	\includegraphics[width=\textwidth]
	{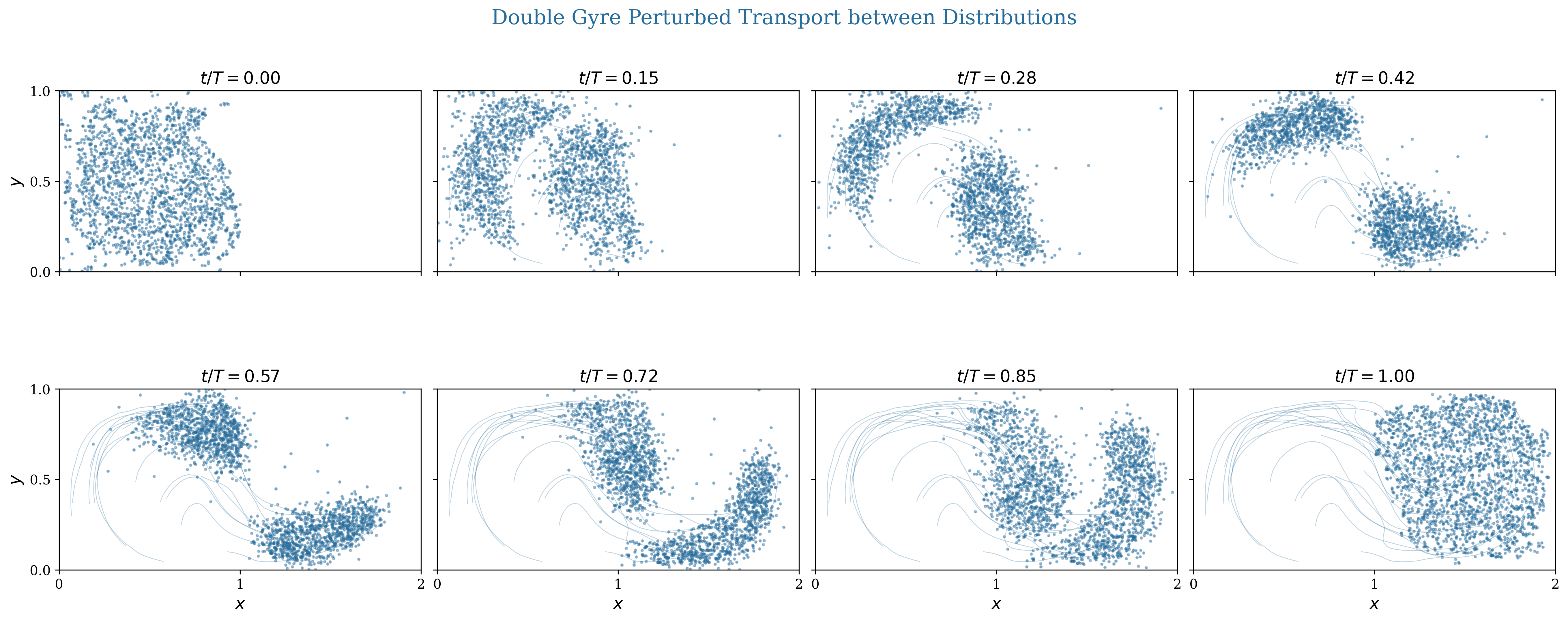}
	\caption{Controlled transport between the two gyre-core
		distributions.}
	\label{fig:double-gyre-distributions}
\end{figure}

Figures~\ref{fig:double-gyre-point} and ~\ref{fig:double-gyre-distributions}
show eight snapshots of the
resulting transport.  Both the solutions show how the optimal transport uses the underlying dynamics to transfer mass, instead of taking straight line trajectories.

\paragraph{Further Reading}
The continuity-equation viewpoint for classical control problems can be
can be found in measure-based formulations of 
control \cite{vinter1993convex,lasserre2008nonlinear, vaidya2008lyapunov,majumdar2014convex}. The Benamou-Brenier formulation of optimal transport \cite{benamou2000computational}, adapted to transport problems with control costs can be found in \cite{chen2016optimal,elamvazhuthi2016optimal,elamvazhuthi2023dynamical,elamvazhuthi2024benamou}. 

\chapter{Flow Matching for Control and Planning}
\label{chap:flow-matching}

The continuity equation from Chapter~\ref{chap:measure-transport} tells us how a autonomous vector field transports a probability law and how some classical control problems can be reframed in this new perspective. In this chapter, we look at an application of the continuity equation introducing {\it Flow matching}, a method introduced in the generative modeling community for sampling from a probability distribution \cite{peluchetti2021nondenoising, lipman2023flow,liu2023rectified,albergo2025stochastic}. It provides us with an approach to start from a family of sampled trajectories, to construct feedback law whose induced probability evolution matches the sampled one.

\section{Converting Open-loop Controls to Feedback Laws}
\label{sec:exact-flow-matching}

Suppose that, for every pair consisting of an initial condition \(x_0\) and a final condition \(x_T\), we have a family of open-loop controls given by a map
\[
S(x_0,x_T) = \bigl(x(t),u(t)\bigr),
\]
such that the system dynamics are
\[
\dot{x}
=
f_0(x)
+
\sum_{i=1}^{m}u_i(t)f_i\bigl(x(t)\bigr)
\equiv f(x,u),
\]
with boundary conditions
\[
x(0)=x_0,
\qquad
x(T)=x_T.
\]

Given a family of open-loop interpolations $S(x_0,x_T)$, can we synthesize a closed-loop control \(u_i(t,x)\) such that
\[
\dot{x}
=
f_0(x)
+
\sum_{i=1}^{m}u_i(t,x)f_i\bigl(x(t)\bigr),
\]
and such that \(x(T)=x_T\) for every \(x_0\in\R^n\)?

In general, this problem is difficult to solve because trajectories generated by different open-loop controls may cross, preventing their exact realization by a single feedback controller.

Although we may not be able to construct a feedback law that reproduces the same individual trajectories as the open-loop controls, we can construct a feedback law that reproduces their probabilistic behavior.

Let $\mu_0$ and $\mu_T$ denote the probability laws of the initial and final conditions, respectively. Define
\[
S_t\bigl(X_0(\omega),X_T(\omega)\bigr)
=
\bigl(X_t(\omega),U_t(\omega)\bigr),
\]
where \(\omega\) is a sample from a probability space.

Let $\mu_t=\Law(X_t)$ denote the probability law induced by the planner. When $\mu_t$ admits a density, we write $\mu_t(\dd x)=\rho_t(x)\dd x$. Our goal is to derive the evolution of $\mu_t$ in terms of the continuity equation. Let $\varphi\in C_c^\infty(\R^n)$
be a smooth, compactly supported test function. Evaluating \(\varphi\) along \(X_t\) gives
\[
\frac{\dd}{\dd t}\varphi\bigl(X_t(\omega)\bigr)
=
\dot{X}_t(\omega)\cdot\nabla \varphi\bigl(X_t(\omega)\bigr).
\]

Integrating with respect to \(\omega\), we obtain
\[
\frac{\dd}{\dd t}
\int_{\R^n}\varphi(x)\mu_t(\dd x)
=
\int_\Omega
\dot{X}_t(\omega)\cdot
\nabla \varphi\bigl(X_t(\omega)\bigr)
\dd\Pbb(\omega).
\]

Using the {\it disintegration theorem}, which decomposes a joint law into a marginal law and conditional laws, for the joint law $(X_t,\dot{X}_t)_\sharp\Pbb$, we write
\[
(X_t,\dot{X}_t)_\sharp\Pbb(\dd x,\dd v)
=
\mu_t(\dd x)\,\eta_{t,x}(\dd v),
\]
where $\eta_{t,x}$ is the conditional distribution of $\dot{X}_t$ given $X_t=x$. Therefore,
\begin{align*}
&\int_\Omega
\dot{X}_t(\omega)\cdot
\nabla \varphi\bigl(X_t(\omega)\bigr)
\dd\Pbb(\omega)
\\
&\qquad=
\int_{\R^n}\int_{\R^n}
v\cdot\nabla \varphi(x)\,
\eta_{t,x}(\dd v)\,
\mu_t(\dd x)
\\
&\qquad=
\int_{\R^n}
\nabla \varphi(x)\cdot
\E\!\left[\dot{X}_t\mid X_t=x\right]
\mu_t(\dd x),
\end{align*}
where $
\E\!\left[\dot{X}_t\mid X_t=x\right]
=
\int_{\R^n}v\,\eta_{t,x}(\dd v),$
sometimes referred to as the barycenter of the measure $\,\eta_{t,x}$ for every $(t,x)$

Because the vector field \(f(x,u)\) depends affinely on \(u\),
\[
\E\!\left[\dot{X}_t\mid X_t=x\right]
=
f\!\left(
x,
\E\!\left[U_t\mid X_t=x\right]
\right).
\]

Define the feedback law
\[
u(t,x)
=
\E\!\left[U_t\mid X_t=x\right].
\]
Then
\[
\frac{\dd}{\dd t}
\int_{\R^n}\varphi(x)\mu_t(\dd x)
=
\int_{\R^n}
\nabla \varphi(x)\cdot f\bigl(x,u(t,x)\bigr)
\mu_t(\dd x).
\]

Formally exchanging differentiation and integration and then applying integration by parts yields the continuity equation that we saw in the first chapter.
\[
\boxed{
\partial_t\rho_t(x)
+
\nabla\cdot
\left(
f\bigl(x,u(t,x)\bigr)\rho_t(x)
\right)
=
0.
}
\]

Thus, by combining the continuity equation with disintegration, we identify a feedback law having the same probabilistic behavior as the open-loop planner. In particular, if $
\mu_T=\delta_{x_\star},$
or equivalently \(X_T(\omega)=x_\star\), then the feedback control law drives the system to \(x_T\).

Next, we want to understand how one can learn the controller \(u(t,x)\). We use a property of \(L^2\) regression toward this end.

It is well known that
\[
f^\star(x)=\E[Y\mid X=x]
\]
minimizes the mean-squared loss
\[
\E\!\left[\lvert Y-f(X)\rvert^2\right].
\]

We can use this principle to learn the feedback controller. Suppose that \(U_t\) is the random open-loop control generated by the planner. Consider the loss
\[
\mathcal{L}(\theta)
=
\int_0^T
\E\!\left[
\left\lVert
U_t-u^\theta(t,X_t)
\right\rVert^2
\right]\dd t.
\]

Its minimizer is the conditional expectation:
\[
u^{\theta^\star}(t,x)
=
\E[U_t\mid X_t=x].
\]

This leads to the following algorithm.

\begin{tcolorbox}[algorithmbox,title=\textbf{Algorithm: Flow Matching for Closed-Loop Synthesis}]
\textbf{Setup.}
Consider the control-affine system
\[
\dot{x}
=
f_0(x)
+
\sum_{i=1}^{m}u_i g_i(x).
\]

\begin{enumerate}[label=\arabic*.]
  \item Choose an initial law $\mu_0$ and a goal state $x_\star$.

  \item Simulate ODE trajectories using the planner
  \[
  (X_t,U_t)=S(X_0,x_\star),
  \]
  where
  \[
  \dot{X}_t
  =
  f_0(X_t)
  +
  \sum_{i=1}^{m}U_i(t)g_i(X_t).
  \]
  Store samples
  \[
  \left\{
  X_{t_k}^{(j)},U_{t_k}^{(j)}
  \right\}_{j=1}^{N}
  \]
  at each time step \(t_k\).

  \item Train a network
  \[
  u_\theta:\R^n\times[0,T]\to\R^m
  \]
  by minimizing the empirical \(L^2\) loss corresponding to
  \[
  \mathcal{L}(\theta)
  =
  \int_0^T
  \E\!\left[
  \left\lVert U_t-u_\theta(X_t,t)\right\rVert^2
  \right]\dd t.
  \]

  \item Apply the learned feedback controller from any initial condition \(x(0)\):
  \[
  \dot{x}(t)
  =
  f_0\bigl(x(t)\bigr)
  +
  \sum_{i=1}^{m}
  u_{\theta,i}\bigl(x(t),t\bigr)
  g_i\bigl(x(t)\bigr).
  \]
  Equivalently,
  \[
  u(t)=u_\theta\bigl(x(t),t\bigr)\in\R^m.
  \]
  As \(t\to T\), the state is expected to approach the goal \(x_T\).
\end{enumerate}
\end{tcolorbox}

\paragraph{Numerical example.}

Consider the following system with two control inputs:
\[
\begin{aligned}
\dot{x}_1(t)&=u_1(t),\\
\dot{x}_2(t)&=u_2(t),\\
\dot{x}_3(t)&=u_1(t)x_2(t).
\end{aligned}
\]
This system is controllable and a canonical example of how oscillatory inputs can be used for path planning of non-holonomic systems of {\it chained form} \cite{jean2014control}. 

An algorithm for steering a system of this chained form from \(x\) to \(y\) proceeds in two stages.

First, we move \((x_1,x_2)\) to \((y_1,y_2)\) using constant controls:
\[
\bigl(u_1(t),u_2(t)\bigr)
=
\left(
\frac{y_1-x_1}{2\pi},
\frac{y_2-x_2}{2\pi}
\right),
\qquad
t\in[0,2\pi].
\]

Second, we use sinusoidal controls to adjust the third coordinate:
\[
\bigl(u_1(t),u_2(t)\bigr)
=
\left(
\sin t,
\frac{y_3-x_3(2\pi)}{\pi}\cos t
\right),
\qquad
t\in(2\pi,4\pi].
\]

During the second stage the first two coordinates return to $(y_1,y_2)$ because the integrals of $\sin t$ and $\cos t$ over $[2\pi,4\pi]$ vanish, while
\[
    x_3(4\pi)-x_3(2\pi)
    =
    \frac{y_3-x_3(2\pi)}{\pi}
    \int_{2\pi}^{4\pi}\sin^2 t\dd t
    =
    y_3-x_3(2\pi).
\]
Hence the resulting trajectory $\omega$ satisfies
\[
\omega(0)=x,
\qquad
\omega(4\pi)=y.
\]

\begin{figure}[htbp]
    \centering
    \includegraphics[width=\textwidth]
        {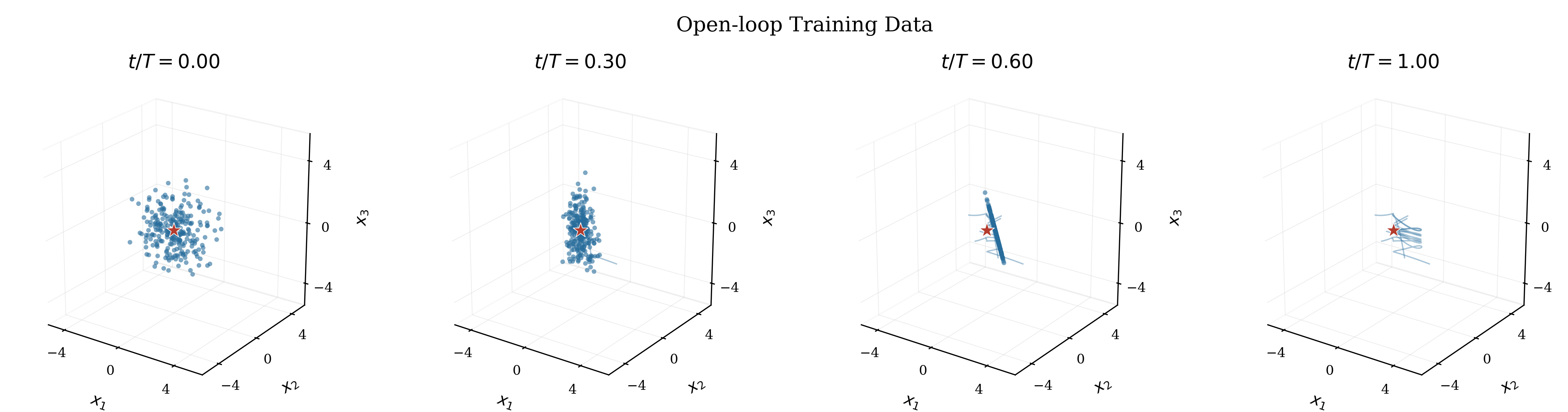}

    \medskip

    \includegraphics[width=\textwidth]
        {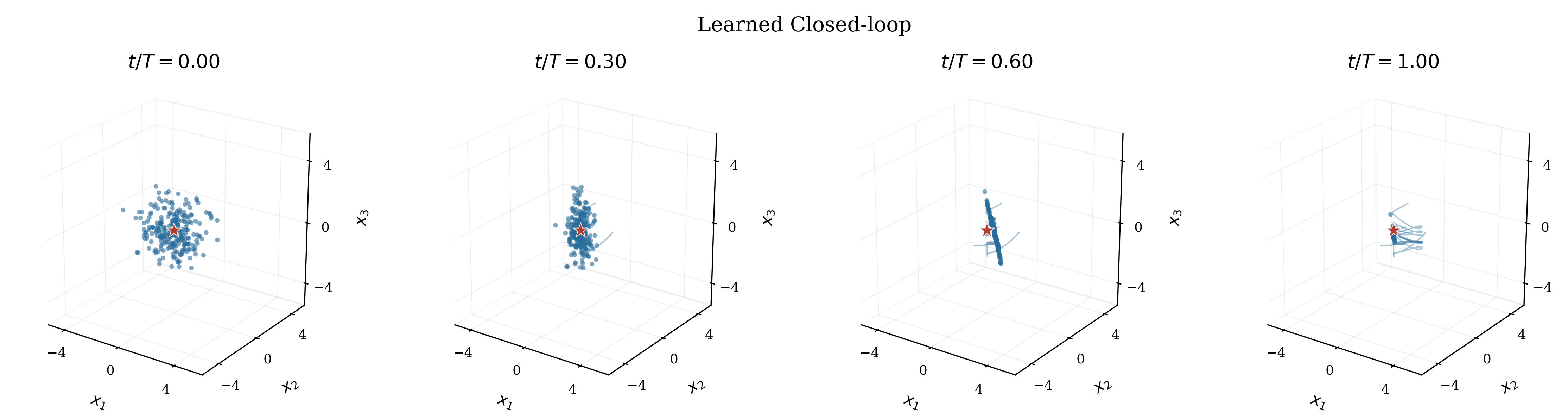}

    \caption{Steering the system
    $\dot{x}_1=u_1$, $\dot{x}_2=u_2$,
    $\dot{x}_3=x_2u_1$ toward the origin.
    Top: open-loop trajectories used to generate data for the flow matching algorithm.
    Bottom: closed-loop trajectories under the learned feedback.
    Each row shows the initial, intermediate, and terminal states
    from left to right. Blue points represent the current states,
    and the red star marks the target. }
    \label{fig:brockett-flow-matching}
\end{figure}

Figure~\ref{fig:brockett-flow-matching} compares the open-loop
trajectories used for training with the trajectories generated
by the learned feedback controller, with the origin as the target.

\section{Noising and Denoising for Feedback Control}
\label{sec:flow-matching-noising}
The previous section assumed that an open-loop planner or interpolator $S(x_0,x_T)$ was already available. We now address the harder question of how to generate a rich family of feasible plans when point-to-point planning is unavailable.

We know from the previous section how flow matching can be used to convert a collection of open-loop control inputs into a feedback control law. This method depends on our ability to solve the open-loop initial--terminal constraint problem
\[
\dot{x}
=
f_0(x)
+
\sum_i u_i(t)g_i(x),
\]
subject to
\[
x(0)=x_0,
\qquad
x(T)=y.
\]

What if we cannot solve this problem analytically for every initial condition \(x_0\) and final condition \(y\)?

In this section, we show how the principle of ``noising and denoising" can be integrated with flow matching. Here, \emph{noising} means deliberately spreading trajectories away from the target, while \emph{denoising} means reversing that transport so that trajectories return toward the target. The idea is to use noisy controls, starting from the goal configuration, to propagate the time-reversed system across its backward reachable set. The corresponding control problem is then to denoise the system from state space back to the goal configuration.

The noising trajectories effectively resolve the initial obstruction that we cannot explicitly solve the point-to-point open-loop trajectory-planning problem. If the noising trajectories are sufficiently rich, then, with positive probability, they explore the entire relevant state space.

\begin{figure}[htbp]
    \centering
    \includegraphics[width=\textwidth]{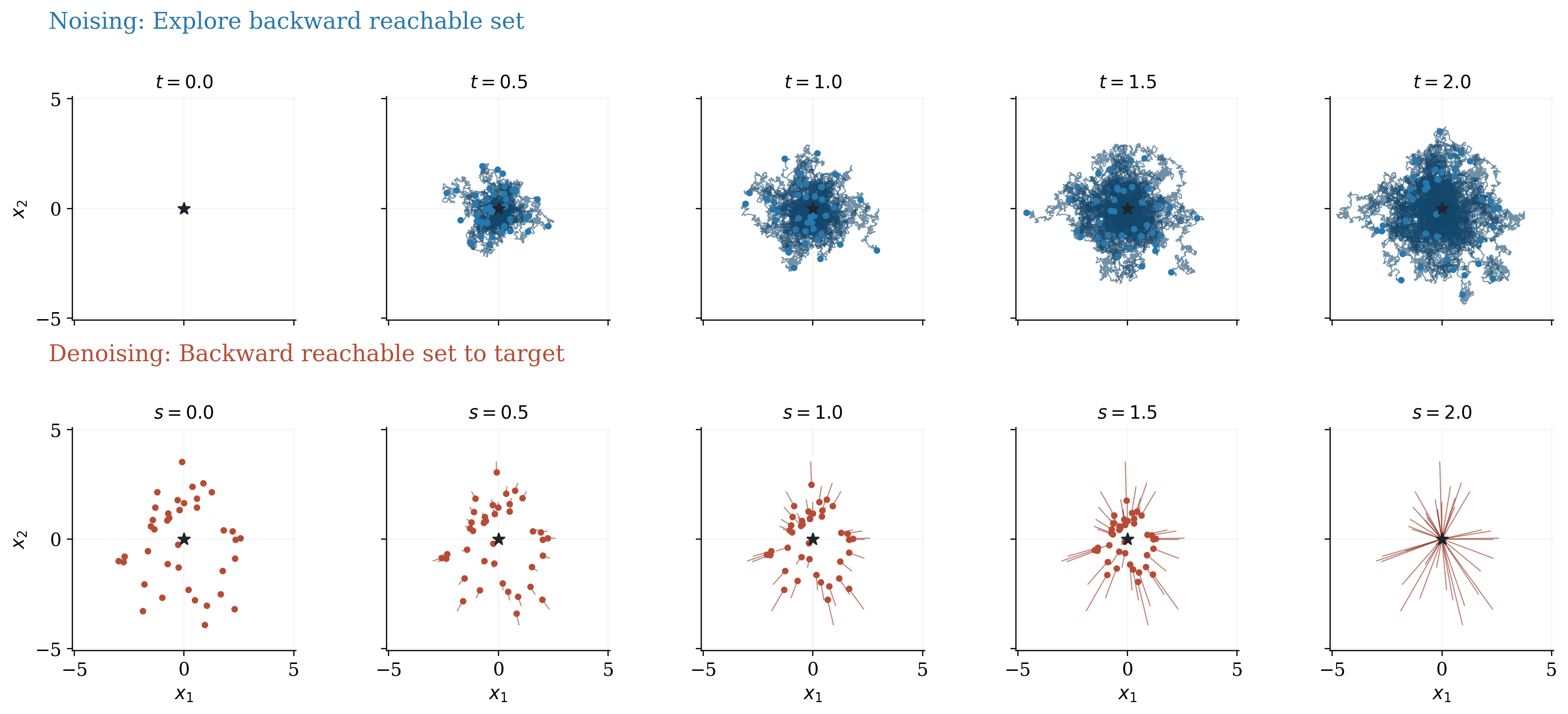}
     Top: Noising explores the backward reachable set from the target, marked as a start.
    Bottom: Deterministic time-reversal or denoising steers the
    noised states back to the target.
    \label{fig:noising-denoising}
\end{figure}

Figure~\ref{fig:noising-denoising} illustrates the two phases:
Noising spreads states from the point target, while a
deterministic flow denoises and transports the resulting
distribution back to the target.

The idea of noising and denoising was originally introduced in the context of SDEs. This is introduced in Chapter \ref{chap:denoising-diffusion}. Here, the aim is to avoid the burden of an stochastic calculus that comes with SDEs. The conceptual idea of noising with flow matching is that one may use \emph{any suitable noising process}, provided that it is consistent with the system dynamics. In the SDE-based framework in Chapter \ref{chap:denoising-diffusion-control}, the standard choice is white noise. However, this choice is not set in any  generative modeling stone. The important requirement is that the noised process can be time-reversed. For ODEs, time reversal is technically easier to derive than for SDEs. Rather than exciting the system with white noise, one may sample controls from any distribution on an admissible control space, provided that the resulting controls explore the reachable set.

The inspiration from denoising models makes flow matching highly versatile. Different choices of control distributions yield different algorithms. One retains the flexibility of denoising diffusion methods for sampling under dynamical and state constraints, while preserving the conceptual simplicity of flow matching.

The abstract flow-matching-based noising--denoising approach is as follows.

\begin{tcolorbox}[algorithmbox,title=\textbf{Algorithm: Noising--Denoising Flow Matching for Closed-Loop Synthesis}]
\textbf{Setup.}
Consider the control-affine system
\[
\dot{x}
=
f_0(x)
+
\sum_{i=1}^{m}u_i g_i(x).
\]

\begin{enumerate}[label=\arabic*.]
  \item Choose a goal state $x_\star$, a horizon $T$, and a
  random control process $U(s)$.

  \item Simulate noising ODE trajectories starting from
  $X_0=x_\star$, where
  \[
  \dot{X}_s
  =
  -f_0(X_s)
  +
  \sum_{i=1}^{m}U_i(s)g_i(X_s).
  \]
  Store samples
  \[
  \left\{
  X_{s_k}^{(j)},U_{s_k}^{(j)}
  \right\}_{j=1}^{N}
  \]
  at each time step \(s_k\).

  \item Train a network
  \[
  u_\theta:\R^n\times[0,T]\to\R^m
  \]
  by minimizing the empirical \(L^2\) loss corresponding to
  \[
  \mathcal{L}(\theta)
  =
  \int_0^T
  \E\!\left[
  \left\lVert U_s-u_\theta(X_s,s)\right\rVert^2
  \right]\dd s.
  \]

  \item Apply the learned feedback controller with time and sign
  reversed, from an initial condition \(x(0)\) in the region
  explored by the terminal noising distribution:
  \[
  \dot{x}(t)
  =
  f_0\bigl(x(t)\bigr)
  -
  \sum_{i=1}^{m}
  u_{\theta,i}\bigl(x(t),T-t\bigr)
  g_i\bigl(x(t)\bigr).
  \]
  Equivalently,
  \[
  u(t)=-u_\theta\bigl(x(t),T-t\bigr)\in\R^m.
  \]
  As \(t\to T\), the state is expected to approach the goal
  \(x_\star\).
\end{enumerate}
\end{tcolorbox}

One possible justification for whether the noising process \(X(t)\) explores the whole space is the following.

Suppose the probability distribution on the set of controls has
{\it support} dense in the space of finite-energy controls
\[
L^2(0,T;\R^m)
:=
\left\{
u:[0,T]\to\R^m 
\;\middle|\;
\int_0^T \lvert u(t)\rvert^2\dd t<\infty
\right\}.
\]
Informally, this means that every $L^2$ neighborhood of every
finite-energy control has positive probability under the
random control law: for every
$u\in L^2(0,T;\R^m)$ and every $\varepsilon>0$,
\[
\mathbb P\!\left(
\lVert U-u\rVert_{L^2(0,T;\R^m)}<\varepsilon
\right)>0.
\]
Under suitable continuity and well-posedness assumptions for the control-to-endpoint map, the support of the terminal noising law is the closure of the endpoints generated by controls in the support of the noising law. If that support is dense in $L^2(0,T;\R^m)$, this terminal support coincides with the closure of the $L^2$-reachable set of the time-reversed system, given by
\begin{equation}
\label{eq:reachableset}
\mathcal R_T^{-}(x_\star)
:=
\left\{
X(T):
\begin{array}{l}
\displaystyle
\dot X(s)=-f_0(X(s))
+\sum_{i=1}^{m}U_i(s)g_i(X(s)),\\[2pt]
X(0)=x_\star,\qquad U\in L^2(0,T;\R^m)
\end{array}
\right\}.
\end{equation}
Thus,
\[
\operatorname{supp}\bigl(\operatorname{Law}(X_T)\bigr)
=
\overline{\mathcal R_T^{-}(x_\star)}.
\]
Equivalently, $\mathcal R_T^{-}(x_\star)$ is the set of initial
states that can be steered to $x_\star$ in time $T$ under the
original system, using the controls $u(t)=-U(T-t)$.

To complete this program one needs to a suitable choice of the noising process. One possible choice of the noising law is to let \(U\) be Brownian motion initialized at zero. The following support property of Brownian motion is useful for ensuring sufficiently rich exploration.

\paragraph{Support property of Brownian } 
Viewing Brownian motion as a random element of path space, the support of associated {\it Wiener measure} is $
C_0([0,T];\R^m)$,
consisting of continuous paths initialized at the origin. Equivalently, every open neighborhood of a path in \(C_0([0,T];\R^m)\) has positive Wiener measure. This property is useful because $C_0([0,T];\R^m)$
is dense in $
L^2(0,T;\R^m).$ Therefore, if the time-reversed system is excited by Brownian motion, it explores the entire backward reachable set \eqref{eq:reachableset} (up to closure).

To remove the restriction that \(U(0)=0\), one may perturb the controls by setting
\[
U(t)=B(t)+I,
\]
where \(B(t)\) is Brownian motion and \(I\) is a time-independent random variable taking values in \(\R^m\). The density property still holds true.

\paragraph{Underactuated system on $SO(3)$.}
The special orthogonal group
\[
SO(3)
=
\left\{
R\in\R^{3\times3}:
R^\top R=I,\ \det R=1
\right\}
\]
is the space of rotations in three dimensions.
Although represented by nine matrix entries, $SO(3)$ is a
three-dimensional manifold.
We test the method on the system
\[
\dot R
=
R\bigl(u_1A_1+u_2A_2\bigr),
\qquad R\in SO(3),
\]
where
\[
A_1=
\begin{pmatrix}
0&0&0\\
0&0&-1\\
0&1&0
\end{pmatrix},
\qquad
A_2=
\begin{pmatrix}
0&0&1\\
0&0&0\\
-1&0&0
\end{pmatrix}.
\]

Although the linearization at the identity is not controllable,
the nonlinear system is controllable. To construct the control law, we initialize the noising process
at the goal configuration $R_\star=I$ and generate trajectories
using bounded controls obtained by applying a $\tanh$ transformation
to Brownian paths. Since the system is driftless, these trajectories
explore its backward reachable set.

\begin{figure}[htbp]
    \centering
    \includegraphics[width=\textwidth]
        {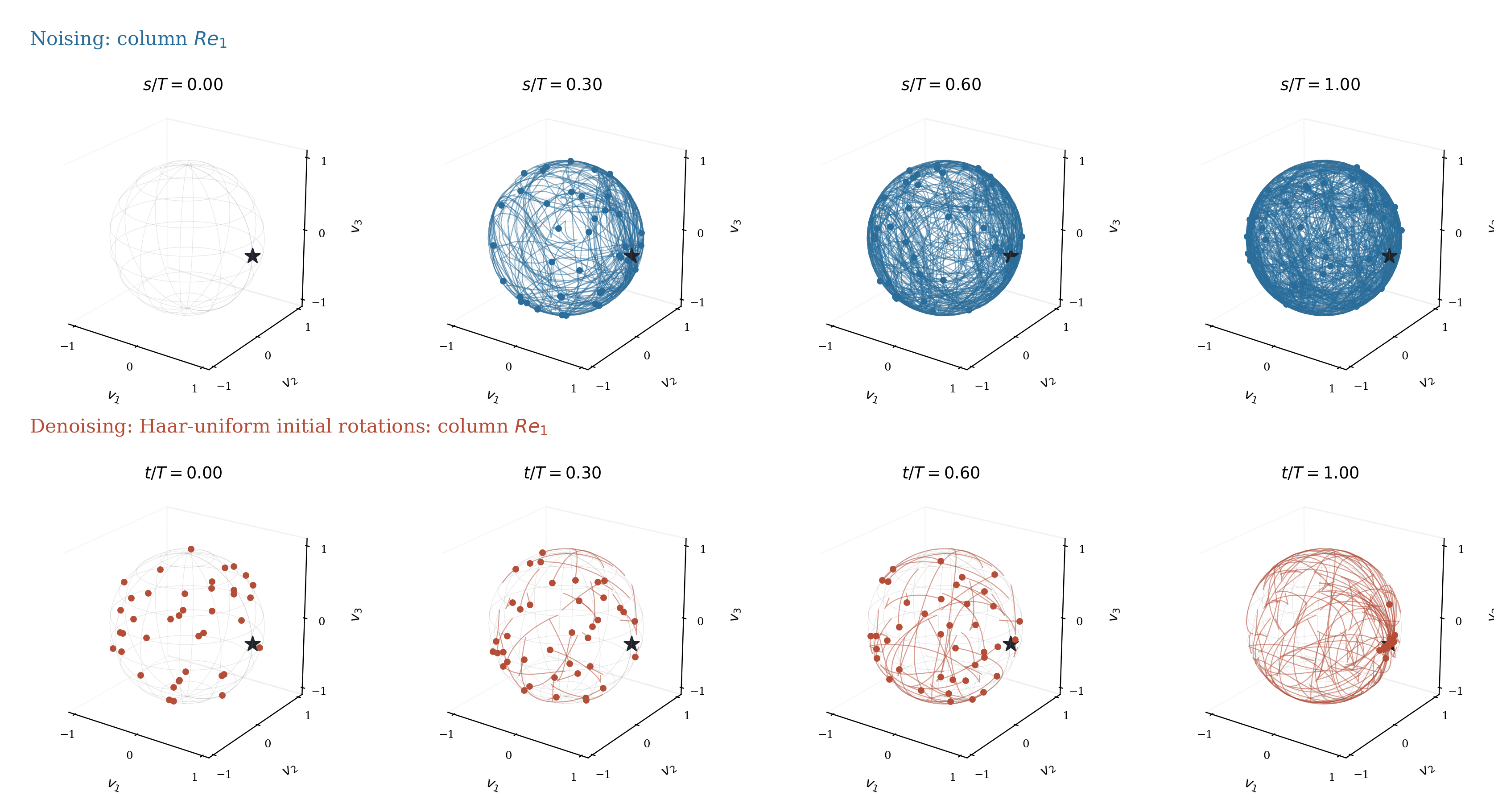}
    \caption{Noising and denoising for an underactuated system
    on $SO(3)$, visualized through the first column $Re_1\in S^2$.
    Top: noising from the identity under bounded Brownian-motion
    controls.
    Bottom: learned denoising from Haar measure.
    Each row progresses from left to right.
    Dots indicate current column vectors, curves show their
    paths, and the star marks $e_1$.}
    \label{fig:so3-noising-denoising}
\end{figure}

We then learn the flow-matched controller by regressing the sampled
controls on the corresponding rotations and times. 
We evaluate this feedback using initial rotations distributed according to the {\it Haar measure}.
Figure~\ref{fig:so3-noising-denoising} illustrates both phases
through the evolution of the first column $Re_1$ on the unit sphere.

The noising--denoising construction suggests that reversing the
control law should steer the state toward $x_\star$. Before
proceeding to the next section, let us provide an informal
clarification of why we expect this to happen.

In practice, we do not know the terminal noising law $\mu_T$
explicitly. We will therefore likely initialize the reverse
system from a different law $\nu_0$. Suppose that both laws
have densities such that
\[
\nu_0(x)\leq C\mu_T(x)
\]
for some finite constant $C$. Thus, the new initialization
does not place mass in regions missed by the noising law,
and its relative weighting is bounded.

To see why this is useful, let $Z_t$ denote the reverse
trajectory under the exact controller. Under the reference
initialization $Z_0\sim\mu_T$, and assuming the flow is
well posed for $t<T$, the time-reversal construction gives
\[
\operatorname{Law}(Z_t)=\mu_{T-t}
\Longrightarrow\delta_{x_\star}
\qquad\text{as }t\uparrow T.
\]
We can change the initial distribution without changing the
controller by reweighting these trajectories:
\[
\frac{\dd\mathbb Q}{\dd\mathbb P}
=
\frac{q_0(Z_0)}{\rho_T(Z_0)}.
\]
Here, $\mathbb P$ is the reference probability measure, while
under $\mathbb Q$ the initial law is $\nu_0$. Every trajectory
still follows the same feedback equation. For every $r>0$,
\[
\mathbb Q\bigl(|Z_t-x_\star|>r\bigr)
\leq
C\,\mathbb P\bigl(|Z_t-x_\star|>r\bigr)
\longrightarrow0.
\]
Thus, the new initial distribution also concentrates at the
target.

In particular, if $\rho_T$ is continuous and strictly positive
near a prescribed initial configuration $x_0$, we may choose
$q_0$ uniformly on an arbitrarily small ball around $x_0$.
This gives an informal justification for planning from a
localized initial distribution near that configuration. 

\section{Integrating Path Planning Algorithms with Flow Matching}

One can combine classical path planning problems from robotics with flow matching. In a typical path planning problem, one is given a closed Euclidean domain with obstacles and the goal is to find a control law that navigates the system to the goal configuration while avoiding obstacles. One way to achieve this synergy is to leverage the principle of noising and denoising.  One could use a noising process that remains confined to allowable regions of state-space, but explores it sufficiently. And then denoising based feedback control could be used to construct a policy that achieves navigation through the domain starting from any initial condition.

\subsection{RRT* Noising-Denoising}

Standard noising processes, such as those constructed by using Brownian motion or white noise (see, for example Chapter \ref{chap:denoising-diffusion}), may explore the state space slowly, and a large number of samples may be required to fill the space. This issue is especially pronounced in nonconvex domains.

It is therefore natural to consider if there are effective noising mechanisms beyond Brownian motion. One option is to use powerful sampling-based path-planning algorithms from robotics \cite{lavalle2006planning}, such as the rapidly exploring random tree family \cite{lavalle2001randomized,karaman2011sampling}.

RRT* is usually formulated to construct trajectories between two points in state space \cite{karaman2011sampling}. Strictly speaking, one could solve a point-to-goal control problem for each sampled initial state and then apply flow matching directly. However, this is wasteful because each run of RRT* requires constructing a new graph from the initial state toward the goal configuration. A more efficient RRT* based noising and denoising algorithm is presented in the following algorithm.

\begin{tcolorbox}[algorithmbox,title=\textbf{Algorithm: RRT*-Based Flow Matching for Closed-Loop Synthesis}]
	\textbf{Setup.}
	Consider the single-integrator system
	\[
	\dot{x}=u,
	\qquad x\in\Omega_{\mathrm{free}}\subset\R^2,
	\]
	where $\Omega_{\mathrm{free}}$ denotes the obstacle-free region.
	
	\begin{enumerate}[label=\arabic*.]
		\item Construct an RRT* tree rooted at the goal configuration
		$x_\star$, using the  tree-construction algorithm.
		
		\item Sample tree nodes $q^{(j)}$ and reverse their root-to-node
		paths. Parameterize each path over
		$t\in[0,1]$, so that
		\[
		X^{(j)}_0=q^{(j)},
		\qquad
		X^{(j)}_1=x_\star.
		\]
		Store state--control samples
		$\{X_{t_k}^{(j)},U_{t_k}^{(j)}\}_{j=1}^{N}$,
		estimating $U^{(j)}_t=\dot X^{(j)}_t$ by finite differences.
		
		\item Train a network
		\[
		u_\theta:\R^2\times[0,1]\to\R^2
		\]
		by minimizing the empirical $L^2$ loss corresponding to
		\[
		\mathcal L(\theta)
		=
		\int_0^1
		\E\!\left[
		\left\lVert U_t-u_\theta(X_t,t)\right\rVert^2
		\right]\dd t.
		\]
		
		\item Apply the learned feedback from a sampled initial
		configuration $x(0)$:
		\[
		\dot{x}(t)=u_\theta\bigl(x(t),t\bigr).
		\]
		The controller is trained to steer the state toward $x_\star$
		as $t\to1$.
	\end{enumerate}
\end{tcolorbox}

We can zoom into the RRT* tree construction algorithm of the first steps as follows.

\begin{tcolorbox}[algorithmbox,title=\textbf{Algorithm: RRT* Tree Construction}]
\textbf{Setup.}
Initialize a tree with root $x_\star$ and cost $c(x_\star)=0$.
The cost $c(q)$ is the length of the tree path from $x_\star$
to node $q$. Choose a maximum extension length $\eta>0$
and a neighborhood radius $r>0$.

Repeat the following steps for a prescribed number of iterations.
\begin{enumerate}[label=\alph*)]
  \item \textbf{Sample.}
  Draw a random point $q_{\mathrm{rand}}$ in free space.

  \item \textbf{Extend.}
  Find the nearest tree node $q_{\mathrm{near}}$ and move toward
  the sample by at most $\eta$:
  \[
  q_{\mathrm{new}}
  =
  q_{\mathrm{near}}
  +
  \min\!\left\{
  1,\frac{\eta}{\lVert q_{\mathrm{rand}}-q_{\mathrm{near}}\rVert}
  \right\}
  (q_{\mathrm{rand}}-q_{\mathrm{near}}).
  \]
  Discard the new point if the segment
  from $q_{\mathrm{near}}$ to $q_{\mathrm{new}}$ intersects an obstacle.

  \item \textbf{Choose a parent.}
  Among $q_{\mathrm{near}}$ and the tree nodes within distance $r$
  of $q_{\mathrm{new}}$, choose a node $q_{\mathrm{parent}}$
  minimizing
  \[
  c(q)+\lVert q_{\mathrm{new}}-q\rVert,
  \]
  subject to the connecting segment being collision-free.
  Add $q_{\mathrm{new}}$ to the tree with this parent and set
  \[
  c(q_{\mathrm{new}})
  =
  c(q_{\mathrm{parent}})
  +
  \lVert q_{\mathrm{new}}-q_{\mathrm{parent}}\rVert.
  \]

  \item \textbf{Rewire.}
  For each nearby node $q$, check whether connecting it through
  $q_{\mathrm{new}}$ gives a shorter collision-free path:
  \[
  c(q_{\mathrm{new}})
  +
  \lVert q-q_{\mathrm{new}}\rVert
  <
  c(q).
  \]
  If so, change its parent to $q_{\mathrm{new}}$ and update
  the costs of $q$ and its descendants.
\end{enumerate}

\textbf{Output.}
A tree of collision-free paths rooted at $x_\star$.
\end{tcolorbox}

In general, designing feedback-stabilizing control laws that avoid obstacles is difficult. The time-reversal-based controller may encounter undesirable equilibria or become trapped near obstacle boundaries. Two heuristics can help mitigate these deadlocks during the reverse process:

\begin{enumerate}[label=\roman*)]
  \item At the boundary, project velocities onto directions that keep the sample inside the free space.

  \item Add a small amount of Gaussian noise to particles near the boundary to help them escape deadlocked configurations.
\end{enumerate}

The motivation for combining noise with projection-based constraint handling is related to projected diffusion processes for path planning. Stochastic perturbations promote exploration, while projection or reflection prevents trajectories from leaving the feasible domain.

\begin{figure}[htbp]
	\centering
	\includegraphics[width=\textwidth]
	{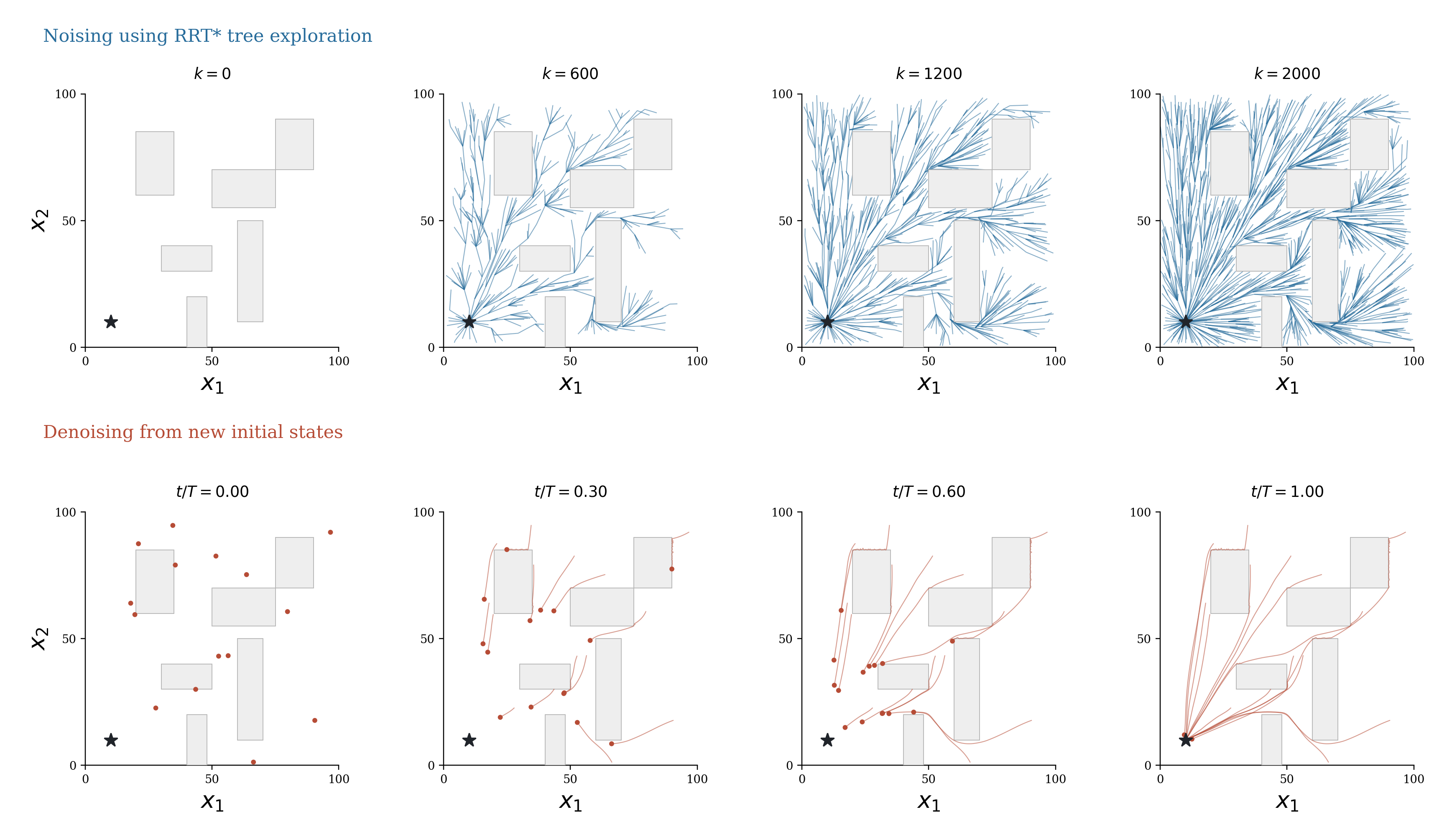}
	\caption{RRT*-based noising and learned denoising.
		Top: tree growth at iterations $0$, $600$, $1200$, and $2000$.
		Bottom: learned return trajectories from fresh initial states
		sampled uniformly in free space, shown at $t/T=0$, $0.3$,
		$0.6$, and $1$.
		Curves show accumulated paths, dots mark current states,
		and the star marks the goal.}
	\label{fig:rrt-noising-denoising}
\end{figure}

The results of an experiment can be seen in 
Figure~\ref{fig:rrt-noising-denoising}. It shows the growth of the
goal-rooted RRT* tree and the learned closed-loop trajectories
from independently sampled free-space initial states.

\subsection{Dijkstra-Based Flow Matching}

The previous section hints at the idea that, if we have a discrete collection of
trajectories, one can turn it into a continuous feedback law that approximately
reproduces the discrete behavior. This viewpoint helps isolate the phenomenon for
which flow matching is useful and keeps us from overselling what it can accomplish
for planning problems. In this setting, flow matching should be treated as a
complementary mechanism for path-planning, rather
than as a replacement.

This idea can be more concretely formalized. Let
\[
V=\{q_1,\ldots,q_N\}\subset \Omega
\]
be a finite discretization of the state space, and suppose that we are given a
discrete feedback law $
\pi:V\rightarrow V.$ Starting from $q_0\in V$, the discrete closed-loop dynamics are
\begin{equation}
	q_{k+1}=\pi(q_k).
\end{equation}
Thus, each initial condition generates a discrete trajectory
\[
q_0\rightarrow q_1\rightarrow q_2\rightarrow\cdots.
\]

The map $\pi$ tells us where to move next, but it does not tell us how to move
between two consecutive states. To obtain trajectories for the continuous-time
system
\begin{equation}
	\dot{x}=f(x,u),
\end{equation}
we therefore assume that a local interpolator is available:
\begin{equation}
	S(p,q)=\big(x(t),u(t)\big),
\end{equation}
where $
S:V\times V\longrightarrow \mathcal{X}\times\mathcal{U}$ and $\mathcal{X}$ denotes a space of state trajectories and $\mathcal{U}$ a space
of admissible control trajectories. The interpolation satisfies
\[
x(0)=p,\qquad x(1)=q,
\]
and, 
\begin{equation}
	\dot{x}(t)=f\big(x(t),u(t)\big).
\end{equation}

Thus, for every transition prescribed by the discrete policy,
\[
q_k\longrightarrow q_{k+1}=\pi(q_k),
\]
the interpolator provides a continuous state-control pair
\[
S\big(q_k,\pi(q_k)\big)
=
\big(x_k(t),u_k(t)\big).
\]
Concatenating these local interpolations gives a continuous realization of the
discrete trajectory.

For example, for the single-integrator system
\begin{equation}
	\label{eq:singleinte}
	\dot{x}=u,
\end{equation}
we can use the linear interpolation
\begin{equation}
	x(t)=(1-t)p+tq,\qquad t\in[0,1],
\end{equation}
with corresponding control
\begin{equation}
	u(t)=q-p.
\end{equation}
For a more general controlled system, $S(p,q)$ may use alternative approaches such as local steering strategies for nonholonomic systems \cite{jean2014control}.

Applying the discrete policy from many initial states therefore gives a family of
continuous state-control trajectories. At this point we are back in the
setting of Section~3.1. If $
(X_t,U_t)$
denotes a randomly sampled state-control pair from this family, then flow matching
constructs the feedback law
\begin{equation}
	u(t,x)=\mathbb{E}[U_t\mid X_t=x],
\end{equation}
which may be learned by minimizing
\begin{equation}
	\mathcal{L}(\theta)
	=
	\int_0^T
	\mathbb{E}
	\left[
	\|U_t-u_\theta(t,X_t)\|^2
	\right]dt.
\end{equation}

The role of flow matching is therefore to combine the discrete policy and the interpolator together generate the
open-loop trajectories,
\[
\pi
\quad+\quad
S
\quad\longrightarrow\quad
\{(X_t,U_t)\},
\]
and flow matching converts this collection into a continuous feedback law $
u_\theta(t,x).$

We now consider a simple example in which the discrete feedback law $\pi$ can be
computed exactly.  Once again, let
\[
V=\{q_1,\ldots,q_N\}\subset\Omega_{\rm free}
\]
be a finite grid of states, and let $x^\star\in V$ denote the goal
configuration. For every $q\in V$, let $\mathcal{N}(q)$ denote the set of neighboring
grid points that can be connected to $q$ without intersecting an obstacle.

Associated with every admissible edge $(p,q)$ is a nonnegative cost $
c(p,q)$
For a Euclidean grid, we take $
c(p,q)=\|p-q\|.$

The objective is to find, for every $q\in V$, the minimum cost of reaching the
fixed goal $x^\star$. Define the discrete value function
\begin{equation}
	J(q)
	:=
	\min_{\substack{
			q_0=q,\;q_K=x^\star\\
			q_{k+1}\in\mathcal{N}(q_k)
	}}
	\sum_{k=0}^{K-1}
	c(q_k,q_{k+1}).
\end{equation}
The value function satisfies the discrete Bellman equation
\begin{equation}
	J(q)
	=
	\min_{p\in\mathcal{N}(q)}
	\left\{
	c(q,p)+J(p)
	\right\},
	\qquad
	J(x^\star)=0.
\end{equation}
Consequently, an optimal discrete feedback law is
\begin{equation}
	\pi(q)
	=
	\arg\min_{p\in\mathcal{N}(q)}
	\left\{
	c(q,p)+J(p)
	\right\}.
\end{equation}

Thus, once $J$ has been computed, the closed-loop discrete dynamics
\[
q_{k+1}=\pi(q_k)
\]
follow a shortest path from any grid state toward the goal.

An important feature of this construction is that we do not solve a separate
point-to-point problem for every possible initial state. Instead, we fix the terminal
condition $x^\star$ and propagate information outward from the goal through the
state space. {\it Dijkstra's algorithm} performs precisely this propagation or {\it noising}. Starting from $
J(x^\star)=0.$
it successively assigns optimal costs to states farther and farther from the goal.

The direction of this computation is opposite to the direction in which the resulting
feedback law is applied. The value information propagates outward from $x^\star$,
while the policy $\pi$ points inward toward $x^\star$. In this sense, a single outward
exploration produces a family of paths that can subsequently be traversed toward the
goal from many different initial conditions. The algorithm for exploration is shown in the following description.

\medskip

\noindent
\begin{tcolorbox}[algorithmbox,title=\textbf{Algorithm: Dijkstra-Based Discrete Feedback}]
	\textbf{Setup.}
	Let
	\[
	V\subset\Omega_{\mathrm{free}}
	\]
	be a finite grid, let $x_\star\in V$ be the goal configuration, and let
	$\mathcal N(q)$ denote the collision-free neighbors of $q$. Associate to
	each admissible edge $(p,q)$ a nonnegative cost $c(p,q)$.
	
	\begin{enumerate}[label=\arabic*.]
		\item Initialize
		\[
		J(x_\star)=0,
		\qquad
		J(q)=+\infty
		\quad\text{for }q\neq x_\star.
		\]
		Set
		\[
		\mathrm{OPEN}=\{x_\star\},
		\qquad
		\mathrm{CLOSED}=\varnothing.
		\]
		
		\item Choose
		\[
		q
		=
		\arg\min_{p\in\mathrm{OPEN}} J(p).
		\]
		
		\item Remove $q$ from $\mathrm{OPEN}$ and add it to
		$\mathrm{CLOSED}$.
		
		\item For every collision-free neighbor
		$p\in\mathcal N(q)$, compute the candidate cost
		\[
		\widetilde J(p)
		=
		c(p,q)+J(q).
		\]
		
		\item If
		\[
		\widetilde J(p)<J(p),
		\]
		update
		\[
		J(p)\leftarrow\widetilde J(p),
		\qquad
		\pi(p)\leftarrow q,
		\]
		and add $p$ to $\mathrm{OPEN}$ if it is not already present.
		
		\item Repeat Steps 2--5 until $\mathrm{OPEN}$ is empty, or until the
		desired portion of the state space has been explored.
	\end{enumerate}
\end{tcolorbox}

The overall construction therefore has two directions. First, information is propagated
outward from the fixed goal to construct a family of feasible paths throughout the
state space to get the feedback law $\pi $. These paths are then traversed in the opposite direction to produce
state--control trajectories terminating at $x^\star$, from which the continuous feedback
law is learned.

\begin{figure}[t]
	\centering
	\includegraphics[width=\textwidth]{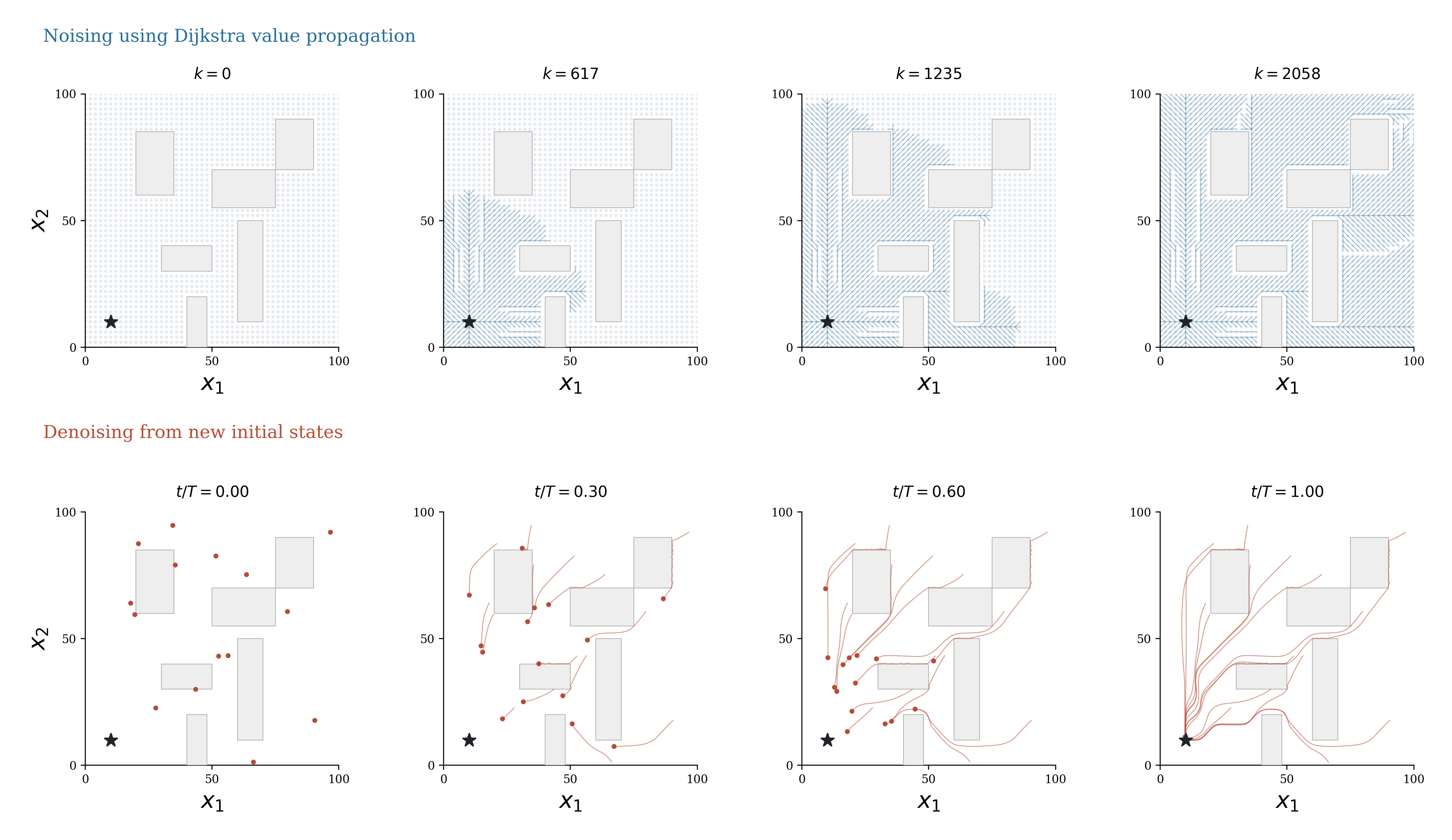}
	\caption{
		Dijkstra-based flow matching.
		Top: outward propagation of the shortest-path tree from the goal configuration at successive stages of the Dijkstra algorithm.
		Bottom: trajectories generated by the learned continuous feedback controller from fresh initial conditions sampled independently in the free space.
		Curves show accumulated trajectories, dots indicate the current states, and the star marks the goal.
	}
	\label{fig:dijkstra_flow_matching}
\end{figure}

Figure~\ref{fig:dijkstra_flow_matching} illustrates the two stages of the construction using Dijkstra's algorithm. The top row shows the outward propagation of the Dijkstra shortest-path tree from the fixed goal, while the bottom row shows trajectories generated by the learned continuous feedback law from fresh initial conditions sampled in free space, assuming single integrator dynamics \eqref{eq:singleinte}.

\section{Optimal Control based Flow Matching via PMP Time-Reversal}
\label{sec:pmp-flow-matching}
\label{sec:pmp-noising-denoising}

The noising mechanisms initially considered choosing some form of randomization strategy. The Dijkstra based approach made us realize that noising could be carried out in a optimality-informed way, by considering a discrete state planning problem. One can execute this idea more directly at the level of the continuous-time continuous state system. 

In this section, we use the Pontryagin's Maximum Principle from optimal control to construct alternative optimality-informed noising schemes. In {\it PMP-based noising}, instead of sampling arbitrary feasible controls, it samples trajectories that also satisfy first-order optimality conditions.

The first section on flow matching rested on the idea that, if one can solve the open-loop control problem for each pair of initial and final conditions, then one can use these solutions to synthesize a closed-loop control. This can potentially be combined with optimal control to solve the optimal feedback synthesis problem.

We recall that there are two classical routes to optimal controls. The first uses the Pontryagin Maximum Principle (PMP); the second uses the Hamilton--Jacobi--Bellman (HJB) equation. The HJB equation directly yields a feedback control, so it does not naturally fit the flow-matching framework. The PMP, by contrast, gives solutions to the open-loop optimal control problem. In principle, one can therefore use the PMP to find open-loop optimal controls and then use flow matching to synthesize a feedback control from multiple samples. The obstacle is that one must first solve a two-point boundary value problem.

\paragraph{The PMP system}

Let us review the boundary value problem that stands in our way. Consider the optimal control problem
\begin{equation}
    \inf_{x,u}\int_0^T L\bigl(x(t),u(t)\bigr)\,\dd t
\end{equation}
subject to
\begin{equation}
    \dot{x}(t)=f\bigl(x(t),u(t)\bigr),
    \qquad x(0)=x_0,
    \qquad x(T)=y.
\end{equation}
where $L : \mathbb{R}^d \times \mathbb{R}^{m} \rightarrow \mathbb{R}$ is the {\it running cost.}

Define the Hamiltonian
\begin{equation}
    H(x,p,a)
    =
    \langle p,f(x,a)\rangle
    +
    L(x,a).
\end{equation}
With the minimization convention used here, the Pontryagin Maximum Principle states that any optimal state--control pair must be accompanied by an adjoint trajectory $p(t)$ satisfying
\begin{align}
    \dot{\omega}(t)
    &=f\bigl(\omega(t),\alpha(\omega(t),p(t))\bigr),\\
    \dot{p}(t)
    &=-\nabla_x f\bigl(\omega(t),\alpha(\omega(t),p(t))\bigr)^{\mathsf T}p(t)
      -\nabla_x L\bigl(\omega(t),\alpha(\omega(t),p(t))\bigr),
\end{align}
with boundary conditions
\begin{equation}
    \omega(0)=x_0,
    \qquad
    \omega(T)=y,
\end{equation}
and minimizing control
\begin{equation}
    \alpha(\omega,p)
    =
    \operatorname*{arg\,min}_{a\in\R^m}
    H(\omega,p,a)
    =
    \operatorname*{arg\,min}_{a\in\R^m}
    \left[
        \langle p,f(\omega,a)\rangle+L(\omega,a)
    \right].
\end{equation}
Here the $x$-derivatives in the adjoint equation are taken with the control argument held fixed at the minimizing value $\alpha(\omega(t),p(t))$.

This system is hard to solve because boundary conditions appear at both endpoints: hence the name \emph{two-point boundary value problem}. A standard numerical approach is direct multiple shooting, in which the time interval is divided into several segments and continuity between segmentwise trajectories is imposed as a nonlinear matching constraint.
Shooting methods suffer some drawbacks. One of them is that they are sensitive to the initial guess for the adjoint state $p(0)$, so poor initialization can lead to divergence or convergence to spurious solutions. 

\paragraph{Swapping the boundary-value problem for an initial-value problem.}

An alternative is to replace the initial condition on the state by a terminal condition on the adjoint state:
\begin{align}
    \dot{\omega}(t)
    &=f\bigl(\omega(t),\alpha(\omega(t),p(t))\bigr),\\
    \dot{p}(t)
    &=-\nabla_x f\bigl(\omega(t),\alpha(\omega(t),p(t))\bigr)^{\mathsf T}p(t)
      -\nabla_x L\bigl(\omega(t),\alpha(\omega(t),p(t))\bigr),
\end{align}
with terminal data
\begin{equation}
    p(T)=p_T,
    \qquad
    \omega(T)=y,
\end{equation}
and
\begin{equation}
    \alpha(\omega,p)
    =
    \operatorname*{arg\,min}_{a\in\R^m}
    \left[
        \langle p,f(\omega,a)\rangle+L(\omega,a)
    \right].
\end{equation}

This system is easily solved backward in time using standard initial value problem solvers. The cost is that information about the initial condition no longer enters the optimality system explicitly. To recover the missing initial-condition information, we invoke the noising--denoising principle developed in the previous section. Treating the backward-in-time solution of the PMP system as a noising process, we observe the following. If $p_T$ is sampled from a distribution with full support on $\R^n$, then integrating the PMP system backward from the goal state can explore the backward-reachable set. If, for every initial condition, the optimal trajectory--control pair satisfies the PMP, then sufficiently rich sampling of $p_T$ should induce correspondingly rich sampling of initial states, provided the map $p_T \mapsto x(t)$ is well behaved. This yields an analogue of the support property from the previous lecture. The only difference is that we now sample terminal adjoint states rather than entire control trajectories.

\textbf{PMP-Based Support Property:}
Suppose every optimal control is represented by the PMP system and the map from terminal adjoint states $p_T$ to the resulting backward-integrated initial states $X_0$ is continuous and sufficiently surjective. If $p_T$ is sampled from a distribution with full support on $\R^n$, then the induced law of $X_0$ has full support on the corresponding backward-reachable set. In particular, if that set is all of $\R^n$, then
\[
    \operatorname{supp}\bigl(\mathcal L(X_0)\bigr)
    =
    \R^n.
\]

Therefore, by sampling a range of terminal adjoint states and integrating backward, we obtain an optimality-informed noising process. Exactly as before, flow matching can then be used to learn a controller that denoises this process as described in the following algorithm.

\begin{tcolorbox}[algorithmbox,title=\textbf{Algorithm: PMP-Based Flow Matching}]
\begin{enumerate}[label=\arabic*.]
    \item Fix the goal state $x_\star$.

    \item Sample
    \[
        p_T\sim\nu_p,
    \]
    where $\nu_p$ has full support on $\R^n$.

    \item Integrate the PMP system backward in time from
    \[
        \omega(T)=x_\star,
        \qquad
        p(T)=p_T.
    \]

    \item Collect the control data
    \[
        U_t=\alpha\bigl(X_t,\rho_t\bigr)
    \]
    across many samples of $p_T$.

    \item Learn a flow-matched controller $u^\theta(t,X_t)$ by minimizing
    \[
        \mathcal{L}(\theta)
        =
        \int_0^T
        \E\!\left[
            \left\lVert U_t-u^\theta(t,X_t)\right\rVert^2
        \right]
        \dd t.
    \]
\end{enumerate}
\end{tcolorbox}

The classical and denoising-based approaches to optimal control can be summarized as follows.

\begin{table}[htbp]
    \centering
    \begin{tabular}{@{}p{0.28\textwidth}p{0.64\textwidth}@{}}
        \toprule
        \textbf{Approach} & \textbf{Procedure}\\
        \midrule
        Classical optimal control
        & Fix initial and final conditions, then solve the PMP-based boundary value problem or the HJB equation.\\[0.4em]
        Denoising-based optimal control
        & Fix the final condition, propagate states through the backward-reachable set using the PMP system, and denoise to synthesize a feedback controller.\\
        \bottomrule
    \end{tabular}
    \caption{Classical and denoising-based routes to feedback synthesis.}
    \label{tab:classical-vs-pmp-denoising}
\end{table}

\paragraph{Unicycle example.}

We test the method on minimum-energy control of the unicycle model:
\begin{equation}
    \dot{x}_1=u_1\cos(x_3),
    \qquad
    \dot{x}_2=u_1\sin(x_3),
    \qquad
    \dot{x}_3=u_2.
\end{equation}

\begin{figure}[htbp]
	\centering
	\includegraphics[width=\textwidth]
	{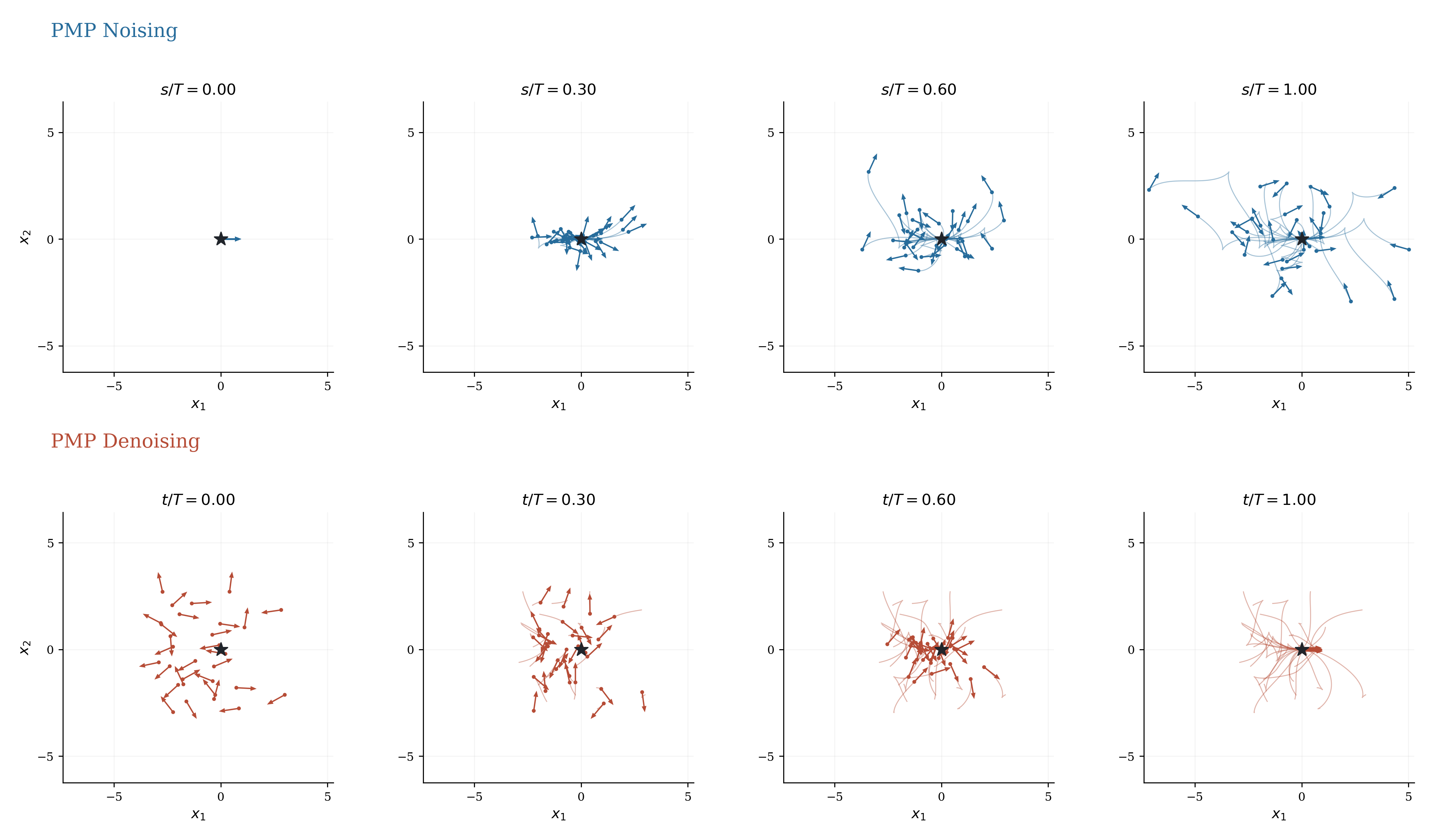}
	\caption{PMP-based noising and learned denoising for the unicycle.
		Top: trajectories spreading from the goal under sampled adjoint
		conditions. Bottom: learned feedback trajectories from uniformly
		sampled initial configurations.
		Curves show paths, while arrows indicate the heading $x_3$,
		and the star marks the target position.}
	\label{fig:unicycle-pmp-noising-denoising}
\end{figure}

Figure~\ref{fig:unicycle-pmp-noising-denoising} illustrates PMP-based
noising and learned denoising for the unicycle. Starting at the goal
configuration, sampled adjoint conditions generate the noising
trajectories. The learned feedback then steers independently sampled
initial configurations toward the goal. Both phases are shown in the
$(x_1,x_2)$ plane, with arrows indicating the heading $x_3$.

\paragraph{What optimality is retained?}

The PMP is a necessary condition for optimality but is not, in general, sufficient. Accordingly, there is no guarantee that the learned controller is globally optimal. What we recover is a feedback controller informed by the optimality structure of the problem. For linear systems with convex running costs and convex admissible control sets, the PMP conditions are sufficient under the standard convexity assumptions. In such cases, optimality of the learned controller can be guaranteed when the flow-matching step exactly recovers the corresponding optimal feedback field.

\paragraph{Further Reading}

The flow matching philosophy as adapted to control problems can be found in \cite{elamvazhuthi2025flow}. Flow matching for stochastic linear systems can be found in \cite{mei2024flow}.

\chapter{How Not to Sample from the Reachable Set}
\label{chap:not-sample-reachable}

As we saw in Chapter~\ref{chap:flow-matching}, an important element of the noising--denoising flow-matching approach to feedback control is the noising step. The noising step propagates randomness through the control channels so that the state law becomes fully supported on the relevant backward-reachable set. There are, however, some naive ways to randomize controls that work naturally in discrete time but fail in the continuous-time limit. In this chapter, we will distract ourselves with these naive strategies. The reason for this indulgence is that this distinction matters beyond flow matching, including areas such as system identification, entropy-based reinforcement learning, and path-integral control.

\section{Discrete-Time Intuition, Continuous-Time Failure}

Consider the discrete-time control system
\begin{equation}
	x_{n+1}=f(x_n,u_n).
\end{equation}
For a fixed horizon $N$, define the reachable set from $x_0$ by
\begin{equation}
	\mathcal R_N(x_0;U)
	:=
	\left\{
	x_N:\
	\begin{array}{l}
		\text{there exists }(u_0,\ldots,u_{N-1})\in U^N\\
		\text{such that }(x_n)_{n=0}^N\text{ is a system trajectory}
	\end{array}
	\right\}.
\end{equation}

Suppose $\pi(\dd u\mid x)$ is a stochastic control law with support equal to $U$, and choose
\begin{equation}
	u_n\sim \pi(\dd u\mid x)
\end{equation}
independently at each time step. Under standard assumptions, the map from the finite control sequence $(u_0,\ldots,u_{N-1})$ to $x_N$ is continuous. Consequently,
\begin{equation}
	\operatorname{supp}\bigl(\mathcal L(x_N)\bigr)
	=
	\overline{\mathcal R_N(x_0;U)}.
\end{equation}
where $(\mathcal L(\cdot)$denotes the law or distribution of the variable $x_N$.
If the reachable set is closed, then the support is exactly $\mathcal R_N(x_0;U)$. Thus, in discrete time, independently sampling a full-support control at every step is a natural way to sample the reachable set.

Now consider the continuous-time system
\begin{equation}
	\dot x=f(x,u).
\end{equation}
A tempting procedure is to discretize it using forward Euler scheme,
\begin{equation}
	x_{k+1}
	=
	x_k+\Delta t\,f(x_k,u_k),
\end{equation}
and again draw the $u_k$ independently from a fixed distribution $\pi(\dd u\mid x_k)$.

The important issue is what happens as $\Delta t\to0$. We know from results in {\it stochastic approximation theory} \cite{borkar2008stochastic} that with fixed-amplitude control samples and under the usual stochastic-approximation assumptions, the random fluctuations average out. The limiting dynamics are
\begin{equation}
	\boxed{
		\dot x
		=
		\int_U f(x,u)\,\pi(\dd u\mid x).
	}
	\label{eq:relaxed-control-average}
\end{equation}
The limit is therefore deterministic! Instead of obtaining many trajectories that explore the reachable set, repeated simulations converge toward the trajectory generated by the averaged vector field.

The conditional probability $\pi(\dd u\mid x)$ appearing in \eqref{eq:relaxed-control-average} is commonly interpreted as a \emph{relaxed control}. Although it resembles a discrete-time stochastic policy, its continuous-time effect is through the averaged velocity field.

The loss of exploration is especially transparent for
\begin{equation}
	\dot x
	=
	f_0(x)
	+
	\sum_{i=1}^m u_i g_i(x).
\end{equation}
The relaxed dynamics are
\begin{equation}
	\dot x
	=
	f_0(x)
	+
	\sum_{i=1}^m
	\E_\pi[u_i\mid x]g_i(x).
\end{equation}
Only the conditional means of the controls enter the dynamics. Two policies with very different variances or entropies can therefore generate exactly the same state trajectory.

For example, if the control distribution is symmetric and  $ \E_\pi[u\mid x]=0,$ then a driftless system does not move at all in the relaxed limit, while a system with drift simply follows $f_0$.

A simple illustration is the controlled Van der Pol oscillator
\begin{align}
	\dot x &= v,\\
	\dot v &= (1-x^2)v-x+u,
\end{align}
with $u\in[-4,4]$ and, for example,
\begin{equation}
	x(0)=1.5,
	\qquad
	v(0)=0.
\end{equation}
If an independent bounded control is sampled at every Euler step, as the time step is reduced, the trajectories increasingly collapse onto the averaged dynamics.

\begin{figure}[htbp]
	\centering
	\includegraphics[width=\textwidth]
	{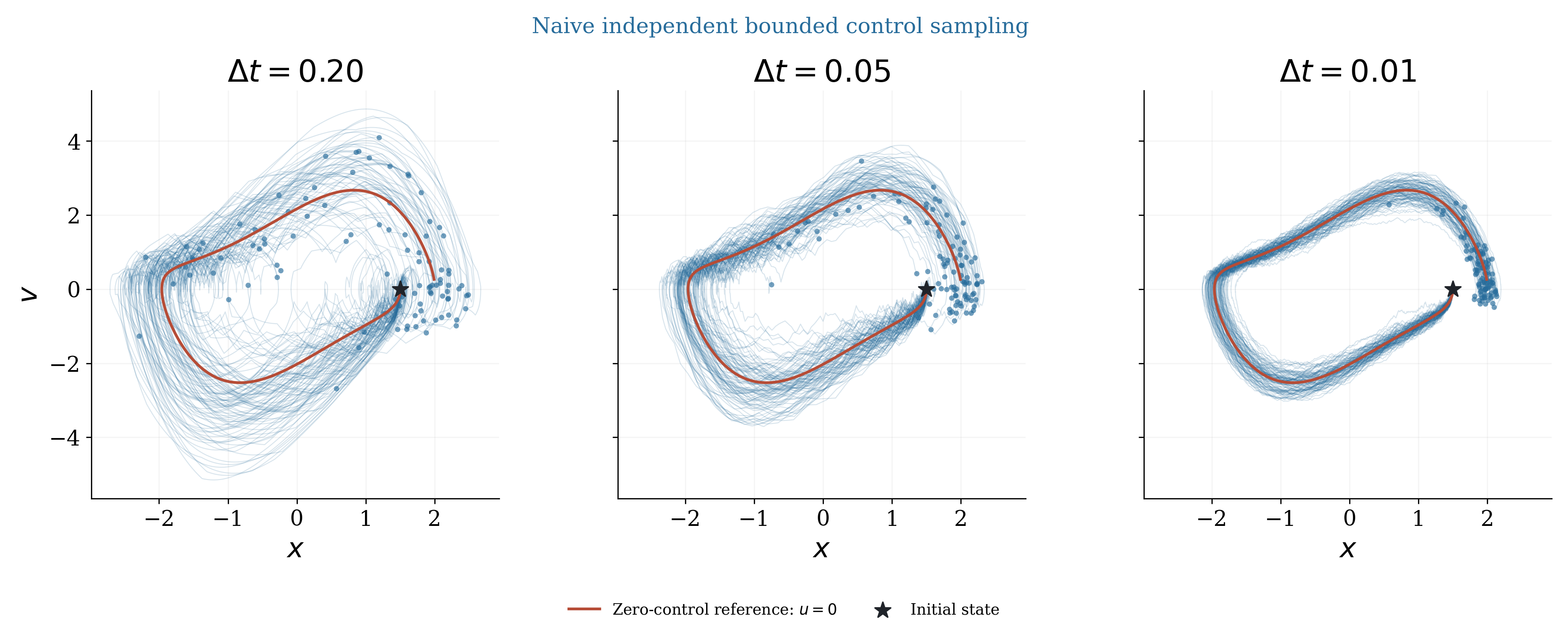}
	\caption{Independent bounded control sampling for the Van der Pol
		oscillator with $\Delta t=0.2$, $0.05$, and $0.01$.
		Blue curves show sampled trajectories, and blue dots mark their
		terminal states. The red curve shows the averaged dynamics
		with $u=0$; the star marks the common initial state.}
	\label{fig:naive-reachable-sampling}
\end{figure}

Figure~\ref{fig:naive-reachable-sampling} illustrates the loss of
exploration caused by increasingly frequent independent control
sampling. We consider the controlled Van der Pol oscillator
\[
\dot x=v,
\qquad
\dot v=(1-x^2)v-x+u,
\]
initialized at $(x(0),v(0))=(1.5,0)$. On each interval of length
$\Delta t$, the control is sampled independently and uniformly
from $[-4,4]$ and held constant. As $\Delta t$ decreases, the
trajectories concentrate around the averaged dynamics, which
corresponds to $u=0$. Thus, over a fixed observation horizon,
more frequent randomization can reduce state-space exploration.

\paragraph{Consequences for control problems.} 
~

\emph{System identification.}

Randomized controls are often used to excite a system when collecting
data for identifying a model for the system. Although the preceding discussion
concerns sampling multiple trajectories, the same loss of excitation
can arise when learning from a single trajectory, as is typical in many works on system identification. If a fixed-amplitude
control is independently resampled on increasingly short intervals,
then, under suitable assumptions, the state approaches the trajectory
of the averaged dynamics over a fixed observation horizon: 
\begin{equation}
	\dot x
	=
	\int_U f(x,u)\,\pi(\dd u\mid x).
	\label{eq:relaxed-control-average}
\end{equation}
Collecting
more samples does not by itself restore the lost state-space
exploration, and identification can become ineffective when the
limiting trajectory does not provide enough information to distinguish
the dynamics within the chosen model class.

\emph{Entropy-regularized reinforcement learning.}

Another instance, where our discussion is relevant is in reinforcement learning. In discrete time, increasing the entropy of a stochastic policy is potentially one way of generating a larger variety of state trajectories in order to encourage exploration in reinforcement learning \cite{haarnoja2018soft}.  The same interpretation does not automatically carry over to continuous-time relaxed controls. Since the continuous-time dynamics see only $
\int_U f(x,u)\,\pi(\dd u\mid x),$
large entropy in control set $U$ does not by itself imply large entropy or broad exploration in state space. One way to address this, is to encourage entropy over the state-space instead of the control set. Another approach to address these issues is to enforce entropy regularization over the space of control trajectories instead.

\section{Random Alternatives}

\paragraph{Propagation using White Noise}
One way to address the issue highlighted in the previous section is to scale the amplitude of the
random controls proportionally to $1/\sqrt{\Delta t}$. For example,
consider the control-affine system
\[
\dot X_t=f_0(X_t)+\sum_{i=1}^{m}u_i(t)g_i(X_t),
\]
and choose piecewise-constant controls
\[
u_i(t)=\frac{1}{\sqrt{\Delta t}}\,\xi_{i,k},
\qquad
t\in[k\Delta t,(k+1)\Delta t),
\]
where the $\xi_{i,k}$ are independent standard Gaussian random
variables. The integrated control increments then have variance
$\Delta t$, allowing their cumulative effect to persist
as $\Delta t\to0$.

Under suitable regularity  assumptions, solving
the controlled ODE with these inputs leads to the
\emph{Stratonovich stochastic differential equation}
\[
\dd X_t
=
f_0(X_t)\dd t
+
\sum_{i=1}^{m}g_i(X_t)\circ\dd W_t^i,
\qquad X_0=x_0,
\]
where the $W^i$ are independent standard Brownian motions.
Thus, the limiting perturbations can be interpreted formally as
white noise acting through the control channels.

The connection with reachability is made precise by the
Stroock--Varadhan support theorem \cite{stroock1972support}.
Under suitable regularity and nonexplosion assumptions,
the support of the trajectory law is the closure, in the
uniform topology on $[0,T]$, of the controlled trajectories
satisfying
\[
\dot x^h_t
=
f_0(x^h_t)
+
\sum_{i=1}^{m}g_i(x^h_t)\dot h_i(t),
\qquad x^h_0=x_0,
\]
where $h(0)=0$, $h$ is continuously differentiable, and
$\dot h\in L^2(0,T;\R^m)$.

The controls $u_i=\dot h_i$ range over all finite-energy
controls. Consequently,
\[
\operatorname{supp}\bigl(\mathcal L(X_T)\bigr)
=
\overline{\mathcal R_T^{L^2}(x_0)},
\]
where $\mathcal R_T^{L^2}(x_0)$ denotes the set of endpoints
reachable at time $T$ using such controls. Thus, every
neighborhood of a reachable endpoint has positive probability.

A severe drawback of this choice of controls is that these are not subject to a common bound. Hence, the white-noise is not compatible
with a fixed bounded control set.

\begin{figure}[htbp]
	\centering
	\includegraphics[width=\textwidth]
	{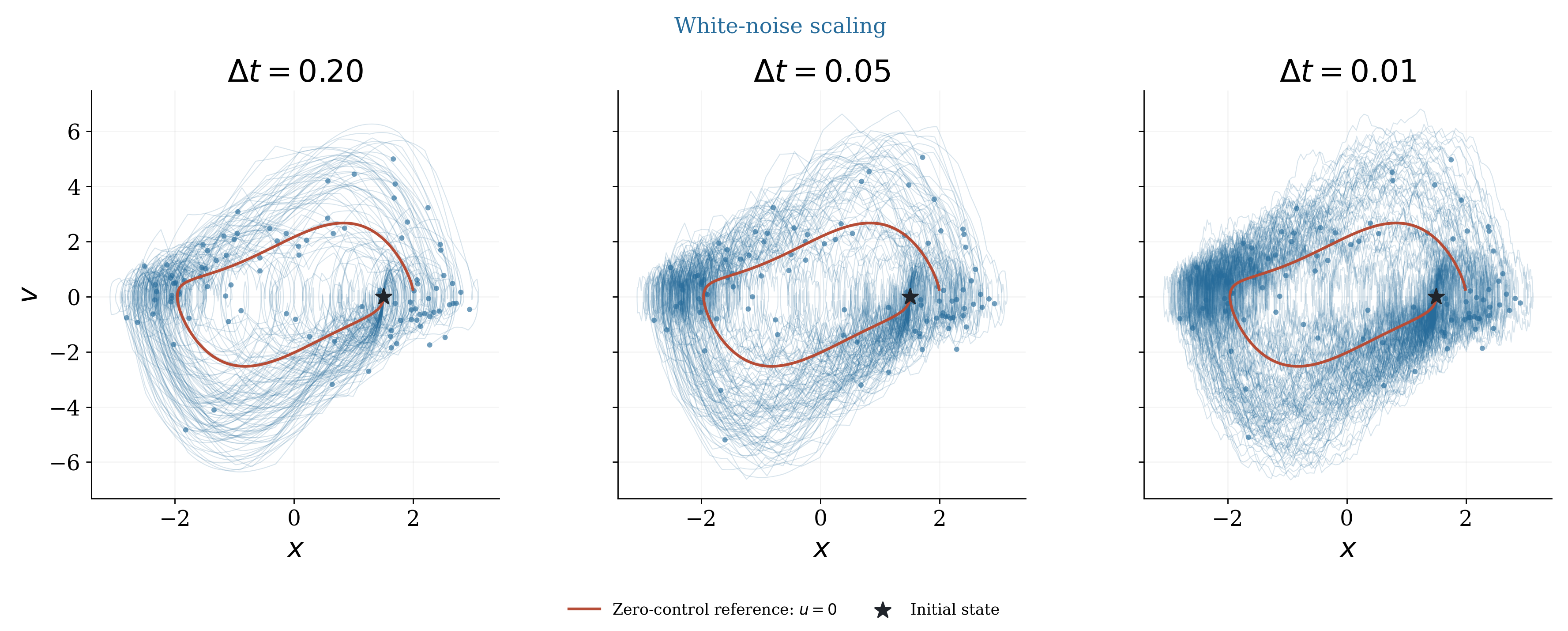}
	\caption{White-noise scaling for the controlled Van der Pol
		oscillator. On each interval of length $\Delta t$, the control
		is $u_k=\sigma\xi_k/\sqrt{\Delta t}$, where
		$\xi_k\sim\mathcal N(0,1)$ are independent and $\sigma=2$.
		Blue curves show sampled trajectories, and blue dots mark
		their terminal states. The red curve shows the zero-control
		dynamics; the star marks the initial state $(1.5,0)$.
		All panels use the same horizon $T=6$ and axis limits.}
	\label{fig:vdp-white-noise}
\end{figure}

 One can see in \ref{fig:vdp-white-noise} that white noise avoids the failure mode of sampling from the same unscaled control set.

\paragraph{Propagation using Brownian Motion} To deal with bounds on control, a more appropriate continuous-time construction is to place a probability measure directly on an admissible control-function space, as we did when performing the noising step using flow matching. For example, let
\begin{equation}
	\mathcal U
	=
	\left\{
	u\in C([0,T];\R^m):
	u(t)\in U\text{ for every }t
	\right\}.
\end{equation}
Suppose we choose a probability measure $\eta$ on $\mathcal U$ having full support. That is, every open neighborhood of every admissible control function has positive probability.

Similar to the choice of white noise, one can provide a support result. Suppose the solution map $
u\mapsto x_T^u$
is continuous in the topology used on $\mathcal U$. If $\eta$ has full support on $\mathcal U$, then
\begin{equation}
	\operatorname{supp}\bigl(\mathcal L(x_T^u)\bigr)
	=
	\overline{\mathcal R_T^{C}(x_0;U)},
	\qquad u\sim\eta.
\end{equation}
Thus a fully supported law on entire control paths induces a fully supported law on the reachable set up to closure.

For unconstrained controls, the {\it Wiener measure} considered in the flow matching chapter provides a useful example. This is the measure induced on the set of paths by Brownian motion. That is, instead of white noise, we are considering integral of white noise. It is known that Brownian motion $W_t$ has full support on
\begin{equation}
	C_0([0,T];\R^m)
	:=
	\{u\in C([0,T];\R^m):u(0)=0\}
\end{equation}
with respect to the uniform topology. Thus Brownian paths can be used as random \emph{control functions}. This should not be confused with white noise, which is formally the time derivative of Brownian motion.

For a closed convex bounded control set $U$, one may take a random variable $Y$ with full support on $U$, an independent Brownian motion $W$, and define
\begin{equation}
	u(t)
	=
	P_U\bigl(Y+\sigma W_t\bigr),
	\qquad \sigma>0,
\end{equation}
where $P_U$ is the Euclidean projection onto $U$. Because pointwise projection onto a closed convex set is continuous and onto, the resulting law has full support on the corresponding space of continuous $U$-valued controls.

\begin{figure}[htbp]
	\centering
	\includegraphics[width=\textwidth]
	{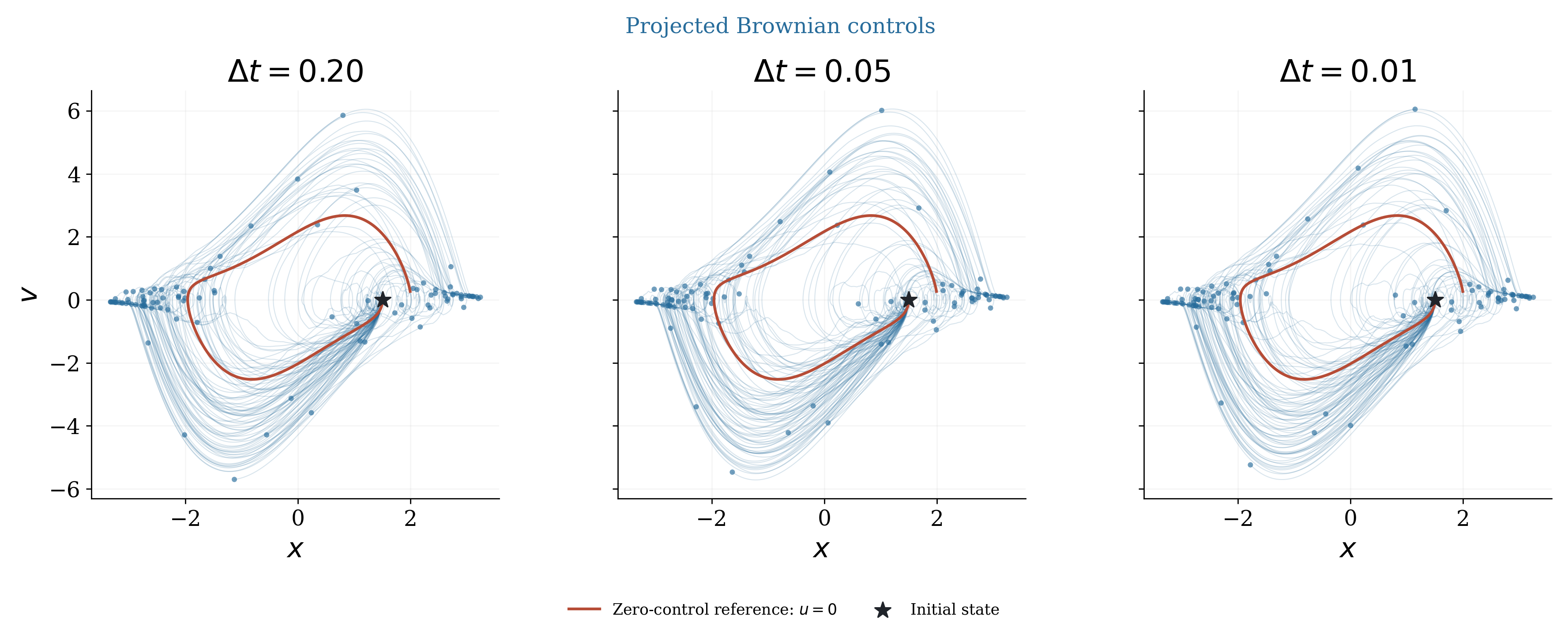}
	\caption{Projected Brownian controls for the Van der Pol
		oscillator. The control is
		$u(t)=P_{[-4,4]}(2W_t)$, approximated by
		piecewise-constant values sampled at intervals of length
		$\Delta t$. Blue curves show sampled state trajectories,
		and blue dots mark their terminal states. The red curve
		shows the zero-control dynamics; the star marks the initial
		state $(1.5,0)$. All panels use the same Brownian realizations,
		horizon $T=6$, and axis limits.}
	\label{fig:vdp-projected-brownian}
\end{figure}

Figure~\ref{fig:vdp-projected-brownian} illustrates the use of
bounded controls obtained by projecting Brownian motion onto
$[-4,4]$. Unlike independent resampling, in this construction, as the discretization is refined, the random excitation does not disappear. Interestingly, the trajectories are much smoother than white noise based propagation.

\paragraph{Piecewise Deterministic Markov Processes.}
Instead of Brownian motion, one can use random switching to construct
a probability measure with full support on the set of admissible
controls taking values in a compact and convex set $U\subset\R^m$. This construction is inspired by the literature on piecewise deterministic Markov processes and their support properties \cite{benaim2015qualitative}.

The idea is to restrict sampled controls to piecewise constant controls. But unlike the failure mode of naive discrete time sampling, we allow the number of switchings and switching times of the controls to be random as well. Toward this end, define a suitable probability law on the number of
switches, their locations, and the control values between switches.
Fix the horizon $T>0$ and a switching rate $\lambda>0$. First sample
\[
N\sim\operatorname{Poisson}(\lambda T),
\]
where $\operatorname{Poisson}(\lambda T)$ denotes the Poisson
distribution with mean $\lambda T$, so that
\[
\mathbb{P}(N=n)
=
e^{-\lambda T}\frac{(\lambda T)^n}{n!},
\qquad n=0,1,2,\ldots.
\]
Conditional on $N=n$, draw $n$ independent uniform samples from
$(0,T)$ and sort them to obtain
\[
0=t_0<t_1<\cdots<t_n<t_{n+1}=T.
\]
Independently sample $a_0,\ldots,a_n$ from a probability measure
$\mu$ with $\operatorname{supp}\mu=U$, and define
\[
u(t)=a_k,
\qquad t\in[t_k,t_{k+1}),\quad k=0,\ldots,n.
\]
The terminal time is therefore fixed, while the number and locations
of the switches are random. This construction induces a probability measure $\nu$ on control
trajectories. Every finite switching schedule can be approximated
with positive probability. Since finite piecewise-constant controls
are dense in $L^1(0,T;U)$, it follows that
\[
\nu\!\left(
\left\{u:\|u-v\|_{L^1(0,T)}<\varepsilon\right\}
\right)>0
\]
for every $v\in L^1(0,T;U)$ and every $\varepsilon>0$.
Although each sampled control has finitely many switches, allowing arbitrarily large values of $N$ is gives us
this full-support property. Here, $L^1(0,T;U)$ is the set of $U$-valued controls that have a finite integral over $(0,T)$.

Equivalently, the control value is refreshed from $\mu$ at the
arrival times of a Poisson process of rate $\lambda$. Between
arrivals, the state evolves according to
\[
\dot X_t=f_0(X_t)+\sum_{i=1}^{m}u_i(t)g_i(X_t).
\]
The joint state--control process is thus a piecewise deterministic
Markov process. Under assumptions ensuring existence through $T$
and continuity of the control-to-endpoint map $E_T(u)=x_T^u$,
the terminal law satisfies
\[
\operatorname{supp}\bigl((E_T)_\#\nu\bigr)
=
\overline{\mathcal R_T^{L^1}(x_0;U)},
\]
where $\mathcal R_T^{L^1}(x_0;U)$ is the reachable set under measurable
$U$-valued controls.

\begin{figure}[htbp]
	\centering
	\includegraphics[width=\textwidth]
	{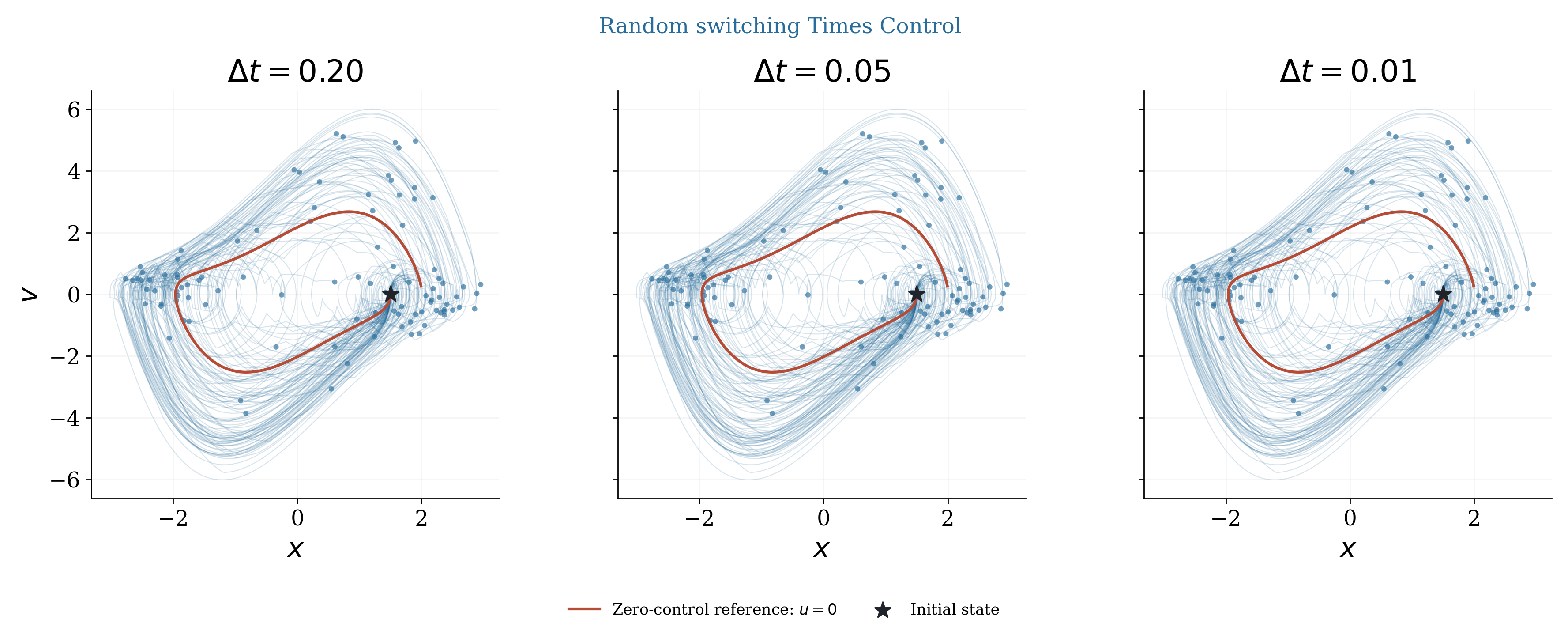}
	\caption{Random-switching controls for the Van der Pol
		oscillator. Switching times follow a Poisson process
		with rate $\lambda=2$, and the initial control and each
		subsequent control value are sampled independently and
		uniformly from $[-4,4]$. Blue curves show sampled state
		trajectories, and blue dots mark their terminal states.
		The red curve shows the zero-control dynamics; the star
		marks the initial state $(1.5,0)$.
		All panels use the same switching times and control values,
		horizon $T=6$, and axis limits.
		Here, $\Delta t$ specifies the time grid discretization.}
	\label{fig:vdp-poisson-switching}
\end{figure}

This construction preserves the control bounds and allows the
numerical integration step to be refined without changing the
switching rate.  Figure~\ref{fig:vdp-poisson-switching} illustrates the use of
random-switching controls for the Van der Pol oscillator.
The control is held constant between random switching times
and independently resampled from $[-4,4]$ at each switch.

\begin{figure}[htbp]
	\centering
	\includegraphics[width=\textwidth]
	{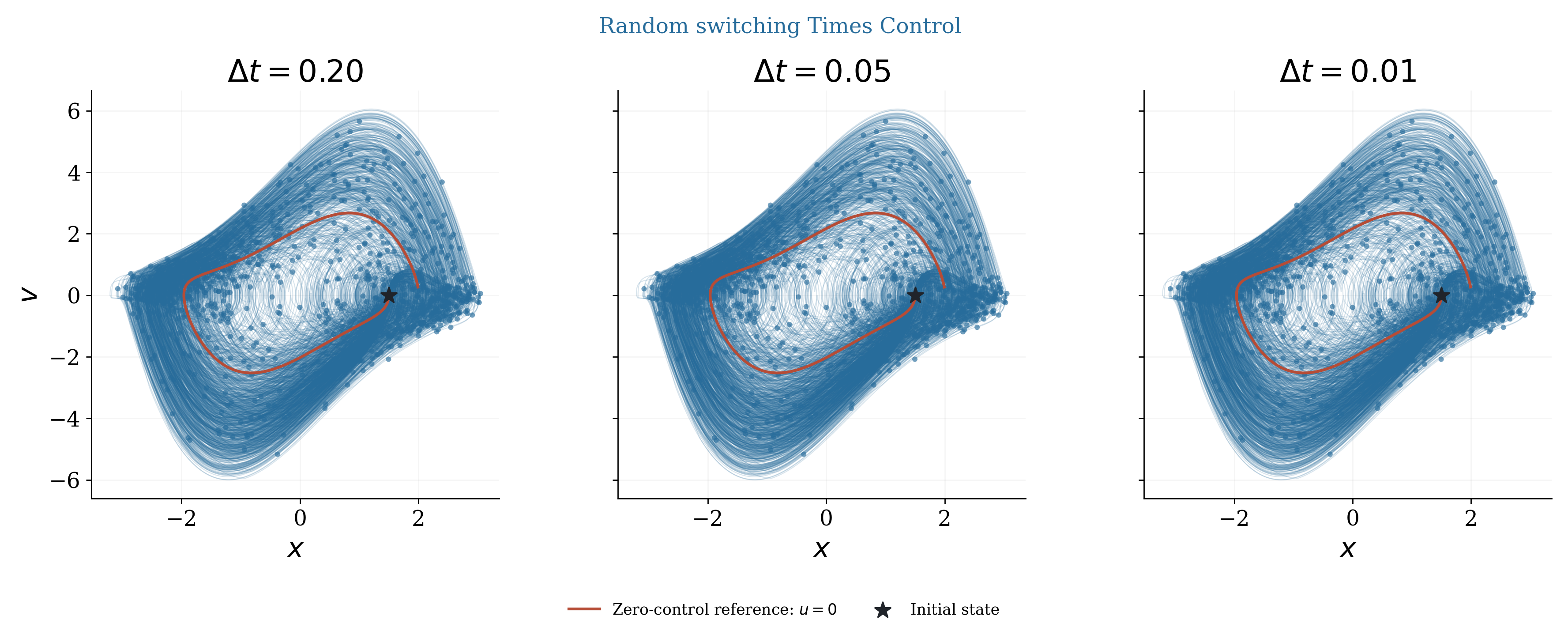}
	\caption{Random-switching controls experiment repeated with $N = 1000$ samples.}
	\label{fig:vdp-poisson-switching2}
\end{figure}

In the plots shown, one must interpret the results carefully. The samples from reachable sets are denoted by the markers, and not the traces of the trajectories. A standard way to improve coverage is to increase the number of samples. Figure~\ref{fig:vdp-poisson-switching2} illustrates this. The number of samples have been increased from $100$ to $1000$ to improve coverage. 

\section{Sampling using Optimal Transport}

The preceding approaches generate a distribution on the reachable
set by randomly sampling controls. Even when this distribution has
full support, it may place very little mass in some regions.
An alternative is to optimize the controls so that their
terminal states spread more uniformly over the reachable set \cite{elamvazhuthi2025uniform}. In this section, we review an approach to solve this problem using {\it optimal transport} \cite{villani2021topics}

Consider once again, the control-affine system
\[
\dot{x}
=
f_0(x)+\sum_{i=1}^{m}u_i g_i(x),
\qquad
x(0)=x_0,
\qquad
u(t)\in U,
\]
where $U \subset \mathbb{R}^m$ is compact and convex. Write
$\mathcal R_T(x_0)$ for its reachable set at time $T$ and assume
that this set has positive, finite volume. Ideally, we would
choose a probability law on admissible controls whose terminal
distribution is
\[
\rho_\star(x)
=
\frac{\mathbf{1}_{\mathcal R_T(x_0)}(x)}
{|\mathcal R_T(x_0)|}.
\]
This gives us an optimal transport problem from the initial point
distribution to the uniform distribution on the reachable set. In general optimal transport problems, the initial and final distributions can be more arbitrary but fixed. The difficulty,
of course, in our problem is that we do not know the reachable set. 

The key observation is that a probability density supported on
$\mathcal R_T(x_0)$ satisfies
\[
\int_{\R^n}\rho(x)^2\,\dd x
\geq
\frac{1}{|\mathcal R_T(x_0)|},
\]
with equality precisely for the uniform density. Thus, minimizing
the squared $L^2$ norm promotes uniformity without requiring an
explicit description of the set or its volume. One could also use other {\it entropy} like functionals, such as the KL-divergence, but we will stick to the $L^2$ norm here.

Therefore, uniform sampling can be approximated by the
{\it regularized} optimal transport problem
\[
\inf_{\mathbb{P}\in\mathcal{P}(\mathcal{U})}
\left\{
\int_{\mathcal{U}}
\left(\int_0^T |u(t)|^2\,\dd t\right)
\dd\mathbb{P}(u)
+
\frac{1}{\lambda}
\int_{\R^n}\rho_T(x)^2\,\dd x
\right\},
\]
where $\mathcal{U}$ is the set of admissible controls and
$\rho_T$ is the density of the terminal law
$(E_T)_\#\mathbb{P}$, with $E_T(u)=x^u(T)$, the end-point map that evaluates the solution at the final time for a given control $u$.

Under suitable assumptions, the terminal laws associated
with minimizers converge weakly to the uniform distribution
on the reachable set as $\lambda\downarrow0$. We call this optimal transport problem regularized when $\lambda >0$, and one recovers the classical optimal transport problem when $\lambda = 0$, which imposes fixed constraint on the terminal law.

To obtain a computational method, consider $N$ controlled
trajectories and their empirical terminal measure
\[
\mu_T^N
=
\frac{1}{N}\sum_{i=1}^{N}\delta_{x_i(T)}.
\]
This measure places mass $1/N$ at each terminal state and
therefore has no density with respect to Lebesgue measure.
In particular, we cannot directly evaluate its squared
$L^2$ density penalty.

The idea is to replace each point mass by a small, smooth
bump centered at that point. The function describing this
bump is called a \emph{kernel}. For example, we
can use the Gaussian kernel
\[
k_\delta(x)
=
\frac{1}{(2\pi\delta^2)^{n/2}}
\exp\!\left(-\frac{|x|^2}{2\delta^2}\right),
\]
where $\delta>0$ controls the width of the bump. This
kernel is nonnegative, symmetric, and integrates to one.

Replacing each point mass by such a bump gives the density
\[
\rho_{T,\delta}^N(x)
=
\frac{1}{N}\sum_{i=1}^{N}
k_\delta\bigl(x-x_i(T)\bigr)
=
(k_\delta*\mu_T^N)(x).
\]
The last expression is called the \emph{convolution} of
the kernel with the empirical measure. Here, it simply
means averaging the bumps centered at the terminal states.

\begin{figure}[htbp]
	\centering
	\includegraphics[width=\textwidth]{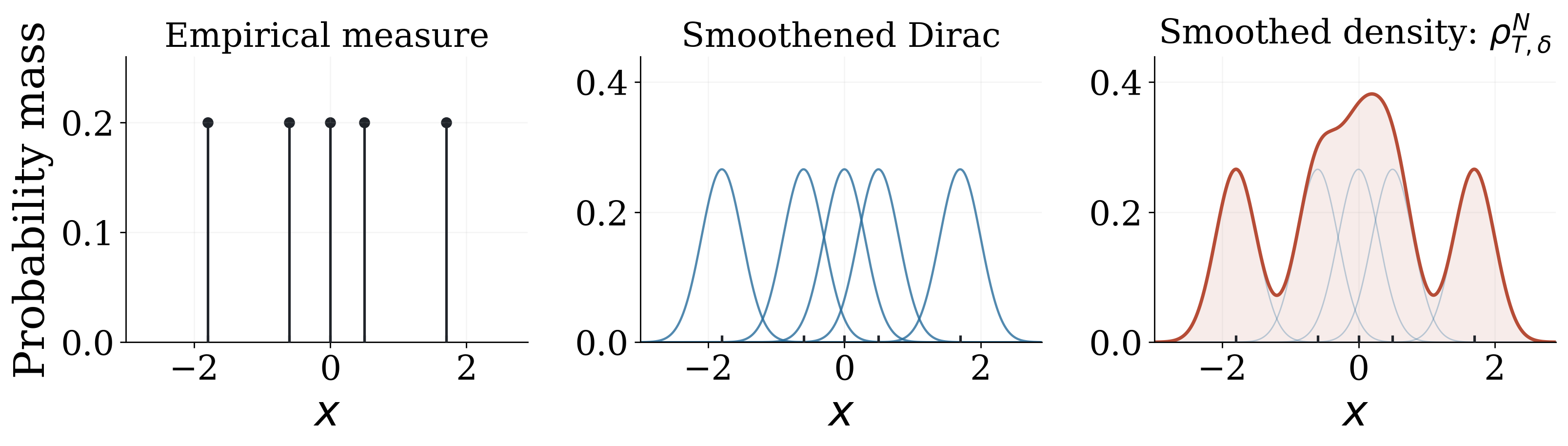}
	\caption{Kernel smoothing of an empirical measure.
		Left: point masses at the sampled terminal states.
		Center: Gaussian bumps centered at these states, each
		with integral $1/N$.
		Right: the resulting density
		$\rho_{T,\delta}^N=k_\delta*\mu_T^N$ (red),
		obtained by summing the weighted bumps (blue).}
	\label{fig:kernel-smoothing}
\end{figure}

Figure~\ref{fig:kernel-smoothing} illustrates the idea of convolving empirical measure with a smooth kernel. Each point mass is replaced by
a Gaussian bump of mass $1/N$. Summing these bumps gives
a smooth probability density.

\begin{figure}[htbp]
	\centering
	\includegraphics[width=\textwidth]{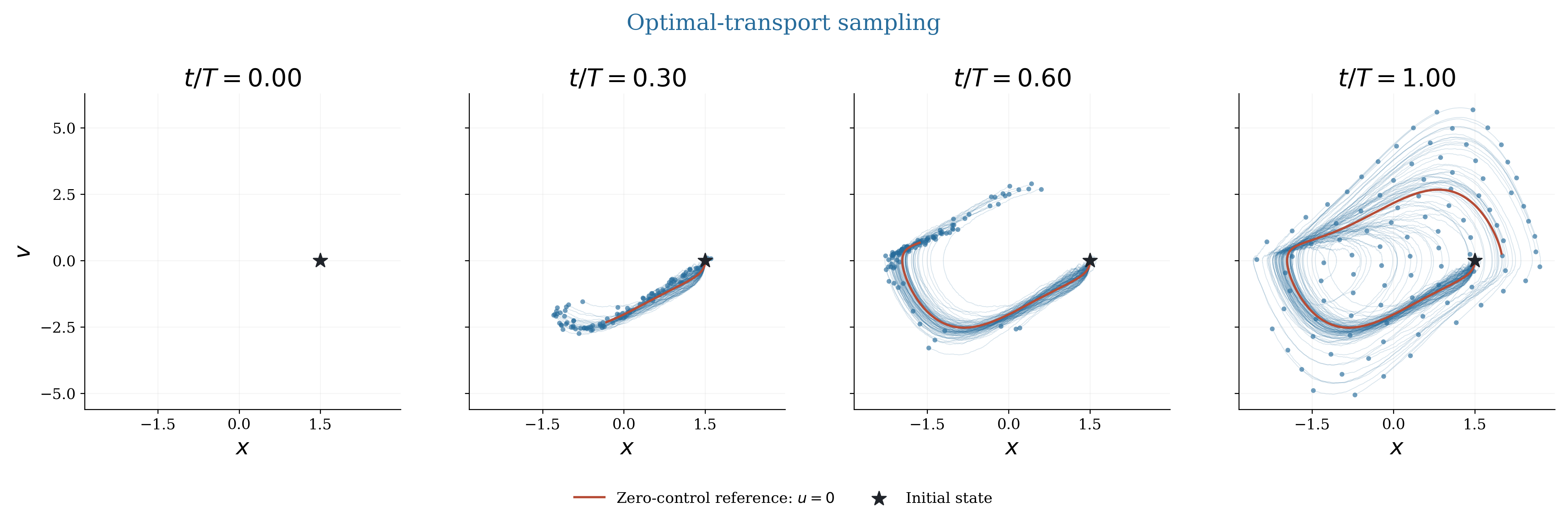}
	\caption{Optimal-transport-based reachable-set sampling for
		the Van der Pol oscillator with $\mu=1$ and $|u|\leq4$.
		Panels show the states at $t/T=0$, $0.3$, $0.6$, and $1$.
		Blue points indicate current states, and faint blue curves
		show their accumulated paths. The red curve shows the
		zero-control trajectory up to the displayed time, and the
		star marks the common initial state $(1.5,0)$.
		The terminal interaction cost promotes dispersion without
		requiring prior knowledge of the reachable set.}
	\label{fig:vdp-ot-sampling}
\end{figure}

We can now evaluate
\[
\int_{\R^n}\bigl|\rho_{T,\delta}^N(x)\bigr|^2\,\dd x.
\]
When terminal states cluster together, their bumps overlap
strongly and this penalty increases. Minimizing it therefore
encourages the terminal states to spread apart, subject to
the dynamics and control constraints.
Writing $K_\delta=k_\delta*k_\delta$, we have
\[
\int_{\R^n}|k_\delta*\mu_T^N(x)|^2\,\dd x
=
\frac{1}{N^2}
\sum_{i,j=1}^{N}
K_\delta\bigl(x_i(T)-x_j(T)\bigr).
\]
This leads to the coupled optimal control problem
\[
\begin{aligned}
	\min_{u_1,\dots,u_N}\quad&
	\frac{1}{N}\sum_{i=1}^{N}
	\int_0^T |u_i(t)|^2\,\dd t
	+
	\frac{1}{\varepsilon N^2}
	\sum_{i,j=1}^{N}
	K_\delta\bigl(x_i(T)-x_j(T)\bigr),\\
	\text{subject to}\quad&
	\dot{x}_i(t)=f\bigl(x_i(t),u_i(t)\bigr),\\
	&
	x_i(0)=x_0,
	\qquad u_i(t)\in U.
\end{aligned}
\]
This problem can be solved using any nonlinear optimization solver. We use Pontryagin's maximum principle \cite{liberzon2011calculus} to design a gradient descent algorithm to construct the updates of the above algorithm.  Figure~\ref{fig:vdp-ot-sampling} illustrates optimal-transport-based
sampling for the controlled Van der Pol oscillator considered
earlier. Starting from $(x(0),v(0))=(1.5,0)$, the controls
$u(t)\in[-4,4]$ are optimized improve the $L^2$ norm of the distribution at $T=6$. Unlike independent
random sampling, this algorithm adjusts the controls according
to the distribution of their endpoints, encouraging a more
uniform spread within the reachable set. Even with $N= 100$, one can see a more uniform coverage of the reachable set than the stochastic approaches considered in the previous section. Interestingly, the trajectories at intermediate trajectory traces follow the mean closely over most of the time horizon. 

However, one should note that the quality of the results depend a lot on the hyperparameters chosen to solve the optimization problem. One key parameter choice is the bandwidth $\delta$ of the kernel function.  Moreover, due to the inherent non-convexity of the problem we are never truly guaranteed to find a good solution.

\paragraph{Further reading.}

The relation between rapidly varying randomized controls and averaged continuous-time dynamics is part of the theory of stochastic approximation \cite{borkar2008stochastic}. The pieceise deterministic process approach to sampling is inspired from \cite{benaim2015qualitative}.
Details about the optimal transport approach to sample from the reachable set can be found in \cite{elamvazhuthi2025uniform}.

\chapter{Controlled Normalizing Flows}
\label{chap:controlled-normalizing-flows}

Flow matching in Chapter~\ref{chap:flow-matching} learns a feedback control law from trajectory data. Before flow matching, one of the popular approaches to sample from a target density $\rho_\star$ was to pose it as an constrained optimization problem. Continuous normalizing flows use the same continuity equation from Chapter~\ref{chap:measure-transport}, but optimize the law of the maximum likelihood \cite{chen2018neural}.

In this chapter, we present a continuous normalizing flow approach for controlling the continuity equation toward a target probability distribution. This method is especially relevant when the desired terminal state is uncertain and is therefore represented by a probability distribution rather than by a single point. When the target distribution concentrates at a single point, the resulting optimization resembles a continuous-time policy-gradient method.

\section{Maximum Likelihood Estimation}
\label{sec:cnf-mle}

Before using dynamics, we first recall the statistical objective that will identify the flow. This is the only new ingredient relative to the continuity-equation viewpoint of Chapter~\ref{chap:measure-transport}.

Let $\rho_\star(x)$ denote the true data density on $\R^n$, and let $\rho_\theta(x)$ be a parameterized model density.

One way to fit the model is to minimize the discrepancy between $\rho_\star$ and $\rho_\theta$ using the {\it Kullback--Leibler (KL) divergence}:
\begin{equation}
    \mathrm{KL}(\rho_\star\,\|\,\rho_\theta)
    :=
    \int_{\R^n}
    \rho_\star(x)
    \log\!\left(
        \frac{\rho_\star(x)}{\rho_\theta(x)}
    \right)
    \dd x.
\end{equation}

Expanding the logarithm gives
\begin{align}
    \mathrm{KL}(\rho_\star\,\|\,\rho_\theta)
    &=
    \int_{\R^n}\rho_\star(x)\log\rho_\star(x)\dd x
    -
    \int_{\R^n}\rho_\star(x)\log\rho_\theta(x)\dd x.
\end{align}

The first term does not depend on $\theta$. Therefore,
\begin{equation}
    \operatorname*{arg\,min}_{\theta}
    \mathrm{KL}(\rho_\star\,\|\,\rho_\theta)
    =
    \operatorname*{arg\,max}_{\theta}
    \int_{\R^n}
    \rho_\star(x)\log\rho_\theta(x)\dd x.
\end{equation}

The expression on the right-hand side is the expected log-likelihood under the data distribution. Hence,
\begin{equation}
    \text{maximizing likelihood}
    \quad\Longleftrightarrow\quad
    \text{minimizing }
    \mathrm{KL}(\rho_\star\,\|\,\rho_\theta).
\end{equation}

Our goal is to obtain a sample-based approximation of the objective. Given independent and identically distributed samples
\[
    x_1,\dots,x_N\sim\rho_\star,
\]
the expected log-likelihood can be approximated by the empirical average:
\begin{equation}
    \max_\theta
    \int_{\R^n}\rho_\star(x)\log\rho_\theta(x)\dd x
    \approx
    \max_\theta
    \frac{1}{N}
    \sum_{i=1}^{N}
    \log\rho_\theta(x_i).
\end{equation}

\section{Continuous Normalizing Flows}
\label{sec:continuous-normalizing-flows}

We first consider the simple control system
\begin{equation}
    \dot{x}(t)=u(t,x(t)).
\end{equation}
Our objective is to choose the vector field $u$ so that the induced flow transports $\rho_0$ to $\rho_\star$. From the last section, we can pose this as likelihood maximization problem.

Let $\rho(t,x)$ denote the density of the state. Assuming that $\rho(t,x(t))>0$, the chain rule gives
\begin{align}
    \frac{\dd}{\dd t}\log\rho(t,x(t))
    &=
    \frac{
        \partial_t\rho(t,x(t))
        +
        u(t,x(t))\cdot\nabla_x\rho(t,x(t))
    }{
        \rho(t,x(t))
    }.
    \label{eq:cnf-log-density-chain-rule}
\end{align}

Positivity of the density is preserved under a sufficiently smooth invertible flow. More generally, if $Z$ has density $\rho_Z$, $X=f(Z)$, and
\[
    f:\R^n\to\R^n
\]
is an invertible $C^1$ map with a $C^1$ inverse, then the change-of-variables formula states that
\begin{equation}
\boxed{  \rho_X(x)
    =
    \rho_Z\!\left(f^{-1}(x)\right)
    \left|
        \det\!\left(
            \frac{\partial f^{-1}(x)}{\partial x}
        \right)
    \right|.}
    \label{eq:cnf-static-change-of-variables}
\end{equation}

Recall that $\rho$ satisfies the continuity equation
\begin{equation}
    \partial_t\rho
    =
    -\nabla\cdot(u\rho).
\end{equation}
Expanding the divergence yields
\begin{equation}
    \partial_t\rho
    =
    -u\cdot\nabla\rho
    -
    \rho\,\nabla\cdot u.
\end{equation}
Substituting this identity into \eqref{eq:cnf-log-density-chain-rule} gives
\begin{equation}
    \frac{\dd}{\dd t}\log\rho(t,x(t))
    =
    -\nabla\cdot u(t,x(t)).
    \label{eq:cnf-instantaneous-change-of-variables}
\end{equation}

Integrating \eqref{eq:cnf-instantaneous-change-of-variables} from $0$ to $T$ gives the continuous change-of-variables formula
\begin{equation}
    \log\rho(T,x(T))
    =
    \log\rho(0,x(0))
    -
    \int_0^T
    \nabla\cdot u(t,x(t))
    \dd t.
    \label{eq:cnf-continuous-change-of-variables}
\end{equation}

Replace $u(t,x)$ by a parameterized vector field $u_\theta(t,x)$. For each data point $x_i\sim\rho_\star$, let $x_i(t)$ solve the terminal-value problem
\begin{equation}
    \dot{x}_i(t)
    =
    u_\theta(t,x_i(t)),
    \qquad
    x_i(T)=x_i,
    \label{eq:cnf-backward-ode}
\end{equation}
integrated backward from time $T$ to time $0$.

The continuous change-of-variables formula gives
\begin{equation}
    \log\rho_\theta(T,x_i)
    =
    \log\rho_0(x_i(0))
    -
    \int_0^T
    \nabla\cdot
    u_\theta(t,x_i(t))
    \dd t.
    \label{eq:cnf-sample-log-likelihood}
\end{equation}

Therefore, the sample-based maximum-likelihood objective is
\begin{equation}
    \max_\theta
    \frac{1}{N}
    \sum_{i=1}^{N}
    \left[
        \log\rho_0(x_i(0))
        -
        \int_0^T
        \nabla\cdot
        u_\theta(t,x_i(t))
        \dd t
    \right].
    \label{eq:cnf-mle-objective}
\end{equation}

Equivalently, one minimizes the negative log-likelihood
\begin{equation}
    \mathcal{L}_{\mathrm{CNF}}(\theta)
    =
    -\frac{1}{N}
    \sum_{i=1}^{N}
    \left[
        \log\rho_0(x_i(0))
        -
        \int_0^T
        \nabla\cdot
        u_\theta(t,x_i(t))
        \dd t
    \right].
\end{equation}

This gives us the following algorithm.

\begin{tcolorbox}[algorithmbox,title=\textbf{Algorithm: Continuous Normalizing Flow Training}]
\begin{enumerate}[label=\arabic*.]
    \item Initialize the parameters $\theta$.

    \item Sample a batch of terminal data points
    \[
        x_i\sim\rho_\star,
        \qquad
        i=1,\dots,N.
    \]

    \item For each $x_i$, solve backward in time:
    \[
        \dot{x}_i(t)
        =
        u_\theta(t,x_i(t)),
        \qquad
        x_i(T)=x_i.
    \]

    \item Along each trajectory, compute
    \[
        s_i
        =
        \int_0^T
        \nabla\cdot
        u_\theta(t,x_i(t))
        \dd t.
    \]

    \item Evaluate the loss
    \[
        \mathcal{L}_{\mathrm{CNF}}(\theta)
        =
        -\frac{1}{N}
        \sum_{i=1}^{N}
        \left[
            \log\rho_0(x_i(0))-s_i
        \right].
    \]

    \item Update the parameters:
    \[
        \theta
        \leftarrow
        \theta
        -
        \eta\nabla_\theta
        \mathcal{L}_{\mathrm{CNF}}(\theta).
    \]

    \item Repeat Steps 2--6 until convergence.
\end{enumerate}
\end{tcolorbox}

In high-dimensional systems, evaluating
\[
    \nabla\cdot u_\theta(t,x)
    =
    \operatorname{tr}
    \left(
        \frac{\partial u_\theta(t,x)}{\partial x}
    \right)
\]
can be computationally expensive. A common alternative is Hutchinson's trace estimator, which expresses the trace as an expectation of random quadratic forms. If $\varepsilon$ is a random vector satisfying
\[
    \E[\varepsilon]=0,
    \qquad
    \E[\varepsilon\varepsilon^{\mathsf T}]=I,
\]
then
\begin{equation}
    \operatorname{tr}(A)
    =
    \E\!\left[
        \varepsilon^{\mathsf T}A\varepsilon
    \right].
\end{equation}
Consequently,
\begin{equation}
    \nabla\cdot u_\theta(t,x)
    =
    \E_\varepsilon
    \left[
        \varepsilon^{\mathsf T}
        \frac{\partial u_\theta(t,x)}{\partial x}
        \varepsilon
    \right],
\end{equation}
which can be approximated using one or more random samples of $\varepsilon$.

\begin{figure}[htbp]
	\centering
	\includegraphics[width=\textwidth]
	{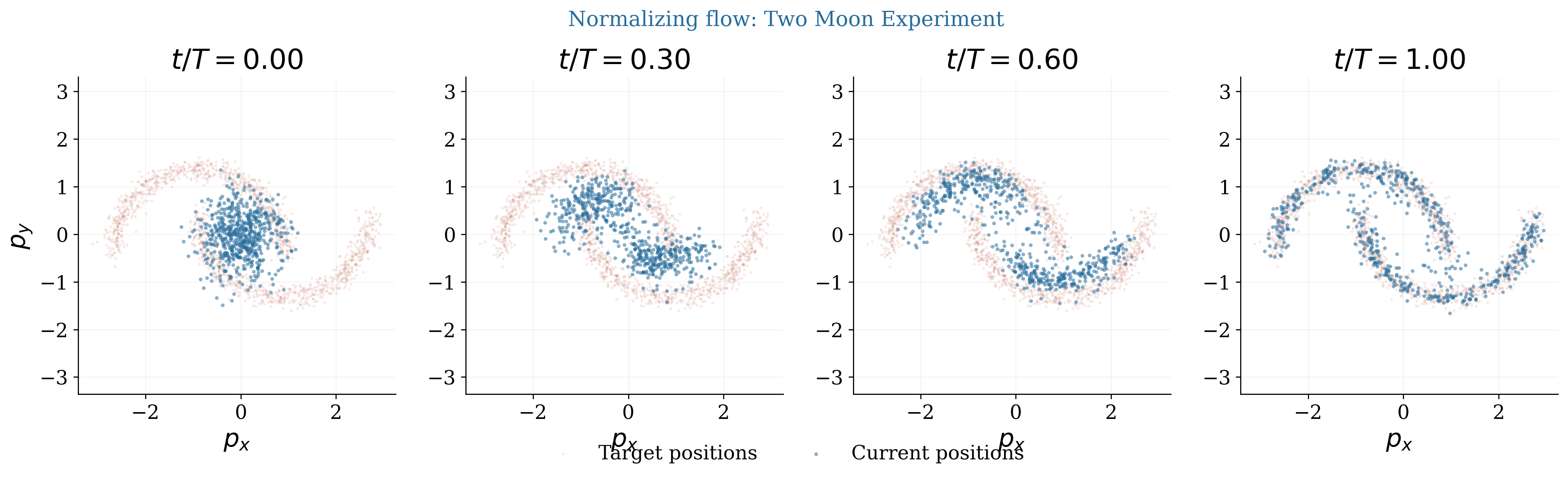}
	\caption{Continuous normalizing flow for the fully actuated
		planar system. Blue points show the evolving model samples,
		initialized from a Gaussian distribution, while faint red
		points show samples from the noisy two-moons target
		distribution. Snapshots are shown at
		\(t/T=0\), \(0.3\), \(0.6\), and \(1\), with identical axis
		limits. At the terminal time, the learned flow approximately
		reproduces the two-moons distribution.}
	\label{fig:integrator-cnf}
\end{figure}

We test the continuous normalizing-flow with initial law as Gaussian, and the target law is the noisy two-moons distribution. The velocity field \(u^\theta(t,x)\) is
trained by maximizing the likelihood of samples from the target
distribution. Figure~\ref{fig:integrator-cnf} shows the resulting
evolution of samples from the initial Gaussian distribution.

\section{Controlled Normalizing Flows}
\label{sec:controlled-normalizing-flow}

The unconstrained flow above does not impose any constraints on how samples can move during their evolution.  We now generalize the construction to a control-affine system
\begin{equation}
    \dot{x}
    =
    f\bigl(x,u^\theta(t,x)\bigr)
    :=
    f_0(x)
    +
    \sum_{j=1}^{m}
    u^\theta_j(t,x)f_j(x).
    \label{eq:control-cnf-system}
\end{equation}

The induced density now satisfies
\begin{equation}
    \partial_t\rho
    +
    \nabla\cdot
    \left(
        f\bigl(x,u^\theta(t,x)\bigr)\rho
    \right)
    =
    0.
\end{equation}
Along a closed-loop trajectory, the log-density therefore evolves according to
\begin{equation}
    \frac{\dd}{\dd t}\log\rho(t,x(t))
    =
    -
    \nabla\cdot
    f\bigl(x(t),u^\theta(t,x(t))\bigr).
    \label{eq:control-cnf-log-density}
\end{equation}

For each terminal sample $x_i\sim\rho_\star$, solve
\begin{equation}
    \dot{x}_i(t)
    =
    f\bigl(x_i(t),u^\theta(t,x_i(t))\bigr),
    \qquad
    x_i(T)=x_i,
    \label{eq:control-cnf-backward-system}
\end{equation}
backward in time. The corresponding model log-density is
\begin{equation}
    \log\rho_\theta(T,x_i)
    =
    \log\rho_0(x_i(0))
    -
    \int_0^T
    \nabla\cdot
    f\bigl(x_i(t),u^\theta(t,x_i(t))\bigr)
    \dd t.
    \label{eq:control-cnf-log-likelihood}
\end{equation}

Thus, the controlled CNF loss is
\begin{equation}
    \mathcal{L}_{\mathrm{control\text{-}CNF}}(\theta)
    =
    -\frac{1}{N}
    \sum_{i=1}^{N}
    \left[
        \log\rho_0(x_i(0))
        -
        \int_0^T
        \nabla\cdot
        f\bigl(x_i(t),u^\theta(t,x_i(t))\bigr)
        \dd t
    \right].
    \label{eq:control-cnf-loss}
\end{equation}
The construction assumes that the reference density $\rho_0$ is available in closed form, so that $\log\rho_0(x_i(0))$ can be evaluated. 
Given these derivations, we have the following algorithm.

\begin{tcolorbox}[algorithmbox,title=\textbf{Algorithm: Controlled Normalizing Flow Training}]
\begin{enumerate}[label=\arabic*.]
    \item Initialize the controller parameters $\theta$.

    \item Sample terminal data points
    \[
        x_i\sim\rho_\star,
        \qquad
        i=1,\dots,N.
    \]

    \item For each terminal sample, solve backward:
    \[
        \dot{x}_i(t)
        =
        f\bigl(x_i(t),u^\theta(t,x_i(t))\bigr),
        \qquad
        x_i(T)=x_i.
    \]

    \item Compute
    \[
        s_i
        =
        \int_0^T
        \nabla\cdot
        f\bigl(x_i(t),u^\theta(t,x_i(t))\bigr)
        \dd t.
    \]

    \item Evaluate
    \[
        \mathcal{L}_{\mathrm{control\text{-}CNF}}(\theta)
        =
        -\frac{1}{N}
        \sum_{i=1}^{N}
        \left[
            \log\rho_0(x_i(0))-s_i
        \right].
    \]

    \item Update
    \[
        \theta
        \leftarrow
        \theta
        -
        \eta
        \nabla_\theta
        \mathcal{L}_{\mathrm{control\text{-}CNF}}(\theta).
    \]

    \item Repeat Steps 2--6 until convergence.
\end{enumerate}
\end{tcolorbox}

We illustrate controlled normalizing flows using the four-dimensional
kinematic vehicle model
\begin{align}
	\dot p_x &= v\cos\theta, &
	\dot p_y &= v\sin\theta,\\
	\dot\theta &= \omega, &
	\dot v &= a,
\end{align}
where the neural feedback controller generates the yaw rate
$\omega\in[-3,3]$ and acceleration $a\in[-3,3]$. This is just the unicycle model that has been turned into a system with drift. Instead of translational motion being affected directly by $v$, we now control $a$. 
The state is initialized from a centered Gaussian distribution
with covariance $
\operatorname{diag}(0.45^2,0.45^2,0.8^2,0.6^2).$
The target position distribution again consists of the two moons distribution. To define the full four-dimensional target law required by
the likelihood objective, we augment these positions with
independent Gaussian heading and speed variables having
standard deviations $0.8$ and $0.6$, respectively.
The heading is treated as an unwrapped real coordinate.

\begin{figure}[htbp]
	\centering
	\includegraphics[width=\textwidth]
	{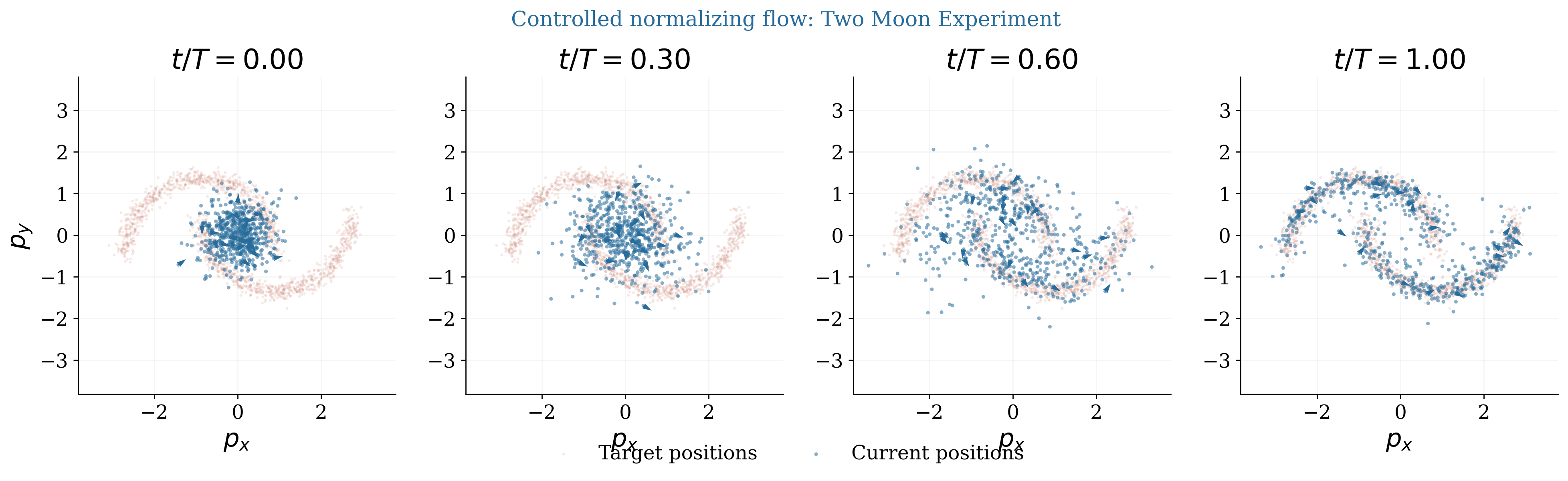}
	\caption{Controlled normalizing flow for a four-dimensional
		kinematic vehicle model. Blue points show current positions,
		and faint red points show samples from the two-moons target
		position distribution. Arrows indicate the heading of a
		subset of vehicles, rather than their velocity.
		Snapshots are shown at $t/T=0$, $0.3$, $0.6$, and $1$,
		with $T=4$ and identical axis limits.
		Only the position coordinates are displayed}
	\label{fig:car-controlled-cnf}
\end{figure}

Figure~\ref{fig:car-controlled-cnf} shows the position marginal
evolving under the learned feedback over the horizon $T=4$.
The controller is trained using the controlled normalizing-flow
likelihood objective, and the displayed trajectories are evaluated
from fresh Gaussian initial samples.

\section{Reachable Set Sampling using Controlled Normalizing Flows}
\label{sec:cnf-reachable-sampling}

The same constructions considered in the previous section can be used for sampling from the reachable set. In the previous chapter, we consider random sampling based strategies to achieve this, and also looked at an optimal transport based method. Similar to the optimal transport based method, in the normalizing flow setting, instead of prescribing a target density $\rho_\star$,
we can optimize for a feedback law to maximize the entropy of the terminal
state distribution. This alternative strategy forces samples generated to more uniformly sample over the reachable set.

Firstly, one can see that this problem is not well posed, for a single initial condition as a deterministic, well-posed feedback flow starting from a single
point generates only one trajectory. We therefore initialize the
state uniformly on a small ball around $x_0$:
\begin{equation}
	\rho_0(x)
	=
	\frac{\mathbf{1}_{B_\varepsilon(x_0)}(x)}
	{|B_\varepsilon(x_0)|},
	\qquad \varepsilon>0,
\end{equation}
where $|B_\varepsilon(x_0)|$ denotes its volume.

Consider the controlled system
\begin{equation}
	\dot{x}(t)
	=
	f\bigl(x(t),u^\theta(t,x(t))\bigr),
	\qquad
	u^\theta(t,x)\in U,
\end{equation}
with $f$ as defined in \eqref{eq:control-cnf-system}.
Let $\rho_\theta(t,x)$ denote the density of the state variable.
Our objective is to maximize its terminal  entropy:
\begin{equation}
	\max_\theta H\bigl(\rho_\theta(T,\cdot)\bigr),
	\qquad
	H(\rho)
	:=
	-\int_{\R^n}\rho(x)\log\rho(x)\,\dd x.
	\label{eq:control-cnf-entropy-objective}
\end{equation}

Under assumptions ensuring a sufficiently smooth invertible flow
and integrability of the expressions below, the continuous
change-of-variables formula gives
\begin{equation}
	\log\rho_\theta(T,x(T))
	=
	\log\rho_0(x(0))
	-
	\int_0^T
	\nabla\cdot
	f\bigl(x(t),u^\theta(t,x(t))\bigr)
	\,\dd t.
\end{equation}
This identity holds for $\rho_0$-almost every initial state;
the discontinuity of the uniform density at the boundary of
the ball does not affect the change-of-variables calculation.

Taking the negative expectation over $x(0)\sim\rho_0$ yields
\begin{align}
	H\bigl(\rho_\theta(T,\cdot)\bigr)
	&=
	H(\rho_0)
	+
	\E_{x(0)\sim\rho_0}
	\left[
	\int_0^T
	\nabla\cdot
	f\bigl(x(t),u^\theta(t,x(t))\bigr)
	\,\dd t
	\right],\\
	H(\rho_0)
	&=
	\log|B_\varepsilon(x_0)|.
\end{align}

The initial entropy does not depend on $\theta$. Therefore,
maximizing terminal entropy is equivalent to maximizing the
expected accumulated divergence along forward trajectories.
Given samples $x_i(0)\sim\rho_0$, we minimize
\begin{equation}
	\mathcal{L}_{\mathrm{entropy}}(\theta)
	=
	-\frac{1}{N}
	\sum_{i=1}^{N}
	\int_0^T
	\nabla\cdot
	f\bigl(x_i(t),u^\theta(t,x_i(t))\bigr)
	\,\dd t.
	\label{eq:control-cnf-entropy-loss}
\end{equation}
Unlike the likelihood construction, this objective uses initial
samples and forward integration. The divergence is taken with respect to the complete closed-loop
vector field. In particular,
\begin{equation}
	\nabla\cdot f\bigl(x,u^\theta(t,x)\bigr)
	=
	\nabla\cdot f_0(x)
	+
	\sum_{j=1}^{m}
	\left[
	\nabla_x u^\theta_j(t,x)\cdot f_j(x)
	+
	u^\theta_j(t,x)\nabla\cdot f_j(x)
	\right].
\end{equation}

There is some intuition as to why this encourages uniform reachable-set sampling. Every generated sample belongs to the reachable set from the
initial region,
\begin{equation}
	\mathcal R_T(B_\varepsilon(x_0);U)
	:=
	\bigcup_{y\in B_\varepsilon(x_0)}
	\mathcal R_T(y;U).
\end{equation}
Among all probability densities supported on this set, the
uniform density maximizes entropy. Writing
$S=\mathcal R_T(B_\varepsilon(x_0);U)$, we have
\begin{equation}
	\mathrm{KL}
	\left(
	\rho\,\middle\|\,
	\frac{\mathbf{1}_S}{|S|}
	\right)
	=
	\log|S|-H(\rho).
\end{equation}
The unknown volume contributes only a constant. Thus, the entropy
objective can be evaluated and optimized without computing the
reachable set. This identity motivates uniform sampling, but {\it does not guarantee}
that the learned flow covers the entire reachable set. 

The implementation of this idea is shown in the following algorithm.

\begin{tcolorbox}[algorithmbox,title=\textbf{Algorithm: Entropy-Maximizing Controlled Normalizing Flow}]
	\begin{enumerate}[label=\arabic*.]
		\item Parameterize the feedback controller by a neural network
		\[
		u^\theta:[0,T]\times\R^n\to U,
		\]
		whose output is confined to the admissible control set $U$
		by construction, and initialize its parameters $\theta$.
		\item Sample initial states
		\[
		x_i(0)\sim\rho_0,
		\qquad
		i=1,\dots,N.
		\]
		
		\item For each initial sample, solve forward:
		\[
		\dot{x}_i(t)
		=
		f\bigl(x_i(t),u^\theta(t,x_i(t))\bigr).
		\]
		
		\item Along each trajectory, compute
		\[
		s_i
		=
		\int_0^T
		\nabla\cdot
		f\bigl(x_i(t),u^\theta(t,x_i(t))\bigr)
		\,\dd t.
		\]
		
		\item Evaluate the loss
		\[
		\mathcal{L}_{\mathrm{entropy}}(\theta)
		=
		-\frac{1}{N}\sum_{i=1}^{N}s_i.
		\]
		
		\item Update the parameters:
		\[
		\theta
		\leftarrow
		\theta
		-
		\eta\nabla_\theta
		\mathcal{L}_{\mathrm{entropy}}(\theta).
		\]
		
		\item Repeat Steps 2--6 until convergence.
	\end{enumerate}
\end{tcolorbox}

We revisit the controlled Van der Pol oscillator from
Chapter~\ref{chap:not-sample-reachable},
\begin{align}
	\dot x &= v,\\
	\dot v &= (1-x^2)v-x+u,
	\qquad u\in[-4,4].
\end{align}
In that chapter, we explored the effect of different randomized
control constructions on state-space exploration. Here, we learn
a feedback controller by maximizing the terminal state entropy.

\begin{figure}[htbp]
	\centering
	\includegraphics[width=\textwidth]
	{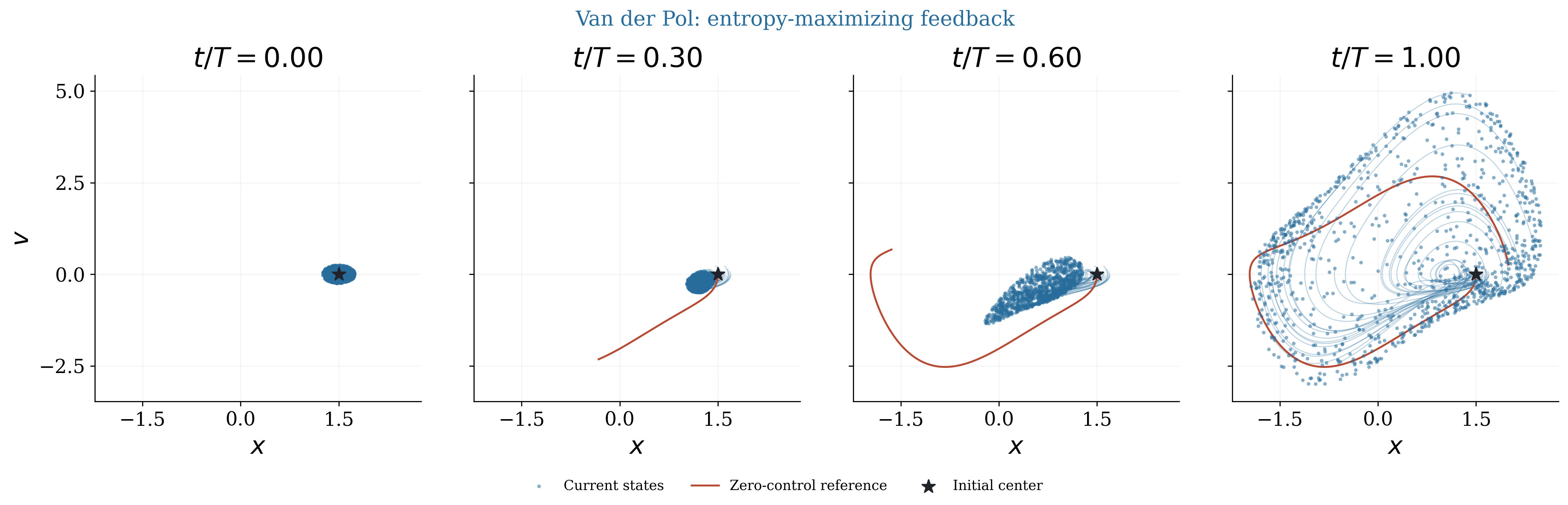}
	\caption{Entropy-maximizing controlled normalizing flow
		for the Van der Pol oscillator.
		Initial states are sampled uniformly from a disk of
		radius $0.25$ centered at $(1.5,0)$.
		Blue points show the current states under the learned
		feedback $u^\theta\in[-4,4]$, and blue curves show
		accumulated paths for a subset of particles.
		The red curve shows the zero-control trajectory from
		the initial center, marked by the star.
		Snapshots are shown at $t/T=0$, $0.3$, $0.6$, and $1$,
		with $T=6$ and identical axis limits.}
	\label{fig:vdp-entropy-cnf}
\end{figure}

We initialize the state uniformly on a disk of radius $0.25$
centered at $(1.5,0)$ and fix the horizon $T=6$. The feedback
$u^\theta(t,x,v)$ is represented by a neural network whose output
is confined to $[-4,4]$.
We train the controller using
\eqref{eq:control-cnf-entropy-loss}, without requiring the
reachable-set boundary or samples from a prescribed target law.
Figure~\ref{fig:vdp-entropy-cnf} shows the resulting evolution
from fresh initial samples under the learned feedback.

Compared with the randomized-control constructions and optimal transport based method in
Chapter~\ref{chap:not-sample-reachable}, the entropy-maximizing
approach learns a deterministic feedback that transports an
initial distribution toward a more skewed terminal law. Though our heuristic suggested that we should get a more uniform distribution over the reachable set, the method is not able to sample from points in the lower regions of the reachable set. One hypothesis is that this is due to the initial distribution being sampled from a ball, instead of a single point. Alternatively, the non-convexity of the training problem is finding sub-optimal solutions that do not achieve a more uniform spread. Lastly, it could just be that the optimal distribution is more skewed towards regions over the attractor. The constraint that the flow is constructed from a feedback law, rather than a general probability measure, could be another obstacle to the method. Excuses, excuses!

 \paragraph{Further Reading}
 
 Normalizing flows were introduced in \cite{chen2018neural}. The application of normalizing flows for reachability analysis has not been considered the literature, as far as the author knows. But it is an interesting alternative to the optimal transport based method introduced in the previous chapter, as it does not involves kernel density estimates, which are harder to scale to high dimensions.

\chapter[Fokker--Planck Equations]{Fokker--Planck Equations\\and Stabilization of Probability Densities}
\label{chap:fokker-planck}

The previous chapters used deterministic flows and the continuity equation to achieve corresponding control and sampling objectives. We now add stochasticity into the picture. While we have already seen stochasticity by exciting controls systems using Brownian paths, we will use its derivative, white noise, as our source of stochasticity. Although individual sample paths may be irregular, their distribution evolves according to a partial differential equation (PDE). The {\it Fokker--Planck PDE} is the stochastic analogue of the continuity equation in this setting.  The Fokker-Planck equation (FPK) describes the deterministic evolution of the law of a {\it stochastic differential equation} and will provide the forward ``noising'' dynamics used later in Chapter~\ref{chap:denoising-diffusion} on diffusion models.

 Beyond diffusion models, understanding the the long-term behavior of these equations have a number of applications in uncertainty shaping, probabilistic control, swarm coordination, ergodic exploration. This gives white noise based perturbations of the system a significant advantage over the randomization in chapter on flow matching. Particularly, one can exactly characterize the long term behavior of the FPK, rather than just talk about the support properties of the process.
 
 While there is a large industry on understanding the long term behavior of Markov processes and their associated Fokker-Planck equations \cite{bakry2014analysis}, we will take only a very small peek into this area. Lastly, our goal will be to solve the inverse stabilization problem as is desired in control: Given the target steady state behavior, can we design control laws that achieve this objective? While this answer is very challenging in the deterministic setting, surprisingly in the stochastic case, one can construct \textbf{closed-form expressions} to solve the problem for some system classes.

\section{Classical Diffusion and Long Term Behavior}
\label{sec:classical-fp}

We start with the simplest setting first to motivate the subsequent developments. Suppose we have the control system,
\[\dot{x} = u(t)\]
Our goal is to design a feedback law $u(x)$, independent of the law of $x(0)$, so that the law of the process $x(t)$ will asymptotically converge to some given target density $\rho_\star(x)$ as $t \rightarrow \infty$. As such, this problem is very challenging to solve. We will consider if we can solve a stochastic version of this problem.

We introduce a {\it reflected diffusion process} \cite{ikeda2014stochastic}. Let $\Omega\subset\R^n$ be a bounded connected domain with sufficiently regular boundary. Let $n_{\mathrm{in}}(x)$ denote the inward-pointing unit normal. Consider the reflected It\^o SDE
\begin{equation}
    \dd X_t
    =
    -\nabla U(X_t)\dd t
    +
    \sqrt{2}\,\dd W_t
    +
    n_{\mathrm{in}}(X_t)\dd L_t,
    \label{eq:reflected-langevin}
\end{equation}
where $W_t$ is a $n$-dimensional Brownian motion and $L_t$ is a nondecreasing boundary local-time process that increases only when $X_t\in\partial\Omega$. The reflection term keeps $X_t$ inside $\overline{\Omega}$. We work in this setting, rather than the standard case  of $\Omega = \mathbb{R}^n$ to avoid issues of compactness and due to which some of the claims made in this chapter might be ``less true'' from a mathematically rigorous standpoint. The equation \eqref{eq:reflected-langevin}  above is sometimes referred to as the {\it Overdamped Langevin equation}.

The choice \[-\nabla U(X_t)\dd t
+
\sqrt{2}\,\dd W_t
+
n_{\mathrm{in}}(X_t)\dd L_t,\] can loosely be thought of as a {\it stochastic feedback law} with
the first term is deterministic component. The second term adds noise, and the third term is a correction process that keeps the system in the domain $\Omega$. 
If
\[
    \Pbb(X_t\in\dd x)=\rho_t(x)\dd x.
\]
That is $\rho_t(x)$ is the law of $X_t$. Then $\rho_t$ satisfies
\begin{equation}
    \partial_t \rho_t
    =
    \Delta \rho_t
    +
    \nabla\cdot\bigl(\rho_t\nabla U\bigr)
    =
    \nabla\cdot\bigl(\nabla \rho_t+\rho_t\nabla U\bigr).
    \label{eq:classical-fp}
\end{equation}
Here, the \emph{Laplacian} is the divergence of the gradient:
$\Delta\phi=\nabla\cdot(\nabla\phi)
=\sum_{i=1}^{n}\partial_{x_i}^2\phi$.
Thus, $\Delta\rho_t$ describes diffusion, while
$\nabla\cdot(\rho_t\nabla U)$ describes transport under
the drift $-\nabla U$.
Due to invariance of the domain $\Omega$, we need to supplement the equation with a boundary condition for the PDE,
\begin{equation}
    n_{\mathrm{out}}\cdot
    \bigl(\nabla \rho_t+\rho_t\nabla U\bigr)
    =
    0
    \qquad\text{on }\partial\Omega,
    \label{eq:zero-flux}
\end{equation}
where $n_{\mathrm{out}}=-n_{\mathrm{in}}$. This condition guarantees conservation of mass: $\int_\Omega \rho_t(x) dx= 1$ for all $t \geq 0$ .

Given our interest in the long-term behavior of the FPK, we can first try to identify the steady states of the FPK.

One can check by direct substitution that the steady state, or
{\it stationary density}, is
\begin{equation}
	\rho_\infty(x)
	=
	\frac{1}{Z}e^{-U(x)},
	\qquad
	Z
	=
	\int_\Omega e^{-U(x)}\dd x.
	\label{eq:gibbs-density}
\end{equation}
This already gives us some information about how to solve the inverse
problem when $\rho_\infty$ is prescribed. If $\rho_\infty$ is a smooth,
strictly positive probability density, it may be made the target
stationary density by choosing
\begin{equation}
	U(x)
	=
	-\log\rho_\infty(x)
	+
	\text{constant}.
\end{equation}

For now, we will take it on faith that this stationary density is
unique and ask whether it is also {\it asymptotically stable}. That
is, does every solution of the Fokker--Planck equation converge to
$\rho_\infty$ as $t\rightarrow\infty$?

When $U\equiv0$, equation \eqref{eq:classical-fp} becomes
\begin{equation}
	\partial_t\rho_t
	=
	\Delta\rho_t
\end{equation}
with {\it Neumann boundary conditions}. On a connected, bounded
domain, the stationary density is
\[
\rho_\infty
=
\frac{1}{|\Omega|}.
\]
The long-time behavior of Fokker--Planck equations is often determined
by certain functional inequalities \cite{bakry2014analysis}. In this regard, the following
{\it Poincar\'e inequality} is known to hold for a large class of
bounded domains $\Omega$:
\begin{equation}
	\boxed{
	\|f-f_\Omega\|_{L^2(\Omega)}^2
	\leq
	C_{\mathrm P}
	\|\nabla f\|_{L^2(\Omega)}^2 },
\end{equation}
with 	\[
\|\phi\|_{L^2(\Omega)},
:=
\left(
\int_\Omega |\phi(x)|^2\dd x
\right)^{1/2}
 \quad f_\Omega
:=
\frac{1}{|\Omega|}
\int_\Omega f(x)\dd x \qquad.\]
From this inequality, one can derive
\begin{equation}
	\|\rho_t-\rho_\infty\|_{L^2(\Omega)}
	\leq
	e^{-t/C_{\mathrm P}}
	\|\rho_0-\rho_\infty\|_{L^2(\Omega)}.
	\label{eq:heat-convergence}
\end{equation}

To see this, define the Lyapunov functional on the space of densities
by
\[
\mathcal F(\rho)
=
\|\rho-\rho_\infty\|_{L^2(\Omega)}^2,
\]
is the standard $L^2$ norm. Differentiating the Lyapunov functional
along a solution and integrating by parts gives
\[
\frac{\dd}{\dd t}\mathcal F(\rho_t)
=
-2\int_\Omega|\nabla\rho_t|^2\dd x
\leq
-\frac{2}{C_{\mathrm P}}\mathcal F(\rho_t).
\]
We recall that if a nonnegative continuously differentiable function
$E$ satisfies
\[
\frac{\dd}{\dd t}E(t)
\leq
-cE(t)
\]
for every $t\geq0$, then
\[
E(t)\leq e^{-ct}E(0).
\]
Applying this result to $\mathcal F(\rho_t)$ and then taking square
roots gives \eqref{eq:heat-convergence}. This additionally guarantees that we have only one stationary density.

One can extend this idea to more general potentials. Define
\begin{equation}
	f_t
	=
	\frac{\rho_t}{\rho_\infty}.
\end{equation}
Using $\nabla\rho_\infty=-\rho_\infty\nabla U$, the Fokker--Planck
equation becomes
\begin{equation}
	\partial_t\rho_t
	=
	\nabla\cdot
	\left(
	\rho_\infty\nabla f_t
	\right).
	\label{eq:weighted-fp}
\end{equation}

For a function $\phi$, define the weighted $L^2$ norms by
\[
\|\phi\|_{L^2(\rho_\infty)}
:=
\left(
\int_\Omega
|\phi(x)|^2\rho_\infty(x)\dd x
\right)^{1/2}
\]
and
\[
\|\phi\|_{L^2(\rho_\infty^{-1})}
:=
\left(
\int_\Omega
\frac{|\phi(x)|^2}{\rho_\infty(x)}\dd x
\right)^{1/2}.
\]
We also define the weighted mean
\[
\langle f\rangle_{\rho_\infty}
:=
\int_\Omega f(x)\rho_\infty(x)\dd x.
\]

Consider the weighted Lyapunov functional
\begin{align}
	\mathcal F(\rho_t)
	&=
	\|f_t-1\|_{L^2(\rho_\infty)}^2 \\
	&=
	\|\rho_t-\rho_\infty\|_{L^2(\rho_\infty^{-1})}^2.
\end{align}
Differentiating along a solution and integrating by parts gives
\begin{equation}
	\frac{\dd}{\dd t}\mathcal F(\rho_t)
	=
	-2
	\int_\Omega
	\rho_\infty|\nabla f_t|^2\dd x.
	\label{eq:weighted-dissipation}
\end{equation}

Suppose that the weighted Poincar\'e inequality
\begin{equation}
	\|f-\langle f\rangle_{\rho_\infty}\|_{L^2(\rho_\infty)}^2
	\leq
	C_{\mathrm{WP}}
	\int_\Omega
	\rho_\infty|\nabla f|^2\dd x
	\label{eq:weighted-poincare}
\end{equation}
holds. Since both $\rho_t$ and $\rho_\infty$ have total mass one,
\[
\langle f_t\rangle_{\rho_\infty}
=
\int_\Omega
\frac{\rho_t}{\rho_\infty}\rho_\infty\dd x
=
1.
\]
Consequently,
\[
\frac{\dd}{\dd t}\mathcal F(\rho_t)
\leq
-\frac{2}{C_{\mathrm{WP}}}
\mathcal F(\rho_t).
\]
Gr\"onwall's inequality therefore gives
\begin{equation}
	\|\rho_t-\rho_\infty\|_{L^2(\rho_\infty^{-1})}
	\leq
	e^{-t/C_{\mathrm{WP}}}
	\|\rho_0-\rho_\infty\|_{L^2(\rho_\infty^{-1})}.
\end{equation}

On a bounded domain, the weighted Poincar\'e inequality follows from
the ordinary Poincar\'e inequality whenever $\rho_\infty$ is bounded
above and bounded away from zero. Suppose that
\[
0<m
\leq
\rho_\infty(x)
\leq
M<\infty
\]
for all $x\in\Omega$.  The constant
$\langle f\rangle_{\rho_\infty}$ minimizes
\[
c
\longmapsto
\int_\Omega |f-c|^2\rho_\infty\dd x.
\]
Hence, we may use the unweighted mean $f_\Omega$ as a competitor. It follows
that
\begin{align}
	\|f-\langle f\rangle_{\rho_\infty}\|_{L^2(\rho_\infty)}^2
	&\leq
	\int_\Omega
	|f-f_\Omega|^2\rho_\infty\dd x \\
	&\leq
	M\int_\Omega |f-f_\Omega|^2\dd x \\
	&\leq
	MC_{\mathrm P}
	\int_\Omega|\nabla f|^2\dd x \\
	&\leq
	\frac{M}{m}C_{\mathrm P}
	\int_\Omega
	\rho_\infty|\nabla f|^2\dd x.
\end{align}
Thus, one may take
\begin{equation}
	C_{\mathrm{WP}}
	\leq
	\frac{M}{m}C_{\mathrm P}.
\end{equation}

For the density $
\rho_\infty(x)
=
\frac{1}{Z}e^{-U(x)},$
we have
\[
\frac{M}{m}
=
\exp\left(
\sup_{x\in\Omega}U(x)
-
\inf_{x\in\Omega}U(x)
\right).
\]
Therefore,
\begin{equation}
	C_{\mathrm{WP}}
	\leq
	C_{\mathrm P}
	\exp\left(
	\sup_{x\in\Omega}U(x)
	-
	\inf_{x\in\Omega}U(x)
	\right).
\end{equation}
\paragraph{Numerical Example}

As a simple two-dimensional example, consider the target density on
$\Omega=[0,1]^2$ given by
\begin{equation}
	\rho_\infty(x,y)
	\propto
	\sum_{k=1}^{4}
	\exp\!\left(
	-\frac{
		(x-\mu_{x,k})^2
		+
		(y-\mu_{y,k})^2
	}{2\sigma^2}
	\right),
\end{equation}
where
\[
(\mu_{x,k},\mu_{y,k})
\in
\left\{
(0.2,0.8),
(0.8,0.2),
(0.2,0.2),
(0.8,0.8)
\right\}.
\]
Choosing
\[
U(x,y)=-\log\rho_\infty(x,y)
\]
makes $\rho_\infty$ the invariant density of the reflected
diffusion.

\begin{figure}[htbp]
	\centering
	\includegraphics[width=\textwidth]
	{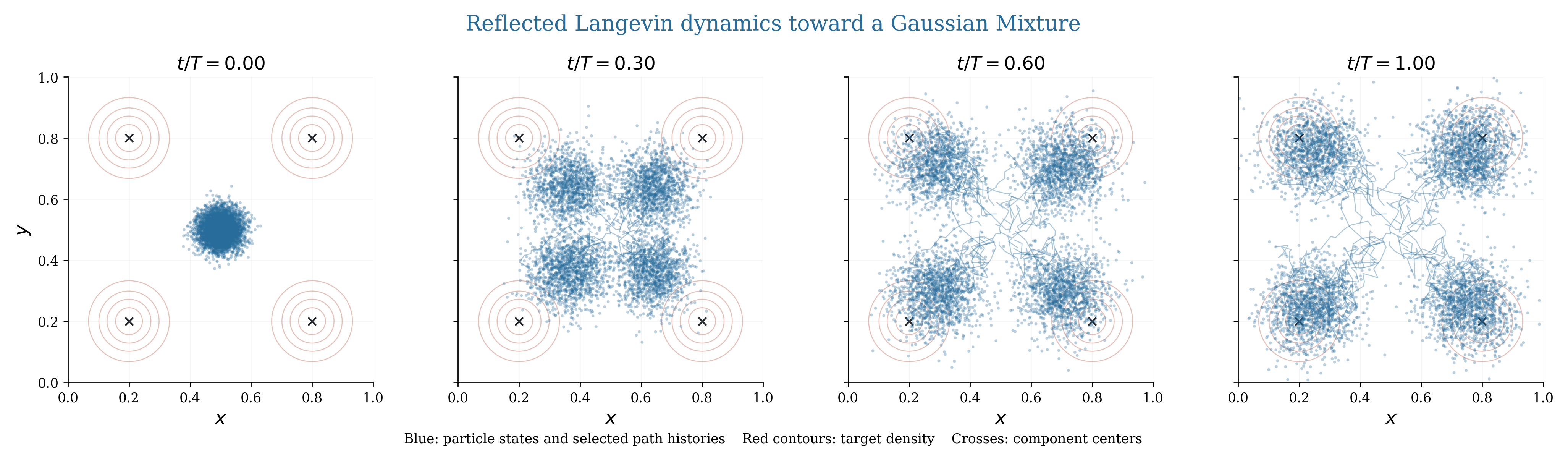}
	\caption{Evolution of reflected diffusion particles toward a
		four-component Gaussian-mixture target on $[0,1]^2$.
		Blue points represent the particle states, red contours represent
		the target density $\rho_\infty$, and the crosses mark the centers of the four
		component. The particles are initialized near
		$(0.5,0.5)$ and separate into the four target
		modes over time.}
	\label{fig:reflected-gmm}
\end{figure}

For the numerical experiment, particles are initialized in a small
neighborhood of the center $(0.5,0.5)$ and propagated using the
Euler--Maruyama method. Figure~\ref{fig:reflected-gmm} shows how the
initial cloud separates into four groups and approaches the modes of
the prescribed Gaussian mixture. The random perturbation spreads the
particles, while the drift $\nabla\log\rho_\infty=-\nabla U$ directs
them toward regions of higher target density.

\section{Nonholonomic Fokker-Planck Equations}
\label{sec:nonholonomic-fp}

The classical Langevin model introduced in the last section can be used to control the distribution of a single integrator $\dot{x} = u(t)$. For robotic systems with nonholonomic or underactuated constraints, diffusion must instead respect the admissible control directions $g_i$ introduced in \eqref{eq:global-control-affine}. In this section, we consider diffusion processes associated with a driftless system,
\begin{equation}
	\dot{x}  = \sum_{i=1}^m g_i(x)u_i,
	\end{equation}
where $g_1,\dots,g_m$ be smooth vector fields and the drift $f_0 \equiv 0$.

Now, consider the reflected Stratonovich SDE associated with a driftless system
\begin{equation}
    \dd X_t
    =
    \sum_{i=1}^{m}
    v_i(X_t)g_i(X_t)\dd t
    +
    \sqrt{2}
    \sum_{i=1}^{m}
    g_i(X_t)\circ\dd W_t^i
    +
    n_{\mathrm{in}}(X_t)\dd L_t.
    \label{eq:nonholonomic-sde}
\end{equation}
Here, $v_i(x)$ can be thought of as deterministic component of the control law and noise acts along only the control vector-fields $g_i$.  We want to describe the Fokker-Planck equation associated with this process. 

Corresponding to these vector-fields,  we define the operator
\begin{equation}
	\mathcal Y_i\varphi
	:=
	g_i\cdot\nabla\varphi.
\end{equation}
The {\it formal adjoint} acting on densities is
\begin{equation}
	\mathcal Y_i^*q
	:=
	-\nabla\cdot(g_i q).
	\label{eq:formal-adjoint}
\end{equation}
Given these operators, in the interior of $\Omega$ the  Fokker--Planck equation is
\begin{equation}
    \partial_t \rho_t
    =
    \sum_{i=1}^{m}
    \mathcal Y_i^*(v_i \rho_t)
    +
    \sum_{i=1}^{m}
    (\mathcal Y_i^*)^2\rho_t.
    \label{eq:nonholonomic-fp}
\end{equation}
One can see that we recover our original FPK if the vector-fields $g_i$ act in coordinate directions. For the upcoming derivation, we would like to express this equation into a more appealing form.

Toward this end let
\begin{equation}
    d_i
    :=
    \nabla\cdot g_i.
\end{equation}
From \eqref{eq:formal-adjoint},
\begin{equation}
    \mathcal Y_i^*
    =
    -\mathcal Y_i-d_i.
    \label{eq:adjoint-decomposition}
\end{equation}

For the sign convention in \eqref{eq:nonholonomic-fp}, define
\begin{equation}
    U_i
    :=
    v_i-d_i.
    \label{eq:Ui-correction}
\end{equation}
Then
\begin{equation}
    (\mathcal Y_i^*)^2q
    +
    \mathcal Y_i^*(v_i q)
    =
    -\mathcal Y_i^*\mathcal Y_i q
    +
    \mathcal Y_i^*(U_i q).
    \label{eq:operator-correction}
\end{equation}
Equivalently, $v_i=d_i+U_i$.
 To see this, using $\mathcal Y_i^*=-\mathcal Y_i-d_i$,
\[
    (\mathcal Y_i^*)^2q
    =
    \mathcal Y_i^*\bigl(-\mathcal Y_i q-d_i q\bigr)
    =
    -\mathcal Y_i^*\mathcal Y_i q
    -
    \mathcal Y_i^*(d_i q).
\]
Therefore,
\begin{align*}
    (\mathcal Y_i^*)^2q
    +
    \mathcal Y_i^*(v_i q)
    &=
    -\mathcal Y_i^*\mathcal Y_i q
    +
    \mathcal Y_i^*\bigl((v_i-d_i)q\bigr)\\
    &=
    -\mathcal Y_i^*\mathcal Y_i q
    +
    \mathcal Y_i^*(U_i q).
\end{align*}
Therefore, we have our equality. 

Thus the Fokker--Planck equation \eqref{eq:nonholonomic-fp} can be written as
\begin{equation}
    \partial_t \rho_t
    =
    -\sum_{i=1}^{m}
    \mathcal Y_i^*\mathcal Y_i \rho_t
    +
    \sum_{i=1}^{m}
    \mathcal Y_i^*(U_i \rho_t).
    \label{eq:transformed-fp}
\end{equation}
To make a smooth, strictly positive density $\rho_\infty$ invariant,
we choose
\begin{equation}
	U_i
	=
	\mathcal Y_i\log\rho_\infty.
	\label{eq:Ui-target}
\end{equation}
Indeed,
\[
U_i\rho_\infty
=
\rho_\infty\mathcal Y_i\log\rho_\infty
=
\mathcal Y_i\rho_\infty,
\]
so each term in \eqref{eq:transformed-fp} vanishes when
$\rho=\rho_\infty$. In terms of the original control law $v_i(x)$,
this gives
\begin{equation}
	\boxed{
		v_i
		=
		\nabla\cdot g_i
		+
		\mathcal Y_i\log\rho_\infty
	}.
	\label{eq:vi-target}
\end{equation}
If $\rho_\infty\propto e^{-V}$, then
\[
v_i
=
\nabla\cdot g_i
-
\mathcal Y_iV.
\]

Once again, we would like to determine whether this control law
asymptotically stabilizes the system toward $\rho_\infty$. Define
\[
f_t
=
\frac{\rho_t}{\rho_\infty}.
\]
Then equation \eqref{eq:transformed-fp} becomes
\begin{equation}
	\partial_t\rho_t
	=
	-\sum_{i=1}^{m}
	\mathcal Y_i^*
	\left(
	\rho_\infty\mathcal Y_i f_t
	\right).
	\label{eq:weighted-horizontal-fp}
\end{equation}
Consider the Lyapunov functional
\begin{equation}
	\mathcal F(\rho_t)
	=
	\int_\Omega
	(f_t-1)^2\rho_\infty\,\dd x
	=
	\|\rho_t-\rho_\infty\|_{L^2(\rho_\infty^{-1})}^2.
\end{equation}
Differentiating along solutions of
\eqref{eq:weighted-horizontal-fp} and integrating by parts gives
\begin{equation}
	\frac{\dd}{\dd t}\mathcal F(\rho_t)
	=
	-2
	\sum_{i=1}^{m}
	\int_\Omega
	\rho_\infty
	|\mathcal Y_i f_t|^2\,\dd x.
	\label{eq:horizontal-dissipation}
\end{equation}

Suppose that the following \emph{horizontal Poincar\'e inequality}
holds:
\begin{equation}
	\boxed{
		\left\|
		f-\langle f\rangle_{\rho_\infty}
		\right\|_{L^2(\rho_\infty)}^2
		\leq
		C_{\mathrm H}
		\sum_{i=1}^{m}
		\|\mathcal Y_i f\|_{L^2(\rho_\infty)}^2
	}.
	\label{eq:horizontal-poincare}
\end{equation}
Since $\rho_t$ and $\rho_\infty$ are probability densities,
\[
\langle f_t\rangle_{\rho_\infty}
=
\int_\Omega f_t\rho_\infty\,\dd x
=
\int_\Omega\rho_t\,\dd x
=
1.
\]
Consequently, \eqref{eq:horizontal-poincare} and
\eqref{eq:horizontal-dissipation} imply
\[
\frac{\dd}{\dd t}\mathcal F(\rho_t)
\leq
-\frac{2}{C_{\mathrm H}}\mathcal F(\rho_t).
\]
Gr\"onwall's inequality therefore gives
\begin{equation}
	\|\rho_t-\rho_\infty\|_{L^2(\rho_\infty^{-1})}
	\leq
	e^{-t/C_{\mathrm H}}
	\|\rho_0-\rho_\infty\|_{L^2(\rho_\infty^{-1})}.
	\label{eq:horizontal-convergence}
\end{equation}

\paragraph{When does a horizontal Poincar\'e inequality hold?}

The usual Poincar\'e inequality reflects the connectivity of the
domain $\Omega$. Similarly, the horizontal Poincar\'e inequality is
related to connectivity in the geometry generated by the control
vector fields $g_1,\ldots,g_m$. That is, {\it is the system globally controllable?} A standard sufficient condition is
that the vector fields satisfy the \emph{bracket-generating
	condition}
\begin{equation}
	\operatorname{Lie}
	\{g_1,\ldots,g_m\}(x)
	=
 \mathbb{R}^d,
	\qquad
	x\in\Omega.
	\label{eq:bracket-generating}
\end{equation}
 Here, $\operatorname{Lie}\{g_1,\ldots,g_m\}(x)$ denotes the linear
span, at $x$, of the vector fields and all of their iterated Lie
brackets, such as
\[
[g_i,g_j],
\qquad
[g_i,[g_j,g_k]],
\qquad \ldots.
\]

For the associated driftless control system
\begin{equation}
	\dot{x}(t)
	=
	\sum_{i=1}^{m}
	u_i(t)g_i\bigl(x(t)\bigr),
\end{equation}
the bracket-generating condition implies, through the
Chow--Rashevskii theorem \cite{agrachev2013control}, that any two points can be connected by a
curve whose velocity lies in the span of the control vector fields,
provided the domain is horizontally connected and the connecting
curve is allowed to remain in $\Omega$. Thus, controllability in the
horizontal geometry plays the role of ordinary Euclidean
connectivity.

This connection can also be seen directly from the horizontal
derivatives. If
\[
\mathcal Y_i f=0,
\qquad
i=1,\ldots,m,
\]
then
\[
[\mathcal Y_i,\mathcal Y_j]f
=
\mathcal Y_i\mathcal Y_jf
-
\mathcal Y_j\mathcal Y_if
=
0,
\]
and the same conclusion holds for every iterated commutator. Under
the bracket-generating condition, these directions span the entire
tangent space. 

Under suitable regularity assumptions on the vector fields and on
the bounded, horizontally connected domain $\Omega$, this
qualitative connectivity property yields a quantitative horizontal
Poincar\'e inequality which is related to embedding results established in \cite{garofalo1996isoperimetric} and one can prove the density stability result \cite{elamvazhuthi2023density}.

\paragraph{Numerical experiment.}

We repeat the preceding Gaussian-mixture experiment for the unicycle
system on
\[
\Omega
=
[0,1]^2\times\mathbb S^1.
\]
Writing the state as $X=(x,y,\theta)$, the admissible vector fields are
\[
g_1(X)
=
\begin{pmatrix}
	\cos\theta\\
	\sin\theta\\
	0
\end{pmatrix},
\qquad
g_2(X)
=
\begin{pmatrix}
	0\\
	0\\
	1
\end{pmatrix}.
\]
The first control moves the unicycle in its current heading direction,
while the second changes its orientation. Both vector fields are
divergence free:
\[
\nabla\cdot g_1
=
\nabla\cdot g_2
=
0.
\]

Let $\bar\rho_\infty(x,y)$ denote the four-component Gaussian-mixture
density from the preceding example. We prescribe the joint target
density
\begin{equation}
	\rho_\infty(x,y,\theta)
	=
	\frac{1}{2\pi}\bar\rho_\infty(x,y),
	\label{eq:unicycle-target-density}
\end{equation}
so that the desired position distribution is the Gaussian mixture,
while the desired heading distribution is uniform on $\mathbb S^1$.

The corresponding horizontal derivatives are
\[
\mathcal Y_1
=
\cos\theta\,\partial_x
+
\sin\theta\,\partial_y,
\qquad
\mathcal Y_2
=
\partial_\theta.
\]
Since $\rho_\infty$ is independent of $\theta$, the feedback law
\eqref{eq:vi-target} becomes
\begin{align}
	v_1(x,y,\theta)
	&=
	\mathcal Y_1\log\rho_\infty(x,y,\theta) \\
	&=
	\cos\theta\,\partial_x\log\bar\rho_\infty(x,y)
	+
	\sin\theta\,\partial_y\log\bar\rho_\infty(x,y),
	\label{eq:unicycle-v1}\\
	v_2(x,y,\theta)
	&=
	\mathcal Y_2\log\rho_\infty(x,y,\theta)
	=
	0.
\end{align}
The resulting diffusive dynamics are
\begin{align}
	\dd x_t
	&=
	\cos\theta_t
	\left(
	v_1(X_t)\,\dd t
	+
	\sqrt{2}\circ\dd W_t^1
	\right),\\
	\dd y_t
	&=
	\sin\theta_t
	\left(
	v_1(X_t)\,\dd t
	+
	\sqrt{2}\circ\dd W_t^1
	\right),\\
	\dd\theta_t
	&=
	\sqrt{2}\circ\dd W_t^2,
	\label{eq:unicycle-langevin}
\end{align}
together with reflecting boundary conditions in the position
coordinates.

Although the noise acts only along the current heading and in the
angular direction, the system explores both position coordinates.
Indeed,
\[
[g_2,g_1]
=
\begin{pmatrix}
	-\sin\theta\\
	\cos\theta\\
	0
\end{pmatrix},
\]
and therefore
\[
\operatorname{span}
\left\{
g_1,g_2,[g_2,g_1]
\right\}
=
\mathbb R^3
\]
at every state. Hence, the unicycle vector fields satisfy the
bracket-generating condition.

\begin{figure}[htbp]
	\centering
	\includegraphics[width=\textwidth]
	{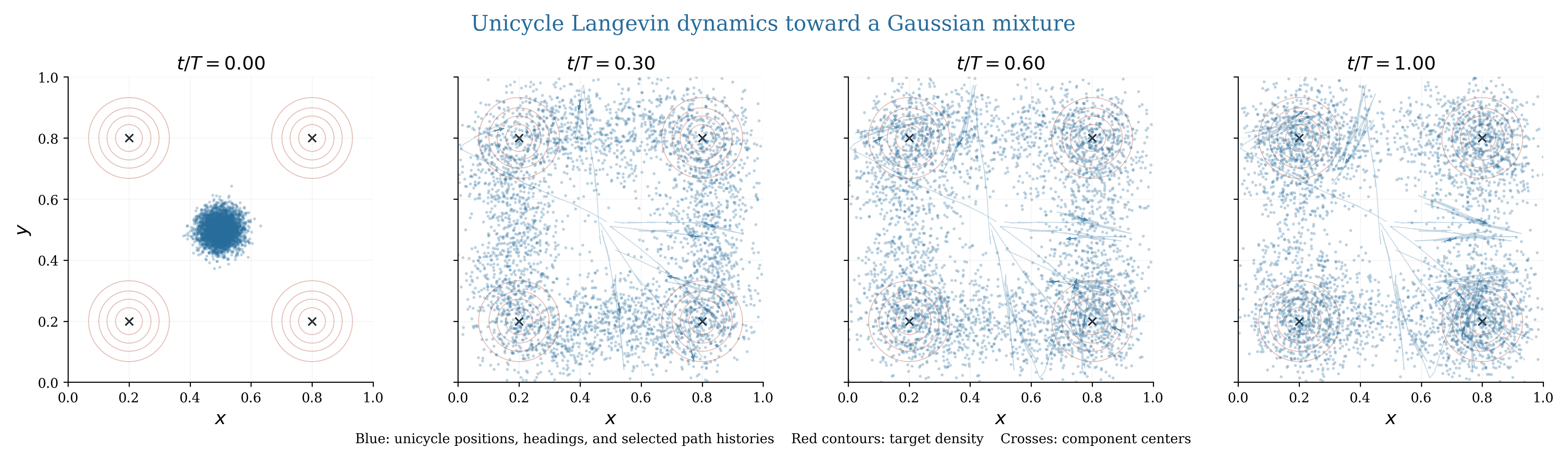}
	\caption{Stochastic gradient flow for the unicycle
		model. Blue points represent the unicycle positions, and the
		arrows indicate the headings of a subset of the particles. Red
		contours show the desired four-component Gaussian-mixture density,
		while the crosses mark its component centers. The positions are
		initialized near $(0.5,0.5)$ with uniformly distributed headings.}
	\label{fig:unicycle-gmm-langevin}
\end{figure}

For the numerical experiment, the positions are initialized in a
small neighborhood of $(0.5,0.5)$, while the headings are sampled
uniformly from $\mathbb S^1$. Figure~\ref{fig:unicycle-gmm-langevin}
shows the evolution of the resulting particle system. The particles
initially move only along their instantaneous heading directions, but
angular diffusion continually changes these directions. Through this
bracket-generating mechanism, the position distribution separates
into four groups and approaches the prescribed Gaussian-mixture
target.

The traces of the sample show the difference between the geometries of the single integrator of the previous section and the unicycle model. In the unicycle system, each instantaneous displacement is
constrained to the current heading. Transverse motion appears only
through changes in orientation and the resulting Lie-bracket
directions. Accordingly, the unicycle paths exhibit a visibly more
constrained path, even though their terminal
position distribution approaches the same four-mode target. In the fully actuated case as seen in Figure \ref{fig:reflected-gmm}, the samples disperse radially, showing motion in both the x- and y-direction proceeds uniformly.

\section{Diffusions for Control Systems with Drift}
\label{sec:fp-with-drift}

The previous section considered diffusions along the control directions. We now restore the drift $f_0$. In this case, we will not be able to arbitarily control the stationary distribution. But nevertheless, there are some important properties that we can identify which will be useful in the next chapter on denoising diffusion. 

Consider the Stratonovich SDE
\begin{equation}
    \dd X_t
    =
    f_0(X_t)\dd t
    +
    \sqrt{2}
    \sum_{i=1}^{m}
    g_i(X_t)\circ\dd W_t^i.
    \label{eq:drifted-sde}
\end{equation}
With
\[
    \mathcal Y_0\varphi
    =
    f_0\cdot\nabla\varphi,
\]
the density satisfies the nonholonomic FPK,
\begin{equation}
    \partial_t \rho_t
    =
    \mathcal Y_0^*\rho_t
    +
    \sum_{i=1}^{m}
    (\mathcal Y_i^*)^2\rho_t.
    \label{eq:drifted-fp}
\end{equation}
Note that here we are not considering a reflected diffusion process. The presence of drift can make things very challenging to answer questions of the system behavior in the setting when the process is confined to a bounded domain $\Omega$. 

We will consider two situations. One, Linear time invariant systems (LTIs), and then even more briefly, the general case.

\subsection{Linear time-invariant systems.}
We consider a linear drift $f_0(x)=Ax$, where
$A\in\R^{n\times n}$, and constant control fields
$g_i(x)=b_i$, where $b_i\in\R^n$. Collecting these fields
as the columns of
\[
B=\begin{bmatrix}b_1 & \cdots & b_m\end{bmatrix}
\in\R^{n\times m},
\]
we can express the SDE in the familiar form
\begin{equation}
    \dd X_t
    =
    AX_t\dd t
    +
    \sqrt{2}B\dd W_t,
\end{equation}
the and the corresponding Fokker--Planck equation is
\begin{equation}
    \partial_t \rho_t
    =
    -\nabla\cdot(Ax\,\rho_t)
    +
    \nabla\cdot
    \left(
        BB^{\mathsf T}\nabla \rho_t
    \right).
\end{equation}
An appealing aspect of the LTI case is that one can almost write an explicit solution for the Fokker-Planck Equation. For this, we need a controllability assumption that  
Define
\begin{equation}
    \Sigma_t
    =
    2
    \int_0^t
    e^{sA}
    BB^{\mathsf T}
    e^{sA^{\mathsf T}}
    \dd s.
\end{equation}
When the LTI system is controllable, $\Sigma_t$ is invertible. The {\it transition kernel} corresponding to this SDE is then given by
\begin{equation}
    K(t,x,y)
    =
    \frac{
        \exp\!\left(
            -\frac{1}{2}
            (x-e^{tA}y)^{\mathsf T}
            \Sigma_t^{-1}
            (x-e^{tA}y)
        \right)
    }{
        (2\pi)^{n/2}
        (\det\Sigma_t)^{1/2}
    },
\end{equation}
The transition kernel denotes the law of the process if the process was initialized at $y$. One can directly substitute this into the Fokker-Planck equation and check that this indeed satisfies the equality. For general initial probability distribution $\rho_0$, the evolution of the law can be expressed as an integral,
\begin{equation}
    \rho_t(x)
    =
    \int_{\R^n}
    K(t,x,y)\rho_0(y)\dd y.
\end{equation}
From this expression, one can even conclude the long-term behavior of the process if $A$ is Hurwitz. In this case, the following infinite-time Grammian is well-defined.
\begin{equation}
    \Sigma_\infty
    =
    2
    \int_0^\infty
    e^{sA}
    BB^{\mathsf T}
    e^{sA^{\mathsf T}}
    \dd s
\end{equation}
and one can formally see that  
solves
\[   \lim_{t \rightarrow \infty} K(t,x,y) =     \rho_\infty(x)
=
\frac{
	\exp\!\left(
	-\frac{1}{2}
	x^{\mathsf T}
	\Sigma_\infty^{-1}x
	\right)
}{
	(2\pi)^{n/2}
	(\det\Sigma_\infty)^{1/2}
}.  \]

This implies, that irrespective of the initial condition $y$, the distribution or the law of the process converges to $\rho_{\infty}$. Therefore, more generally,

\begin{equation}
\lim_{ t \rightarrow \infty} 	\rho_t(x)
	=
	\int_{\R^n}
	\rho_\infty(x)\rho_0(y)\dd y = 	\rho_\infty(x)
\end{equation}

\begin{equation}
    A\Sigma_\infty
    +
    \Sigma_\infty A^{\mathsf T}
    +
    2BB^{\mathsf T}
    =
    0.
\end{equation}

\paragraph{Numerical Example.}
We compare a fully actuated damped single integrator,
\[
\dd X_t=-X_t\dd t+\sqrt{2}\dd W_t,
\qquad X_t\in\R^2,
\]
with a damped oscillator driven by noise only in its second
coordinate,
\[
\dd X_t^1=X_t^2\dd t,
\qquad
\dd X_t^2=(-2X_t^1-X_t^2)\dd t+\sqrt{2}\dd W_t.
\]
Both drift matrices are Hurwitz, and both systems are controllable.
Their stationary laws are centered Gaussians with covariances
\[
\Sigma_\infty^{\mathrm{int}}=I_2,
\qquad
\Sigma_\infty^{\mathrm{osc}}
=
\begin{pmatrix}
	1/2 & 0\\
	0 & 1
\end{pmatrix},
\]
respectively.

\begin{figure}[htbp]
	\centering
	\includegraphics[width=\textwidth]
	{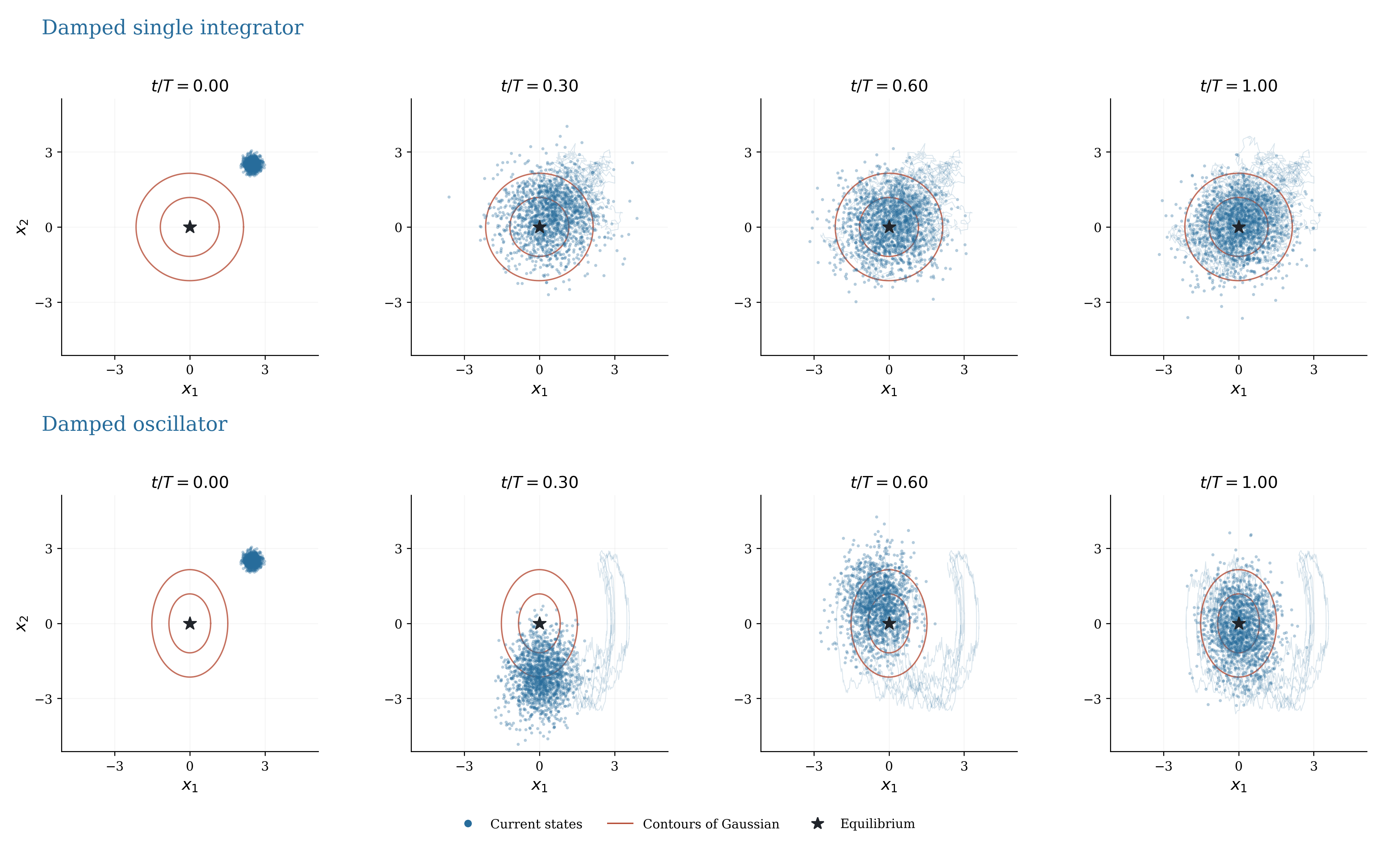}
	\caption{Convergence to Gaussian stationary densities.
		Top: fully actuated damped single integrator.
		Bottom: damped oscillator with noise in its second coordinate.
		Blue points show current states, and faint blue curves trace
		selected particles. The star marks the equilibrium.}
	\label{fig:lti-gaussian-convergence}
\end{figure}

Figure~\ref{fig:lti-gaussian-convergence} shows the evolution from
the same initial Gaussian cloud, with mean $(2.5,2.5)^{\mathsf T}$
and covariance $0.15^2 I_2$. For the damped integrator, the cloud
moves toward the origin and becomes isotropic. For the oscillator,
the cloud rotates and deforms as noise propagates from the second
coordinate to the first, eventually approaching an anisotropic
Gaussian. Thus, both processes forget their initial distribution,
but their dynamics determine different stationary laws.

\subsection{General Nonlinear Systems with Drift}

For nonlinear systems, explicit transition kernels are generally
unavailable, and long-term convergence typically requires additional
estimates, such as the Lyapunov and Poincar\'e inequalities discussed
earlier. We therefore return to the support properties considered
in Chapter~\ref{chap:not-sample-reachable}. These connect the
diffusion to the reachable sets of the corresponding control system.

Consider \eqref{eq:drifted-sde} with $X_0=x_0$. Recall that
$\mathcal R_T^{L^2}(x_0)$ denotes the set of states reachable at
time $T$ under finite-energy controls in
\eqref{eq:global-control-affine}.

	The
	Stroock--Varadhan support theorem \cite{stroock1972support} gives
	\begin{equation}
		\operatorname{supp}\bigl(\mathcal L(X_T)\bigr)
		=
		\overline{\mathcal R_T^{L^2}(x_0)}.
		\label{eq:endpoint-support}
	\end{equation}

Thus, every neighborhood of a reachable state
has positive probability. In particular, if the reachable set is
dense in $\R^n$ at time $T$, the transition law has full support.

One can additionally consider when the density of the process is smooth. This requires a bracket-generating condition, similar to what we saw tin the previous sections. More precisely, define
\[
\mathcal V_0=\{g_1,\dots,g_m\},
\qquad
\mathcal V_{k+1}
=
\mathcal V_k
\cup
\bigl\{
[h,V]:
h\in\{f_0,g_1,\dots,g_m\},\;
V\in\mathcal V_k
\bigr\}.
\]
The \emph{parabolic H\"ormander condition} is
\begin{equation}
	\boxed{
		\operatorname{span}
		\left\{
		V(x):V\in\bigcup_{k=0}^{\infty}\mathcal V_k
		\right\}
		=
		\R^n,
		\qquad x\in\R^n.
	}
	\label{eq:parabolic-hormander}
\end{equation}
This condition implies  what is known as {\it strong accessibility} of the
associated control system. That is, the reachable set at each sufficiently
small fixed positive time has nonempty interior. It does not,
however, imply global controllability.

Without this condition, a density need not exist. For example,
consider the system $\dot{x}_1=u$, $\dot{x}_2=0$.
Replacing the control by white noise gives
$\dd X_t^1=\sqrt{2}\dd W_t$ and $\dd X_t^2=0$.
Starting from a fixed initial condition, the law remains supported
on a line and therefore has no density with respect to
two-dimensional Lebesgue measure.

Even when a smooth density exists, full support does not guarantee
its strict positivity at every point, a property useful in diffusion models. Positivity at a particular
state can be established through a regularity condition on the
control-to-endpoint map. Recall that this map is
$E_T(u)=x_T^u$, with the initial condition $x_0$ fixed.
	Under above parabolic H\"ormander condition one has the 
	{\it Ben Arous--L\'eandre positivity criterion} 
	\cite{arous1991decroissance} that if there exists
	$u\in L^2(0,T;\R^m)$ such that
	\[
	E_T(u)=y,
	\qquad
	DE_T(u):L^2(0,T;\R^m)\rightarrow\R^n
	\quad\text{is surjective}.
	\]
	Then $\rho_t(y)>0$ for the process initialized at $x_0$.

Here, surjectivity of the derivative of the endpoint map $DE_T(u)$ is checking if 
small perturbations of the chosen control can change the terminal
state in every direction. Equivalently, the control system
linearized along $x_t^u$ is controllable over $[0,T]$.

The support theorem therefore concerns approximate reachability,
while this positivity criterion requires a control reaching the
specified state for which the linearized endpoint map is onto.

\section{Brockett's Condition and Stochastic Distribution Stabilization}

The density-shaping results developed throughout this chapter provide
a different perspective on the feedback stabilization problem in nonlinear control. Rather than requiring every trajectory to
converge to a desired equilibrium, we have sought to stabilize the
\emph{probability distribution} of the state. This section revisits
that distinction in light of the topological obstructions to smooth
deterministic stabilization. The stochastic approach changes the
stabilization objective, allowing us to obtain convergence of
distributions even when smooth deterministic stabilization of a
point is impossible.

\paragraph{The classical stabilization problem}
Consider the control system
\[
\dot{x}=f(x,u)
\]
and a desired equilibrium $x_\star$. Informally, the goal is to
construct a time-independent feedback $u=k(x)$ such that
$f(x_\star,k(x_\star))=0$ and the resulting closed-loop system
\[
\dot{x}=f(x,k(x))
\]
has an asymptotically stable equilibrium at $x_\star$. This requires
two properties: trajectories starting sufficiently close to
$x_\star$ remain close for all future times, and trajectories
starting in some neighborhood of $x_\star$ converge to it:
\[
\lim_{t\to\infty}x(t)=x_\star.
\]
These are the properties of \emph{stability} and
\emph{local attractivity}, respectively.

A natural question is whether controllability guarantees the
existence of such a smooth feedback. Brockett's necessary condition
shows that the answer is no.  a system may be controllable and yet
admit no smooth, time-independent feedback that asymptotically
stabilizes the desired equilibrium. More fundamentally, smooth asymptotic stabilization of nonholonomic
systems is obstructed by topology and the following theorem was proved by Brockett in \cite{brockett1983asymptotic}.

\begin{theorem}[Brockett's Necessary Condition]
	Suppose the control system
	
	\[
	\dot{x}
	=
	f(x,u)
	\]
	
	admits a continuous time-invariant feedback law that asymptotically
	stabilizes the origin. Then the map
	
	\[
	(x,u)
	\mapsto
	f(x,u)
	\]
	
	must be locally surjective onto a neighborhood of the origin.
\end{theorem}

Many nonholonomic systems fail this condition. Consequently, even when
the underlying control system is globally controllable, there does not
exist a smooth time-invariant feedback law that asymptotically
stabilizes the origin.

We want to see the results of the previous sections from the perspective of this obstruction to stabilization. First, we start with the following ambitious goal.  Consider the single-integrator system 
\begin{equation}
    \dot{x}=u,
\end{equation}

where $x\in\R^n$.

One of the simpler (but a serious overkill) approaches  to achieve a stable equilibrium point for this system is to construct a smooth
potential function $
V:\R^n\rightarrow\R,$

whose minimum is the desired equilibrium, and choose the feedback law

\begin{equation}
    u=-\nabla V(x).
\end{equation}

The resulting closed-loop system

\begin{equation}
    \dot{x}
    =
    -\nabla V(x)
\end{equation}

is the \text{gradient system} or \text{gradient flow}.

Along every solution,

\begin{align}
    \frac{\dd}{\dd t}V(x(t))
    &=
    \nabla V(x)\cdot\dot{x}\\
    &=
    -\|\nabla V(x)\|^2
    \le0.
\end{align}

Thus the potential decreases monotonically along every trajectory.
Under suitable assumptions, every solution converges to the set of
critical points of $V$. Consequently, stabilization may be viewed as an
optimization problem.  Now, can we construct such a system that naturally respects the geometry of a driftless control system

\begin{equation}
    \dot{x}
    =
    \sum_{i=1}^{m}
    u_i g_i(x),
\end{equation}
The hope is that we can construct a \textbf{nonholonomic or horizontal gradient flow}. 

Towards this ambition, recall the directional derivatives $
\mathcal Y_i
=
g_i\cdot\nabla.$ A natural extension of gradient descent is obtained by choosing
\begin{equation}
    u_i
    =
    -\mathcal Y_iV.
\end{equation}
The resulting dynamics become
\begin{equation}
    \dot{x}
    =
    -
    \sum_{i=1}^{m}
    (\mathcal Y_iV)\,
    g_i(x).
\end{equation}
Defining the horizontal gradient

\[
\nabla_HV
=
\sum_{i=1}^{m}
(\mathcal Y_iV)g_i,
\]

we get the \textbf{nonholonomic or horizontal gradient flow}

\[
\boxed{
\dot{x}
=
-\nabla_HV.}
\]

Moreover,

\begin{align}
\frac{\dd}{\dd t}V(x(t))
&=
-
\sum_{i=1}^{m}
|\mathcal Y_iV|^2
\le0.
\end{align}

Thus, just as in the Euclidean case, the potential is non-increasing
along every trajectory.

Unfortunately, this analogy with classical gradient descent is not very useful. Unlike the Euclidean gradient, the horizontal gradient may vanish even though the full gradient does not. Consequently, the equilibrium set

\[
\{x:\nabla_HV(x)=0\}
\]

is generally much larger than the set of critical points of $V$.

A classical example is the nonholonomic integrator
\begin{equation}
\begin{aligned}
\dot{x} &= u_1,\\
\dot{y} &= u_2,\\
\dot{z} &= xu_2-yu_1.
\end{aligned}
\end{equation}

Although this system satisfies the Lie Algebra Rank Condition and is
globally controllable, Brockett's theorem implies that no smooth static
feedback controller can asymptotically stabilize the origin.

Take for example, the quadratic potential
\[
V(x,y,z)=\frac12(x^2+y^2+z^2),
\]
illustrating this limitation explicitly. The control fields are
\[
g_1=(1,0,-y)^{\mathsf T},
\qquad
g_2=(0,1,x)^{\mathsf T},
\]
so the horizontal derivatives are
\[
\mathcal Y_1V=x-yz,
\qquad
\mathcal Y_2V=y+xz.
\]
Consequently,
\[
|\nabla_HV|^2
=(x-yz)^2+(y+xz)^2
=(1+z^2)(x^2+y^2).
\]
Thus,
\[
\{\nabla_HV=0\}
=
\{(0,0,z):z\in\R\},
\qquad
\{\nabla V=0\}
=
\{(0,0,0)\}.
\]
The horizontal gradient vanishes along the entire $z$-axis,
although the full gradient vanishes only at the origin.

Choosing the horizontal-gradient feedback
\[
u_1=-\mathcal Y_1V=-(x-yz),
\qquad
u_2=-\mathcal Y_2V=-(y+xz)
\]
gives the closed-loop dynamics
\[
\dot{x}=-x+yz,
\qquad
\dot{y}=-y-xz,
\qquad
\dot{z}=-z(x^2+y^2).
\]
Every point on the $z$-axis is therefore an equilibrium, preventing
this feedback from asymptotically stabilizing the origin. Thus, deterministic gradient-based stabilization is fundamentally
limited for nonholonomic systems.

\paragraph{A Stochastic Alternative}

The density-shaping framework developed in the previous sections
suggest a different path around the topological obstruction. Instead of requiring every trajectory to converge to an equilibrium, we
can seek convergence of the probability density toward a prescribed
invariant distribution.

Choose a target density

\begin{equation}
\rho_\infty(x)
=
\frac{1}{Z}e^{-V(x)},
\end{equation}

where $V$ is the same potential function used in the deterministic
gradient system.

From the previous section, the transformed drift coefficient is

\begin{equation}
U_i
=
\mathcal Y_i\log \rho_\infty
=
-\mathcal Y_iV,
\end{equation}

and the corresponding feedback law is

\begin{equation}
v_i
=
\nabla\cdot g_i
-
\mathcal Y_iV.
\end{equation}

This gives the reflected Stratonovich stochastic differential equation

\begin{equation}
\dd X_t
=
\sum_{i=1}^{m}
v_i(X_t)g_i(X_t)\,\dd t
+
\sqrt{2}
\sum_{i=1}^{m}
g_i(X_t)\circ\dd W_t^i
+
n_{\mathrm{in}}(X_t)\,\dd L_t.
\end{equation}

The associated Fokker--Planck equation is precisely the one analyzed in
the previous sections. Consequently,
\[
\lim_{t\rightarrow\infty}
\rho_t
=
\rho_\infty,
\]

and, under the controllability and Poincar\'e assumptions developed
earlier, the convergence is exponential in the weighted
$L^2$ norm.

Therefore, although deterministic smooth feedback may fail to stabilize
a nonholonomic system to an equilibrium, stochastic feedback can still
stabilize the probability distribution around that equilibrium.

This interpretation provides a stochastic analogue of gradient descent.
Instead of driving every trajectory toward a minimum of the potential,
the stochastic dynamics drive the entire probability density toward the
Gibbs distribution $
\rho_\infty
\propto
e^{-V},$ which becomes increasingly concentrated near the minima of $V$.

The resulting notion of stabilization is weaker than deterministic
asymptotic stabilization, but it bypasses Brockett's obstruction while
retaining a meaningful notion of convergence.

\paragraph{Numerical example}

We illustrate the effect of stochastic perturbations for stochastic stabilization using the
unicycle control system
\[
\dot{x}=u_1\cos\theta,
\qquad
\dot{y}=u_1\sin\theta,
\qquad
\dot{\theta}=u_2,
\]
where $(x,y)\in\R^2$ is the position, $\theta\in S^1$ is the
heading, and $u_1,u_2\in\R$ are the signed translational velocity
and angular velocity, respectively. The control vector fields are
\[
g_1=(\cos\theta,\sin\theta,0)^{\mathsf T},
\qquad
g_2=(0,0,1)^{\mathsf T}.
\]

Choose the position-dependent potential
\[
V(x,y)=\frac{K}{2}(x^2+y^2).
\]
Since the control fields are divergence free, the
horizontal-gradient feedback is
\[
v_1=-\mathcal Y_1V=-K(x\cos\theta+y\sin\theta),
\qquad
v_2=-\mathcal Y_2V=0.
\] 
Keeping $K=10$ fixed, we vary the diffusion strength $D$ in
\begin{align}
	\dd x_t
	&=
	\cos\theta_t
	\left(v_1\,\dd t+\sqrt{2D}\circ\dd W_t^1\right),\\
	\dd y_t
	&=
	\sin\theta_t
	\left(v_1\,\dd t+\sqrt{2D}\circ\dd W_t^1\right),\\
	\dd\theta_t
	&=
	\sqrt{2D}\circ\dd W_t^2.
\end{align}
The experiment is performed on $\R^2\times S^1$, without
reflection (though the developed theory strictly holds true only in a bounded domain). All particles start at $(x,y,\theta)=(2,0,\pi/2)$.
At this configuration, the heading is perpendicular to the
direction toward the origin, so the horizontal gradient vanishes.

\begin{figure}[htbp]
	\centering
	\includegraphics[width=\textwidth]
	{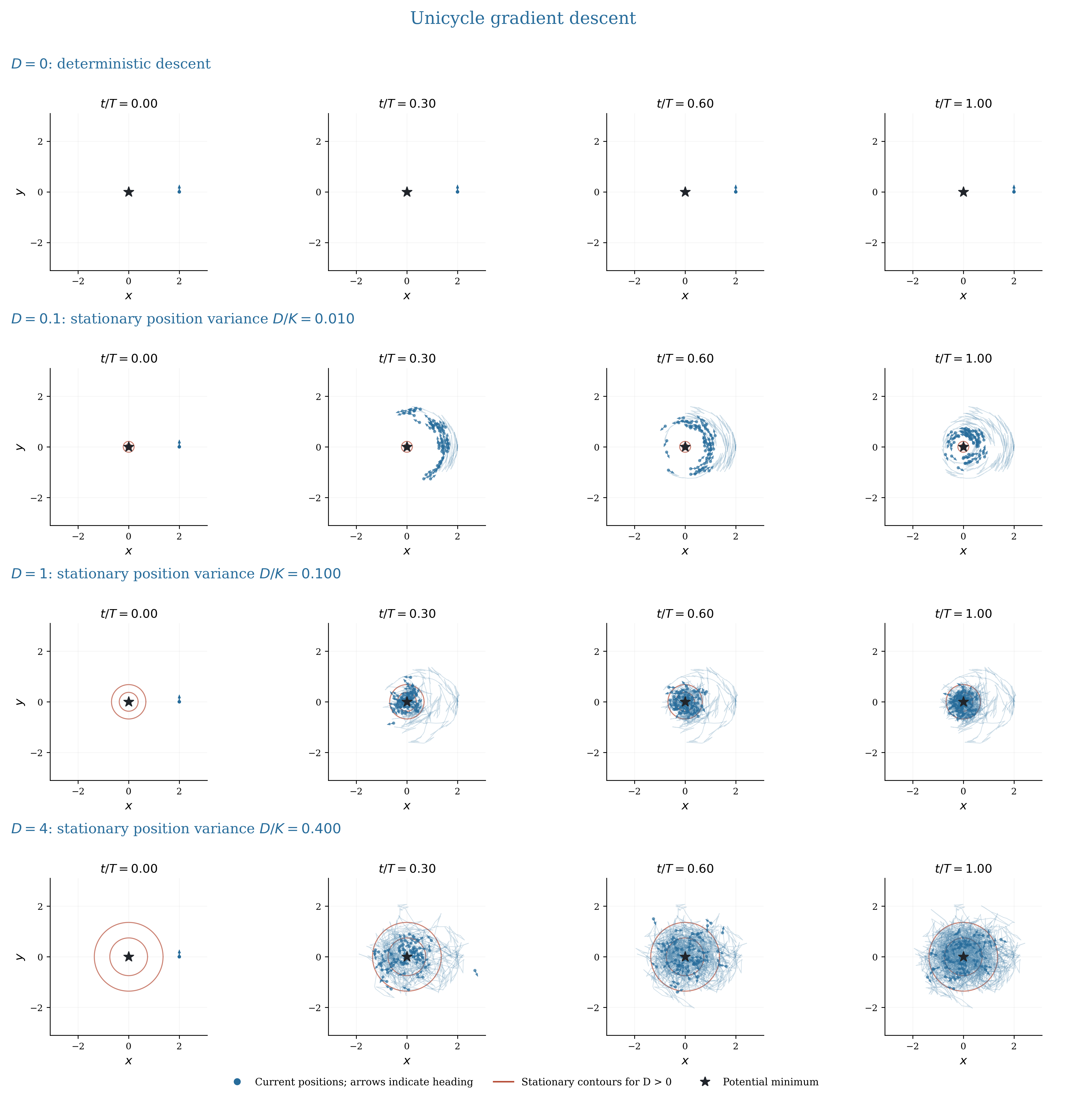}
	\caption{Unicycle horizontal-gradient descent with fixed
		potential strength $K=10$ and increasing diffusion strength.
		Rows correspond to $D=0$, $0.1$, $1$, and $4$; columns show
		successive times. Blue points indicate positions, arrows show
		headings, and faint curves trace selected particles.
		For $D>0$, red circles enclose $50\%$ and $90\%$ of the
		stationary position probability. }
	\label{fig:unicycle-diffusion-comparison}
\end{figure}

For $D=0$, the particles remain at their initial configuration.
For $D>0$, fluctuations in heading allow motion in previously
unavailable directions. The invariant position density is
\[
\rho_{\infty,D}^{xy}(x,y)
=
\frac{K}{2\pi D}
\exp\!\left(-\frac{K(x^2+y^2)}{2D}\right),
\]
while the invariant heading distribution is uniform on $S^1$.
Figure~\ref{fig:unicycle-diffusion-comparison} compares
$D=0$, $0.1$, $1$, and $4$ over the horizon $T=8$.
Increasing the diffusion strengthens exploration but also broadens
the stationary distribution. One can see the trade-off from using too much diffusion here and how noise is not all rosy. Conversely, a smaller positive
diffusion gives a more concentrated target, while slowing the
exploration of headings. 

\paragraph{Further Reading}

Stability properties of Fokker-Planck equations and the associated inequalities can be found in \cite{bakry2014analysis}. The solutions to the stabilizing nonholonomic Fooker-Planck equations can be found in \cite{elamvazhuthi2023density,elamvazhuthi2025score}. More about topological obstructions for stabilization as establish by Brockett \cite{brockett1983asymptotic} can be found in \cite{jongeneel2023topological}.

\chapter{Feedback Control Using Denoising Diffusion}
\label{chap:denoising-diffusion}
\label{chap:denoising-diffusion-control}

The noising--denoising approach developed in the machine-learning
literature \cite{song2020score} provides a way to generate new
samples from a data distribution using available samples.
This chapter adapts that idea to control and planning problems.
Instead of generating data, we seek a feedback controller that
steers a system toward a target configuration.

The construction combines the noising--denoising viewpoint of
Chapter~\ref{chap:flow-matching} with the Fokker--Planck dynamics
of Chapter~\ref{chap:fokker-planck}. Historically,  noising-denoising viewpoint to sampling \cite{song2020score} came before flow matching. The new(or older) ingredient is the
\emph{score}, which was the precursor to the conditional-mean
controller used in flow matching. The score allows us to rewrite
the Fokker--Planck equation as a continuity equation and associate
it with a deterministic ODE, known sometimes as the {\it probabilistic flow ODE}. Reversing this ODE gives a
feedback law for path planning and control.

\section{Reversing Reflected Diffusions for Path Planning}
\label{sec:denoising-general}

Suppose a robot evolves in a bounded, connected domain
$\Omega\subset\R^n$ with sufficiently smooth boundary. Your food delivery robot bringing your chicken briyani might disagree with this assumption, but so far they don't have voting rights. We first
consider the fully actuated system
\[
\dot{x}=u.
\]
Given an initial configuration $x_0\in\Omega$ and a target
$x_\star\in\Omega$, we seek a feedback law $u=u(t,x)$ whose
trajectory $\gamma$ remains in $\Omega$ and satisfies
\begin{equation}
	\gamma(0)=x_0,
	\qquad
	\gamma(T)=x_\star.
\end{equation}

As in Chapter~\ref{chap:flow-matching}, the construction has two
phases:
\begin{enumerate}[label=\arabic*.]
	\item A forward noising phase that explores $\Omega$ through
	diffusion, starting from the target.
	\item A reverse denoising phase that uses the noising
	information to construct feedback trajectories toward
	the target.
\end{enumerate}
We formulate this construction at the level of probability laws
and then use the resulting velocity field as a feedback controller.

\paragraph{Forward noising}

For the fully actuated system, we use reflected Brownian motion:
\begin{equation}
	\dd X_t^{\rm f}
	=
	\sqrt{2D}\,\dd W_t
	+
	n_{\rm in}(X_t^{\rm f})\,\dd L_t,
	\qquad
	X_0^{\rm f}\sim\mu_\star,
	\label{eq:lecture8-forward-reflected}
\end{equation}
where $D>0$ is the diffusion coefficient, $n_{\rm in}$ is the
inward unit normal, and $L_t$ is the boundary local time.
The forward process starts from the target law $\mu_\star$.
For a point target, we take $\mu_\star=\delta_{x_\star}$.
For $t>0$, write its law as let $\rho_t^{\rm f}(x)$ be the law of the process.
As discussed in Chapter~\ref{chap:fokker-planck}, this density
satisfies the heat equation with Neumann boundary conditions:
\begin{align}
	\partial_t\rho_t^{\rm f}
	&=
	D\Delta\rho_t^{\rm f},\\
	n_{\rm out}\cdot\nabla\rho_t^{\rm f}
	&=0
	\qquad\text{on }\partial\Omega.
\end{align}
Since $\Omega$ is bounded and connected, the density converges
to the uniform density $1/|\Omega|$ as $t\to\infty$.

Define the \emph{score} by
\begin{equation}
	s(t,x)
	:=
	\nabla\log\rho_t^{\rm f}(x).
\end{equation}
Using $\nabla\rho_t^{\rm f}=\rho_t^{\rm f}s(t,\cdot)$,
we can rewrite the heat equation as
\begin{equation}
	\partial_t\rho_t^{\rm f}
	=
	\nabla\cdot\left(D\rho_t^{\rm f}s(t,\cdot)\right).
\end{equation}
Comparing this with the continuity equation
\[
\partial_t\rho=-\nabla\cdot(v\rho),
\]
we identify the \emph{probability-flow ODE}
\begin{equation}
	\dot X_t=-D\,s(t,X_t).
	\label{eq:probability-flow-heat}
\end{equation}
Its velocity field transports the forward densities
$\rho_t^{\rm f}$ according to the same continuity equation.
The Neumann boundary condition implies
$s(t,x)\cdot n_{\rm out}(x)=0$ on $\partial\Omega$, so the
velocity is tangent to the boundary.

We must note here that the diffusion and the probability-flow ODE remain very different at the individual-trajectory level. The two agree only at the level of the probability laws. Loosely speaking, we have the systems that are identical macroscopically, but very different microscopically.

\paragraph{Reverse denoising}

Reversing \eqref{eq:probability-flow-heat} gives the feedback law
\begin{equation}
	\boxed{
		u(t,x)=D\,s(T-t,x),}
	\label{eq:reverse-heat-score}
\end{equation}
If the reverse process starts with density $\rho_T^{\rm f}$,
then its density at reverse time $t<T$ is
$\rho_{T-t}^{\rm f}$. Consequently, its law approaches
$\mu_\star$ as $t\uparrow T$. If the target is a point, this means
convergence to $\delta_{x_\star}$.

When $T$ is sufficiently large, $\rho_T^{\rm f}$ is approximately
uniform. Initializing uniformly therefore approximates the exact
time-reversal.

\paragraph{Learning the score}

We would like to learn $s(t,x)$ without first estimating
$\rho_t^{\rm f}$. Consider the mean-squared score loss
\begin{equation}
	\mathcal L_{\rm SM}(\theta)
	=
	\E_{X\sim\mu}
	\left[
	\left\|s_\theta(X)-\nabla\log\rho(X)\right\|^2
	\right].
\end{equation}
As, such this loss cannot be evaluated in a sample-based way due to the unknown $-\nabla\log\rho(X)$. To arrive at a sample-based derivation, we expand the square,
\begin{equation}
	\mathcal L_{\rm SM}(\theta)
	=
	\E_{X\sim\mu}\left[\|s_\theta(X)\|^2\right]
	-
	2\int_\Omega s_\theta(x)\cdot\nabla\rho(x)\dd x
	+
	C,
\end{equation}
where $C$ is independent of $\theta$. If the boundary term
vanishes, integration by parts yields
\begin{equation}
	\mathcal L_{\rm SM}(\theta)
	=
	\E_{X\sim\mu}
	\left[
	\|s_\theta(X)\|^2
	+
	2\nabla\cdot s_\theta(X)
	\right]
	+
	C.
	\label{eq:hyvarinen-score-loss}
\end{equation}
Thus, the score can be learned from samples without evaluating
the density or its gradient.

\begin{remark}[Boundary terms]
	The integration-by-parts formula is
	\[
	\int_\Omega s_\theta\cdot\nabla\rho\,\dd x
	=
	-\int_\Omega\rho\,\nabla\cdot s_\theta\,\dd x
	+
	\int_{\partial\Omega}
	\rho\,s_\theta\cdot n_{\rm out}\,\dd S.
	\]
	Hence \eqref{eq:hyvarinen-score-loss} requires the final term
	to vanish. One way to ensure this is to impose
	$s_\theta\cdot n_{\rm out}=0$ on $\partial\Omega$, consistent
	with the boundary condition satisfied by the exact score.
\end{remark}

Returning to the time-dependent problem, we use the
score-matching objective
\begin{equation}
	\mathcal L_{\rm SM}(\theta)
	=
	\int_0^T
	\E_{X\sim\rho_t^{\rm f}}
	\left[
	\|s_\theta(t,X)\|^2
	+
	2\nabla\cdot s_\theta(t,X)
	\right]\dd t,
	\label{eq:time-dependent-score-loss}
\end{equation}
with the boundary treatment described above. The divergence is
taken with respect to the state variable. As in flow matching,
the integral is approximated by sampling times and states from
the forward trajectories. The learned feedback is then
\[
u_\theta(t,x)=D\,s_\theta(T-t,x).
\]

The following algorithm expands on how this estimation is used to construct a control law.
\begin{tcolorbox}[algorithmbox,title=\textbf{Algorithm: Feedback Control Using Denoising Diffusion}]
	\textbf{Setup.}
	Consider the fully actuated system $\dot{x}=u$ on $\Omega$.
	Choose a target law $\mu_\star$, a diffusion coefficient
	$D>0$, and a horizon $T$.
	For a point target, take $\mu_\star=\delta_{x_\star}$.
	
	\begin{enumerate}[label=\arabic*.]
		\item Generate forward noising trajectories using
		\eqref{eq:lecture8-forward-reflected}, initialized
		from $\mu_\star$. Store samples
		\[
		\left\{X_{t_k}^{\mathrm f,(j)}\right\}_{j=1}^{N}
		\]
		at sampled times $t_k\in(0,T)$.
		
		\item Initialize a score network
		\[
		s_\theta:(0,T]\times\overline{\Omega}\to\R^n,
		\]
		satisfying
		$s_\theta(t,x)\cdot n_{\rm out}(x)=0$
		on $\partial\Omega$.
		
		\item Train the network by minimizing the empirical
		loss corresponding to
		\[
		\mathcal L_{\rm SM}(\theta)
		=
		\int_0^T
		\E_{X\sim\mu_t^{\rm f}}
		\left[
		\|s_\theta(t,X)\|^2
		+
		2\nabla\cdot s_\theta(t,X)
		\right]\dd t.
		\]
		
		\item From the desired initial configuration $x_0$,
		apply the learned feedback
		\[
		u_\theta(t,x)=D\,s_\theta(T-t,x)
		\]
		and integrate
		\[
		\dot{x}(t)=u_\theta(t,x(t)),
		\qquad x(0)=x_0,
		\qquad 0\leq t<T.
		\]
	\end{enumerate}
	
	With the exact score and initial law $\mu_T^{\rm f}$,
	the reverse law approaches $\mu_\star$ as $t\uparrow T$.
\end{tcolorbox}

\paragraph{Numerical example}

We consider the system confined in the square
$[-1,1]^2$, with three circular obstacles removed. 
The target is $x_\star=(0,-0.7)^{\mathsf T}$.

Starting from the target, we simulate a reflected Brownian
approximation with diffusion coefficient $D=0.5$ over the
horizon $T=3$. A neural network is trained to learn the score using the
score-matching objective \eqref{eq:time-dependent-score-loss},
with its output constrained to be tangent to the boundaries.
We then apply the learned feedback
\[
u_\theta(t,x)=D\,s_\theta(T-t,x)
\]
from independently sampled, uniformly distributed free-space
initial configurations.

\begin{figure}[htbp]
	\centering
	\includegraphics[width=\textwidth]
	{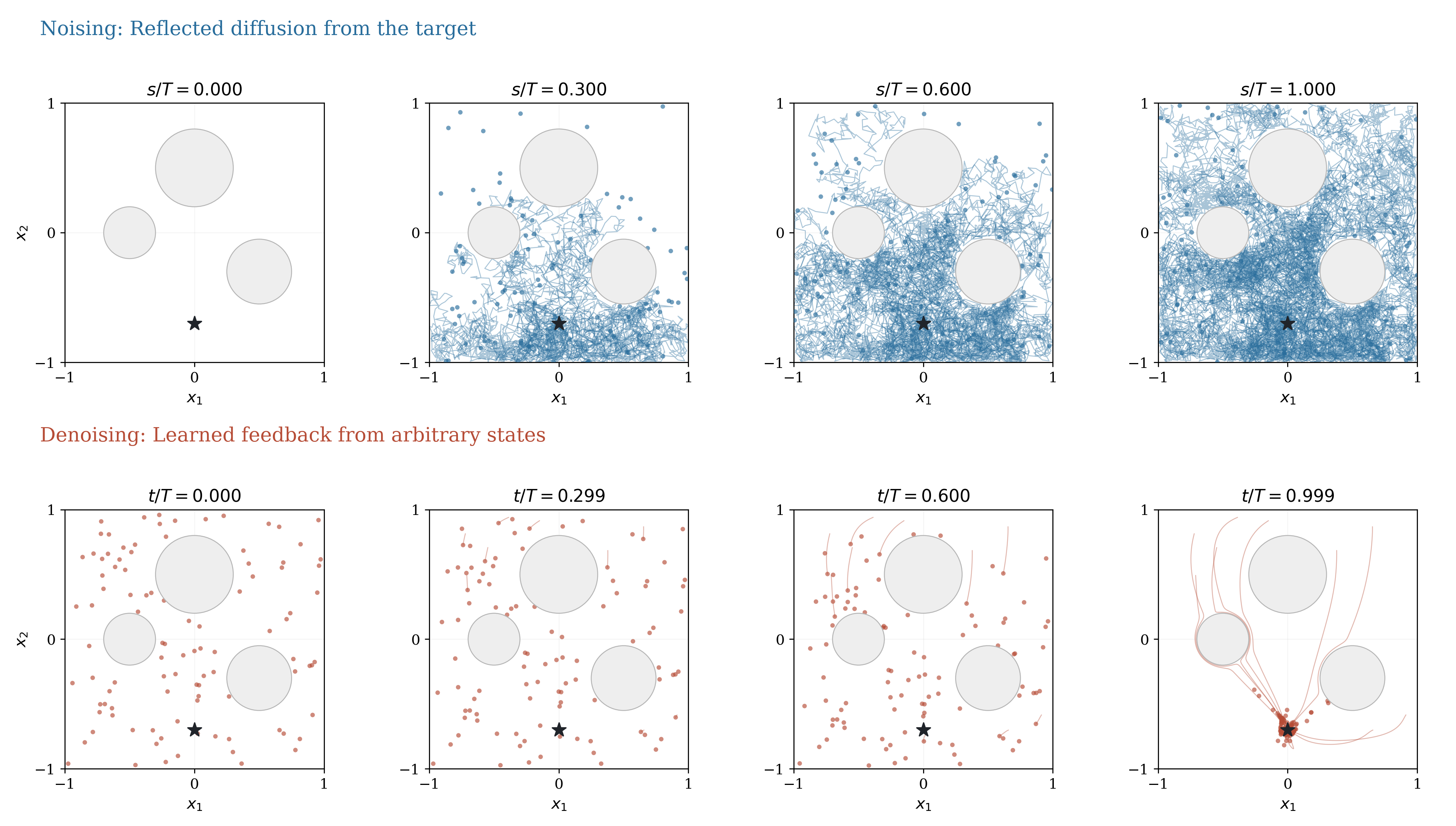}
	\caption{Diffusion-based feedback planning around circular
		obstacles. Top: forward noising from the target.
		Bottom: deterministic denoising under the learned score
		feedback, initialized uniformly in free space.
		Points show current states, and curves trace selected
		particles up to the displayed time. Gray disks represent
		obstacles, and the star marks the target.}
	\label{fig:obstacle-score-matching}
\end{figure}

Figure~\ref{fig:obstacle-score-matching} compares the forward
noising trajectories with the resulting deterministic
denoising trajectories. The forward diffusion spreads from the
target around the obstacles, while the learned feedback uses
this information to guide states back toward the target.

\paragraph{Changing the initial distribution.}

For a fixed target point, the initial distribution for the reverse process does not need to be the uniform distribution. One can ask, why should the same feedback work when initialized near a
prescribed configuration. For a connected domain with sufficiently
smooth boundary, the solution of the heat equation is strictly positive
at positive times. In particular,
$\rho_T^{\rm f}(x)>0$ throughout $\Omega$.

Let $Z_t$ denote the reverse probability-flow process. If the initial law is $\rho_T(x)$, its law approaches $\delta_{x_\star}$ as $t\uparrow T$.
Now choose another initial density $q_0$. We can change the
probability assigned to trajectories by defining
\[
\frac{\dd\mathbb Q}{\dd\mathbb P}
=
\frac{q_0(Z_0)}{\rho_T^{\rm f}(Z_0)}.
\]
Under $\mathbb Q$, the initial density is $q_0$, but the
trajectories still satisfy the same equation
\[
\dot Z_t=D\,s(T-t,Z_t).
\]

Suppose that $q_0/\rho_T^{\rm f}\leq C$ for some finite constant
$C$. Then, for every $r>0$,
\[
\mathbb Q\bigl(|Z_t-x_\star|>r\bigr)
\leq
C\,\mathbb P\bigl(|Z_t-x_\star|>r\bigr)
\longrightarrow0.
\]
Thus, the concentration at the target is preserved under this
change of initialization. This provides an informal justification for planning
from a localized initial density near any configuration.

\section{Denoising for Driftless Control-Affine Systems}

The fully actuated probability-flow ODE does not incorporate dynamical constraints of a control system. We now jump into the same score-based 
for driftless control systems of the form,
\begin{equation}
	\dot x
	=
	\sum_{i=1}^m u_i g_i(x).
	\label{eq:lecture8-driftless-system}
\end{equation}
using the {\it horizontal} operators introduced in Chapter~\ref{chap:fokker-planck}. Unlike the fully actuated case, noise can only be injected along the admissible control directions $g_i$.

Recall the directional operators $\mathcal Y_i$ and $\mathcal Y_i^*$ from \eqref{eq:global-directional-operators}, and write
\begin{equation}
    d_i:=\nabla\cdot g_i.
\end{equation}
We revist the reflected SDE
\begin{equation}
    \dd X_t^{\rm f}
    =
    \sum_{i=1}^m v_i(X_t^{\rm f})g_i(X_t^{\rm f})\dd t
    +
    \sqrt{2}\sum_{i=1}^m g_i(X_t^{\rm f})\circ\dd W_t^i
    +
    n_{\rm in}(X_t^{\rm f})\dd L_t.
    \label{eq:lecture8-nonholonomic-forward}
\end{equation}
Its Fokker--Planck equation is
\begin{equation}
    \partial_t \rho_t
    =
    \sum_{i=1}^m(\mathcal Y_i^*)^2\rho_t
    +
    \sum_{i=1}^m\mathcal Y_i^*(v_i \rho_t).
\end{equation}

As derived earlier, define
\begin{equation}
    U_i:=v_i-d_i.
\end{equation}
Then
\begin{equation}
    \partial_t \rho_t
    =
    -\sum_{i=1}^m\mathcal Y_i^*\mathcal Y_i \rho_t
    +
    \sum_{i=1}^m\mathcal Y_i^*(U_i \rho_t).
    \label{eq:lecture8-transformed-fp}
\end{equation}
If a smooth strictly positive noise density $\rho_{\rm noise}$ is desired, choose
\begin{equation}
    U_i
    =
    \mathcal Y_i\log \rho_{\rm noise},
    \qquad
    v_i
    =
    d_i+\mathcal Y_i\log \rho_{\rm noise}.
\end{equation}
For the uniform distribution, $U_i=0$ and $v_i=d_i$.

\paragraph{The control-constrained probability-flow ODE.}

Define the generalized or horizontal score
\begin{equation}
    s_i(t,x)
    :=
    \mathcal Y_i\log \rho_t(x)
    =
    \frac{\mathcal Y_i \rho_t(x)}{\rho_t(x)}.
\end{equation}
Equation \eqref{eq:lecture8-transformed-fp} can be written as the continuity equation associated with
\begin{equation}
    \dot X_t
    =
    \sum_{i=1}^m
    \left[
        -s_i(t,X_t)+U_i(X_t)
    \right]
    g_i(X_t).
\end{equation}
Therefore the time-reversed deterministic dynamics are
\begin{equation}
    \boxed{
    \dot Z_t
    =
    \sum_{i=1}^m
    \left[
        s_i(T-t,Z_t)-U_i(Z_t)
    \right]
    g_i(Z_t).
    }
    \label{eq:lecture8-reverse-nonholonomic}
\end{equation}
This ODE is automatically compatible with the nonholonomy of the system by respecting the control constraints. It gives the feedback law
\begin{equation}
    u_i(t,x)
    =
    s_i(T-t,x)-U_i(x).
\end{equation}
In the particularly useful uniform-noise case, this reduces to
\begin{equation}
    u_i(t,x)
    =
    s_i(T-t,x).
\end{equation}

\paragraph{Horizontal score matching.}

A network $s_\theta=(s_{\theta,1},\ldots,s_{\theta,m})$ can be trained to approximate the horizontal score. Expanding the $L^2$ regression loss and using the formal adjoint gives, up to a parameter-independent constant,
\begin{equation}
    \boxed{
    \mathcal L_{\rm HSM}(\theta)
    =
    \int_0^T
    \E_{\rho_t}
    \left[
        \sum_{i=1}^m
        \left(
            |s_{\theta,i}(t,X)|^2
            +
            2\nabla\cdot\bigl(g_i(X)s_{\theta,i}(t,X)\bigr)
        \right)
    \right]\dd t.
    }
    \label{eq:horizontal-score-matching-loss}
\end{equation}
As in the Euclidean case, boundary terms must be handled consistently on bounded domains, but we will shove it under the rug, risking scorn of the careful reader or their LLM.

\paragraph{Unicycle example.}

We revisit the unicycle control problem,
\begin{align}
    \dot x_1&=u_1\cos x_3,\\
    \dot x_2&=u_1\sin x_3,\\
    \dot x_3&=u_2.
\end{align}

\begin{figure}[htbp]
	\centering
	\includegraphics[width=\textwidth]
	{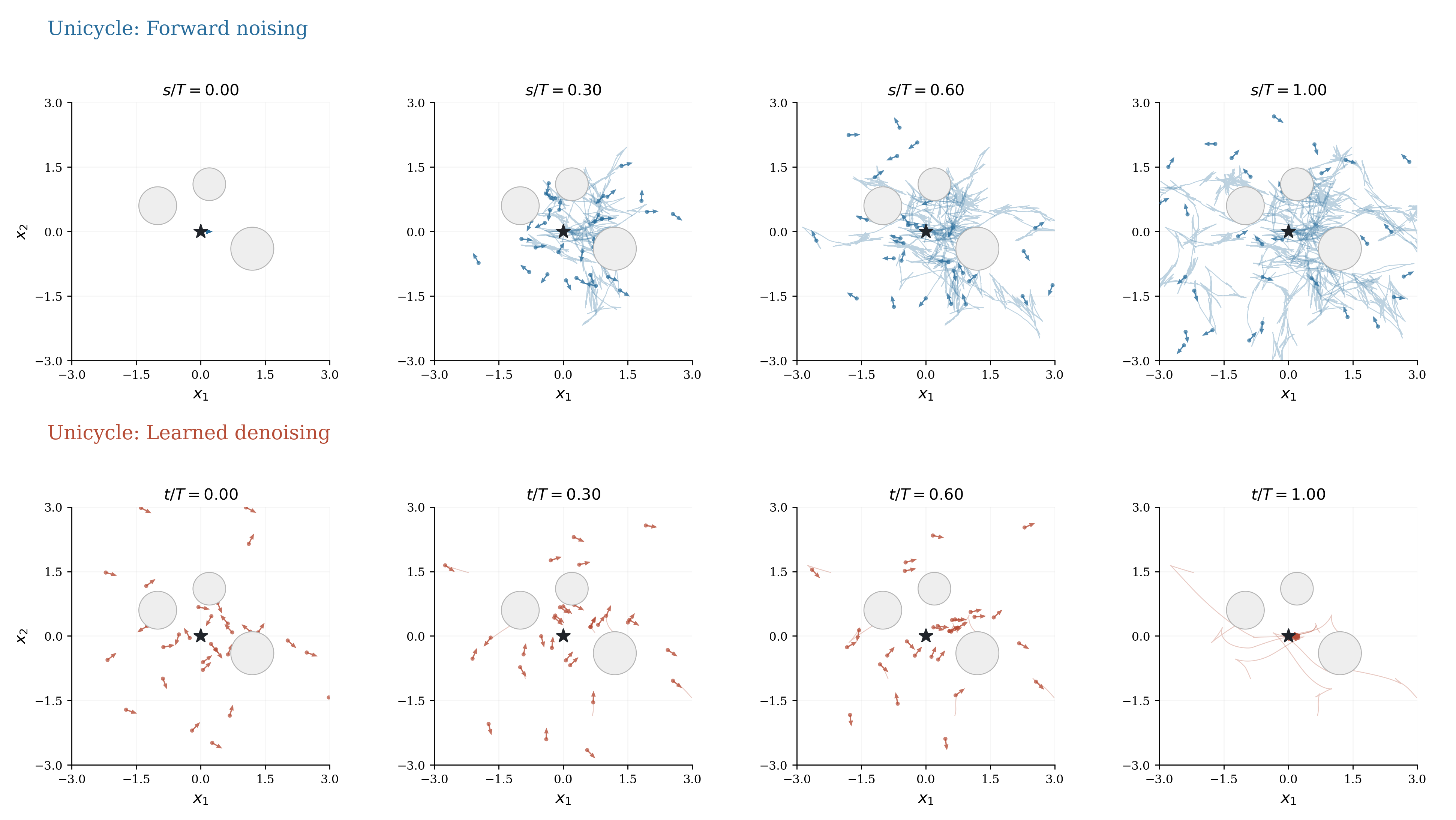}
	\caption{Noising and denoising for the unicycle.
		Top: forward diffusion initialized at the target.
		Bottom: deterministic denoising using the score network
		selected by the lowest validation loss.
		Points show positions, arrows indicate headings, and faint
		curves trace selected trajectories. The star marks the
		target position; gray disks represent obstacles.}
	\label{fig:unicycle-noising-denoising}
\end{figure}

Figure~\ref{fig:unicycle-noising-denoising} shows the denoising based controller in action for this system. Starting from the target
configuration, the forward diffusion explores the state space
through the admissible control directions. The learned score
then defines a deterministic feedback law used in the denoising
phase to steer independently initialized configurations toward
the target. The arrows indicate the heading, while the traces
show how the nonholonomic constraint shapes the trajectories.

\paragraph{Further reading.}

The material in this section follows the score-matching diffusion approach to feedback control developed by Karthik Elamvazhuthi, Darshan Gadginmath, and Fabio Pasqualetti, \emph{Score Matching Diffusion Based Feedback Control and Planning of Nonlinear Systems}, arXiv:2504.09836. Closely related stochastic time-reversal approaches include work by Yuhang Mei, Amirhossein Taghvaei, and Ali Pakniyat. For score-based generative modeling, a standard reference is Yang Song, Jascha Sohl-Dickstein, Diederik P. Kingma, Abhishek Kumar, Stefano Ermon, and Ben Poole, \emph{Score-Based Generative Modeling through Stochastic Differential Equations}.

\section{Denoising Linear Time-Invariant Systems  \\ (Just for fun. Don't do this at Home!)} 
\label{sec:denoising-lti}

In this section, we specialize diffusion-based feedback control to linear time-invariant systems. Since classical stabilization methods for LTIs are already extremely well developed, this mostly a just-for-fun section, as the probability laws and scores can be computed explicitly. More generally, it also shows how the method can be used to stabilize a LTI to multimodal target equilibria or distributions, an objective that is inherently nonlinear.

Consider
\begin{equation}
    \dot x
    =
    Ax+Bu,
    \label{eq:lti-control-denoising}
\end{equation}
where $A\in\R^{n\times n}$ and $B\in\R^{n\times m}$, with usual assumptions on controllability.

To denoise using the original drift $Ax$, use the time-reversed drift in the forward process:
\begin{equation}
    \dd X_t
    =
    -AX_t\dd t
    +
    \sqrt{2}B\dd W_t,
    \qquad
    X_0\sim\mu_\star.
    \label{eq:lti-forward-noising}
\end{equation}
We recall the density evolves according to,
\begin{equation}
    \rho_t(x)
    =
    \int_{\R^n}K(t,x,y)\mu_\star(\dd y),
\end{equation}
where
\begin{equation}
    K(t,x,y)
    =
    \frac{
        \exp\!\left[
            -\frac12
            (x-e^{-tA}y)^{\mathsf T}
            \Sigma_t^{-1}
            (x-e^{-tA}y)
        \right]
    }{
        (2\pi)^{n/2}(\det\Sigma_t)^{1/2}
    },
    \label{eq:lti-noising-kernel}
\end{equation}
and
\begin{equation}
    \boxed{
    \Sigma_t
    =
    2\int_0^t
    e^{-sA}BB^{\mathsf T}e^{-sA^{\mathsf T}}\dd s.
    }
    \label{eq:lti-noising-gramian}
\end{equation}
For a controllable pair $(A,B)$, $\Sigma_t$ is positive definite for every $t>0$.

In this case, the Fokker--Planck equation specializes to,
\begin{equation}
    \partial_t \rho_t
    =
    \nabla\cdot(Ax\,\rho_t)
    +
    \nabla\cdot\left(BB^{\mathsf T}\nabla \rho_t\right).
\end{equation}
The corresponding probability-flow ODE is
\begin{equation}
    \dot x
    =
    -Ax
    -
    BB^{\mathsf T}\nabla\log \rho_t(x).
\end{equation}
Reversing time yields
\begin{equation}
    \dot x
    =
    Ax
    +
    BB^{\mathsf T}\nabla\log \rho_{T-t}(x).
\end{equation}
This has the form \eqref{eq:lti-control-denoising} with feedback control
\begin{equation}
    u(t,x)
    =
    B^{\mathsf T}\nabla\log \rho_{T-t}(x)
    \label{eq:lti-denoising-feedback}
\end{equation}

\subsection{Single Target Equilibrium and Minimum Energy Control}
In order to steer to the origin we take
\begin{equation}
	\mu_\star=\delta_0.
\end{equation}
Then
\begin{equation}
	\rho_t(x)
	=
	\frac{
		\exp\!\left(-\frac12x^{\mathsf T}\Sigma_t^{-1}x\right)
	}{
		(2\pi)^{n/2}(\det\Sigma_t)^{1/2}
	}.
\end{equation}
The score simplifies to
\begin{equation}
		\nabla\log \rho_t(x)
		=
		-\Sigma_t^{-1}x.	
	\label{eq:lti-gaussian-score}
\end{equation}

and \begin{equation}
	\boxed{u(t,x)
		=
		-B^{\mathsf T}	\Sigma_{T-t}^{-1}x}
\end{equation}
Hence, in this case the score has no density dependence and we get a time-dependent state-feedback law. The appearance of the controllability Gramian is reminiscent of minimum-energy control.

We can compare our denoising controller with the minimum energy controller. Lets recall the minimum-energy
control problem for the system 
\begin{equation}
	\begin{aligned}
		\inf_{u}\quad&
		\frac{1}{2}\int_t^T |u(s)|^2\,\dd s,\\
		\text{subject to}\quad&
		\dot{x}(s)=Ax(s)+Bu(s),\\
		&x(t)=x,\qquad x(T)=0.
	\end{aligned}
\end{equation}
If $(A,B)$ is controllable, the corresponding minimum-energy feedback is
\begin{equation}
		u_{\mathrm{ME}}(t,x)
		=
		-B^{\mathsf T}e^{\tau A^{\mathsf T}}
		W_\tau^{-1}e^{\tau A}x,
		\qquad \tau=T-t.
\end{equation}
where $W_\tau$ is the \emph{controllability Gramian},
\begin{equation}
	W_\tau
	=
	\int_0^\tau
	e^{sA}BB^{\mathsf T}e^{sA^{\mathsf T}}\,\dd s.
\end{equation}
We can rewrite this using the Grammian of the time-reversed system and introduce the {\it backward Gramian}
\begin{equation}
	G_\tau
	=
	\int_0^\tau
	e^{-sA}BB^{\mathsf T}e^{-sA^{\mathsf T}}\,\dd s,
\end{equation}
to obtain
\begin{equation}
	u_{\mathrm{ME}}(t,x)
	=
	-B^{\mathsf T}G_{T-t}^{-1}x.
\end{equation}

It is interesting to note that,
$\Sigma_\tau=2G_\tau$. Therefore, the denoising
controller above satisfies
\begin{equation}
	u(t,x)=\frac{1}{2}u_{\mathrm{ME}}(t,x).
\end{equation}
Thus, the two feedback laws involve the same Gramian but have
different gains. Lets look at what happens to the system for the single integrator
$A=0$, $B=I$. In this case,
\begin{equation}
	u(t,x)=-\frac{x}{2(T-t)},
	\qquad
	u_{\mathrm{ME}}(t,x)=-\frac{x}{T-t}.
\end{equation}
Starting from $x(0)=x_0$, their trajectories are respectively
\begin{equation}
	x(t)=\sqrt{1-\frac{t}{T}}\,x_0,
	\qquad
	x_{\mathrm{ME}}(t)=\left(1-\frac{t}{T}\right)x_0.
\end{equation}
Both approach the origin, but for $x_0\neq0$ the denoising
control has infinite energy over $[0,T)$, while the
minimum-energy control is constant and has finite energy.

\subsection{More General Objectives}

For a more general target distribution $\mu_\star$, we want to
transport to the distribution itself, rather than merely steer the
state toward a single point. The reverse process should therefore
be initialized from a distribution close to the terminal noising
law $\mu_T$. This requires the additional assumption is that $-A$ is Hurwitz.
As discussed in Chapter~\ref{chap:fokker-planck}, the noising
process with drift $-Ax$ then converges, for a controllable pair
$(A,B)$, to a Gaussian distribution
$\mathcal N(0,\Sigma_\infty)$, where
\begin{equation}
	\Sigma_\infty
	=
	2\int_0^\infty
	e^{-sA}BB^{\mathsf T}e^{-sA^{\mathsf T}}\,\dd s.
\end{equation}
For sufficiently large $T$, this Gaussian provides an
approximation to $\mu_T$ from which we can readily sample. While denoising diffusion gives us a way to transport an arbitrary (sufficiently regular) distribution to a target distribution, for linear systems, we are restricted to initializing at $ N(0,\Sigma_\infty)$.

The simplest, but still interesting generalization, is to consider a sum of Dirac as the target,
\begin{equation}
    \rho_0
    =
 \frac12\delta_a
    +
    \frac12\delta_b,
\end{equation}
with $a,b\in\ker A$, so that the target points are equilibria of the uncontrolled system. Then
\begin{equation}
    \rho_t(x)
    =
    \frac12K(t,x,a)
    +
    \frac12K(t,x,b).
\end{equation}
Since $e^{-tA}a=a$ and $e^{-tA}b=b$, the score is
\begin{equation}
    \nabla\log \rho_t(x)
    =
    -
    \frac{
        \Sigma_t^{-1}(x-a)K(t,x,a)
        +
        \Sigma_t^{-1}(x-b)K(t,x,b)
    }{
        K(t,x,a)+K(t,x,b)
    }.
    \label{eq:lti-mixture-score}
\end{equation}
The reverse feedback remains
\begin{equation}
    u(t,x)
    =
    B^{\mathsf T}\nabla\log \rho_{T-t}(x).
\end{equation}
Thus the same construction produces a time-varying feedback law whose terminal distribution is bimodal.

\paragraph{Numerical illustration.}

Consider the four-dimensional system
\begin{align}
	\dot x_1 &= x_2,\\
	\dot x_2 &= -x_1+\frac{3}{2}x_2+u_1,\\
	\dot x_3 &= x_4,\\
	\dot x_4 &= -x_3+\frac{3}{2}x_4+u_2.
\end{align}
The pair $(A,B)$ is controllable, and the eigenvalues of $-A$
are $-3/4\pm i\sqrt{7}/4$, each repeated twice. Thus, $-A$
is Hurwitz, and the backward noising process
\[
\dd X_s=-AX_s\,\dd s+\sqrt{2}B\,\dd W_s
\]
converges to a stationary Gaussian distribution.

\begin{figure}[htbp]
	\centering
	\includegraphics[width=\textwidth]
	{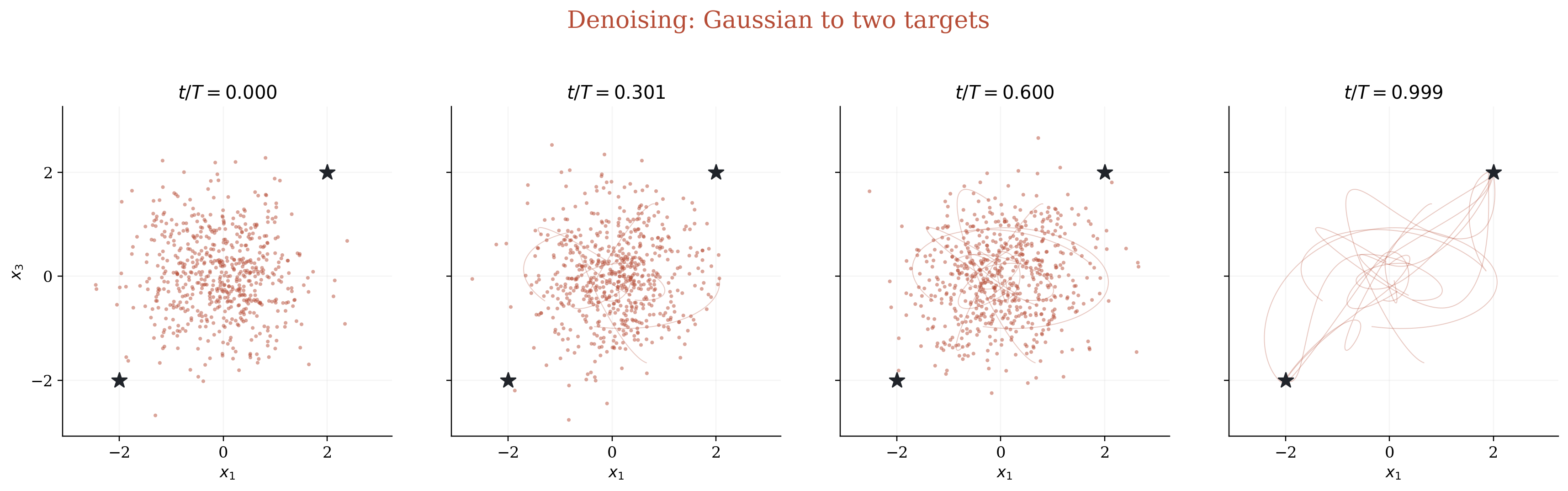}
	\caption{Noising and denoising for a four-dimensional
		controllable linear system. The plot shows eterministic denoising from 
	Gaussian samples.
		Points show the position coordinates $(x_1,x_3)$,
		faint curves trace selected trajectories, and stars mark
		the target positions. }
	\label{fig:lti-denoising-examples}
\end{figure}

Figure~\ref{fig:lti-denoising-examples} shows the evolution
in the position coordinates $(x_1,x_3)$ over a horizon $T=6$.
During denoising, independent samples
from this Gaussian separate and concentrate near the two
target configurations. This illustrates steering toward two targets in finite time.

\paragraph{More general target distributions.}

For a general target law $\mu_\star$, the integral against $K(t,x,y)$ may be expensive to evaluate, and the score need not be available in closed form. One can instead simulate \eqref{eq:lti-forward-noising} and learn
\begin{equation}
    s_\theta(t,x)
    \approx
    \nabla\log \rho_t(x)
\end{equation}
from the simulated samples using score matching, as we did for driftless control systems.

The reverse controller is then
\begin{equation}
    u_\theta(t,x)
    =
    B^{\mathsf T}s_\theta(T-t,x).
    \label{eq:lti-learned-score-control}
\end{equation}
The general algorithm for denoising an LTI system is as follows.
\begin{tcolorbox}[algorithmbox,title=\textbf{Algorithm: Score Matching for LTI Feedback Control}]
	\textbf{Setup.}
	Consider the controllable linear system
	\[
	\dot{x}=Ax+Bu,
	\qquad x\in\R^n,\quad u\in\R^m.
	\]
	Choose a target law $\mu_\star$ and a horizon $T>0$.
	For a point target, take $\mu_\star=\delta_{x_\star}$.
	
	\begin{enumerate}[label=\arabic*.]
		\item Generate forward noising trajectories using
		\[
		\dd X_t^{\rm f}
		=
		-AX_t^{\rm f}\dd t+\sqrt{2}B\dd W_t,
		\qquad X_0^{\rm f}\sim\mu_\star.
		\]
		Denote their density by $\rho_t^{\rm f}$ and store samples
		\[
		\left\{X_{t_k}^{\mathrm f,(j)}\right\}_{j=1}^{N}
		\]
		at positive time steps $t_k\leq T$.
		
		\item Initialize a score network
		\[
		s_\theta:[0,T]\times\R^n\to\R^n.
		\]
		
		\item Train the network by minimizing the empirical
		loss corresponding to
		\[
		\mathcal L_{\rm SM}(\theta)
		=
		\int_0^T
		\E_{X\sim\mu_t^{\rm f}}
		\left[
		\|s_\theta(t,X)\|^2
		+
		2\nabla\cdot s_\theta(t,X)
		\right]\dd t,
		\]
		where $\mu_t^{\rm f}(\dd x)=\rho_t^{\rm f}(x)\dd x$.
		The divergence is taken with respect to the state
		variable, with sufficient decay to justify
		integration by parts.
		
		\item From the desired initial configuration $x_0$,
		apply the learned feedback
		\[
		u_\theta(t,x)
		=
		B^{\mathsf T}s_\theta(T-t,x)
		\]
		and integrate
		\[
		\dot{x}(t)
		=
		Ax(t)+Bu_\theta(t,x(t)),
		\qquad x(0)=x_0,
		\qquad 0\leq t<T.
		\]
	\end{enumerate}
	
	With the exact score and initial law $\mu_T^{\rm f}$,
	the reverse law is $\mu_{T-t}^{\rm f}$ and approaches
	$\mu_\star$ as $t\uparrow T$.
	If $-A$ is Hurwitz and $T$ is sufficiently large,
	the stationary Gaussian law of the noising process
	provides an approximation to $\mu_T^{\rm f}$.
\end{tcolorbox}

\section{Denoising Nonlinear Systems with Drift}
\label{sec:denoising-drift}

We now consider a general control-affine system with drift,
\begin{equation}
    \dot x
    =
    f_0(x)
    +
    \sum_{i=1}^m u_i(t)g_i(x).
    \label{eq:drift-control-system-denoise}
\end{equation}
At this level of generality, constructing a forward process that converges to a simple universal noise distribution is difficult. For this reason, steering to general target distributions is challenging, just as in the case of LTIs. Nevertheless, when the target is a single point, a denoising construction is possible when the time-reversed dynamics explores the relevant state space.

\paragraph{Forward Noising}

The noising process must propagate outward from the target while reversing the uncontrolled drift. This leads to the control system
\begin{equation}
    \dot x
    =
    -f_0(x)
    +
    \sum_{i=1}^m u_i(t)g_i(x).
    \label{eq:time-reversed-controlled-system}
\end{equation}

If this system is globally reachable, then noise injected through the control directions can spread probability mass over the relevant backward-reachable region. We have already stated the corresponding support and exploration properties in Chapter \ref{chap:fokker-planck}. This means our noising process should be the following,
\begin{equation}
	\boxed{
		\dd X_t
		=
		-f_0(X_t)\dd t
		+
		\sqrt{2}\sum_{i=1}^m
		g_i(X_t)\circ\dd W_t^i.
	}
	\label{eq:drift-first-noising-sde}
\end{equation}
	
The corresponding Fokker--Planck equation is
\begin{equation}
	\partial_t\rho_t
	=
	-\mathcal Y_0^*\rho_t
	+
	\sum_{i=1}^m(\mathcal Y_i^*)^2\rho_t.
\end{equation}
Wherever $\rho_t>0$, we can rewrite this as
\begin{equation}
	\partial_t\rho_t
	=
	-\mathcal Y_0^*\rho_t
	+
	\sum_{i=1}^m
	\mathcal Y_i^*
	\left(
	\frac{\mathcal Y_i^*\rho_t}{\rho_t}\rho_t
	\right).
\end{equation}
Since $\mathcal Y_i^*\rho=-\nabla\cdot(g_i\rho)$,
this is the continuity equation associated with the
probability-flow ODE
\begin{equation}
	\boxed{
		\dot X_t
		=
		-f_0(X_t)
		+
		\sum_{i=1}^m
		\frac{\mathcal Y_i^*\rho_t}{\rho_t}(X_t)
		g_i(X_t).
	}
\end{equation}
Using the {\it adjoint score}
\[
\frac{\mathcal Y_i^*\rho_t}{\rho_t}
=
-\nabla\cdot g_i-\mathcal Y_i\log\rho_t,
\]
and reversing time, we obtain a feedback law for the
original control-affine system:
\begin{equation}
	\boxed{
		\dot x
		=
		f_0(x)+\sum_{i=1}^m u_i(t,x)g_i(x),
		\qquad
		u_i(t,x)
		=
		\nabla\cdot g_i(x)
		+
		\mathcal Y_i\log\rho_{T-t}(x).
	}
\end{equation}
For divergence-free control fields, this reduces to
$u_i(t,x)=\mathcal Y_i\log\rho_{T-t}(x)$,
recovering the horizontal-score feedback from the
driftless case.

\paragraph{Divergence correction}
The divergence term in the control  $\nabla\cdot g_i(x)$ might be unpleasant to some users, though the author cannot explain why. If this is the case, one can alternatively consider the following SDE as the noising process,
Writing $d_i=\nabla\cdot g_i$, consider the corrected
noising process
\begin{equation}
	\boxed{
		\dd X_t
		=
		-f_0(X_t)\dd t
		+
		\sum_{i=1}^m d_i(X_t)g_i(X_t)\dd t
		+
		\sqrt{2}\sum_{i=1}^m g_i(X_t)\circ\dd W_t^i.
	}
	\label{eq:drift-corrected-noising-sde}
\end{equation}
Since $\mathcal Y_i^*\rho=-\mathcal Y_i\rho-d_i\rho$,
the divergence correction gives the Fokker--Planck equation
\begin{equation}
	\partial_t\rho_t
	=
	-\mathcal Y_0^*\rho_t
	-
	\sum_{i=1}^m\mathcal Y_i^*\mathcal Y_i\rho_t.
	\label{eq:drift-corrected-fp}
\end{equation}
Again, wherever $\rho_t>0$, this can be written as
\begin{equation}
	\partial_t\rho_t
	=
	-\mathcal Y_0^*\rho_t
	-
	\sum_{i=1}^m
	\mathcal Y_i^*
	\left(
	\rho_t\,\mathcal Y_i\log\rho_t
	\right).
\end{equation}
The corresponding probability-flow ODE is therefore
\begin{equation}
	\dot X_t
	=
	-f_0(X_t)
	-
	\sum_{i=1}^m
	\mathcal Y_i\log\rho_t(X_t)\,g_i(X_t).
\end{equation}
Reversing time gives the feedback law for the original
control-affine system:
\begin{equation}
	\boxed{
		\dot x
		=
		f_0(x)+\sum_{i=1}^m u_i(t,x)g_i(x),
		\qquad
		u_i(t,x)
		=
		\mathcal Y_i\log\rho_{T-t}(x).
	}
	\label{eq:nonlinear-drift-reverse-ode}
\end{equation}

For the reader, who found the divergence term unpleasant, they will note that there is conservation of frustration here.

\paragraph{Reverse Process Initialization}

As in the flow-matching construction, exact reversal assumes
that the reverse process is initialized from the terminal
noising law $\mu_T^{\rm f}$. In practice, we may instead
initialize from another law $\nu_0$, for example a uniform
distribution on a certain domain.

Suppose these laws have densities
\[
\mu_T^{\rm f}(\dd x)=\rho_T^{\rm f}(x)\dd x,
\qquad
\nu_0(\dd x)=q_0(x)\dd x,
\]
and assume that
\begin{equation}
	q_0(x)\leq C\rho_T^{\rm f}(x)
\end{equation}
for some finite constant $C$. Applying the same change-of-measure
argument used for flow matching, the reverse law $\nu_t$
satisfies
\begin{equation}
	\nu_t(E)\leq C\mu_{T-t}^{\rm f}(E)
\end{equation}
for every measurable set $E$, provided the exact feedback
defines a well-posed flow.

For a point target $\mu_\star=\delta_{x_\star}$, take
$E=\{x:|x-x_\star|\geq r\}$, where $r>0$. Then
\begin{equation}
	\nu_t\bigl(\{x:|x-x_\star|\geq r\}\bigr)
	\leq
	C\mu_{T-t}^{\rm f}
	\bigl(\{x:|x-x_\star|\geq r\}\bigr)
	\longrightarrow0
\end{equation}
as $t\uparrow T$. Thus, changing the initialization does not
prevent convergence.

For the more
general control-affine diffusions, the support, smoothness,
and positivity conditions discussed in
Chapter~\ref{chap:fokker-planck} determine where such positivity and bounds
holds. If $\rho_T^{\rm f}$ is continuous and strictly positive
on a compact set containing the support of a bounded density
$q_0$, the required bound follows.

In particular, this justifies initialization from a small
distribution any configuration where the terminal noising
density is positive. A single prescribed initial point is
not covered directly by this change-of-measure argument,
since its law is a Dirac measure. However, in practice if one is completely at which exact point the process is initialized, this method can still be used as a heuristic.

All of this discussion, gives us the following algorithm.

\begin{tcolorbox}[algorithmbox,title=\textbf{Algorithm: Denoising Control for Systems with Drift}]
	\textbf{Setup.}
	Consider the control-affine system
	\[
	\dot{x}
	=
	f_0(x)+\sum_{i=1}^m u_i g_i(x).
	\]
	Choose a point target, take $x_\star$.
	Write $d_i=\nabla\cdot g_i$ and choose
	\[
	c_i=
	\begin{cases}
		0, & \text{uncorrected noising},\\
		d_i, & \text{divergence-corrected noising}.
	\end{cases}
	\]
	
	\begin{enumerate}[label=\arabic*.]
		\item Generate forward noising trajectories using
		\[
		\dd X_t^{\rm f}
		=
		\left(
		-f_0(X_t^{\rm f})
		+
		\sum_{i=1}^m c_i(X_t^{\rm f})g_i(X_t^{\rm f})
		\right)\dd t
		+
		\sqrt{2}\sum_{i=1}^m
		g_i(X_t^{\rm f})\circ\dd W_t^i,
		\qquad X_0^{\rm f} = x_\star.
		\]
		Denote their density by $\rho_t^{\rm f}$ and store samples
		\[
		\left\{X_{t_k}^{\mathrm f,(j)}\right\}_{j=1}^{N}
		\]
		at positive time steps $t_k\leq T$.
		
		\item Initialize a network
		\[
		s_\theta:[0,T]\times\Omega\to\R^m
		\]
		to approximate the horizontal score
		\[
		\left(
		\mathcal Y_1\log\rho_t^{\rm f},
		\ldots,
		\mathcal Y_m\log\rho_t^{\rm f}
		\right).
		\]
		
		\item Train the network by minimizing the empirical
		loss corresponding to
		\[
		\mathcal L_{\rm SM}(\theta)
		=
		\int_0^T
		\E_{X\sim\mu_t^{\rm f}}
		\left[
		\|s_\theta(t,X)\|^2
		+
		2\sum_{i=1}^m
		\left(
		\mathcal Y_i s_{\theta,i}(t,X)
		+
		d_i(X)s_{\theta,i}(t,X)
		\right)
		\right]\dd t,
		\]
		\item From the desired initial configuration $x_0$,
		apply the learned feedback
		\[
		u_{\theta,i}(t,x)
		=
		s_{\theta,i}(T-t,x)+d_i(x)-c_i(x)
		\]
		and integrate
		\[
		\dot{x}(t)
		=
		f_0(x(t))
		+
		\sum_{i=1}^m
		u_{\theta,i}(t,x(t))g_i(x(t)),
		\qquad x(0)=x_0,
		\qquad 0\leq t<T.
		\]
	\end{enumerate}
	
	With the exact horizontal score, the reverse law is $\mu_{T-t}^{\rm f}$
	and approaches $\mu_\star$ as $t\uparrow T$.
\end{tcolorbox}

\paragraph{Inverted-pendulum example.}

Consider the controlled pendulum
\begin{align}
    \dot\theta
    &=
    \omega,\\
    \dot\omega
    &=
    \alpha\sin\theta
    -
    \gamma\omega
    +u,
\end{align}
where $\alpha=g/\ell$ and $\gamma$ is a damping coefficient. The drift and control vector field are
\begin{equation}
    f_0(\theta,\omega)
    =
    \begin{pmatrix}
        \omega\\
        \alpha\sin\theta-\gamma\omega
    \end{pmatrix},
    \qquad
    g_1
    =
    \begin{pmatrix}
        0\\1
    \end{pmatrix}.
\end{equation}
Since $g_1$ is constant,
\begin{equation}
    \nabla\cdot g_1=0,
\end{equation}
so no divergence correction is needed. The forward noising process is therefore
\begin{align}
    \dd\theta_t
    &=
    -\omega_t\dd t,\\
    \dd\omega_t
    &=
    \left(
        -\alpha\sin\theta_t
        +
        \gamma\omega_t
    \right)\dd t
    +
    \sqrt{2}\,\dd W_t.
    \label{eq:pendulum-forward-noising}
\end{align}
Noise acts only through angular acceleration, while the reversed drift couples this excitation into the angle coordinate.

\begin{figure}[htbp]
	\centering
	\includegraphics[width=\textwidth]
	{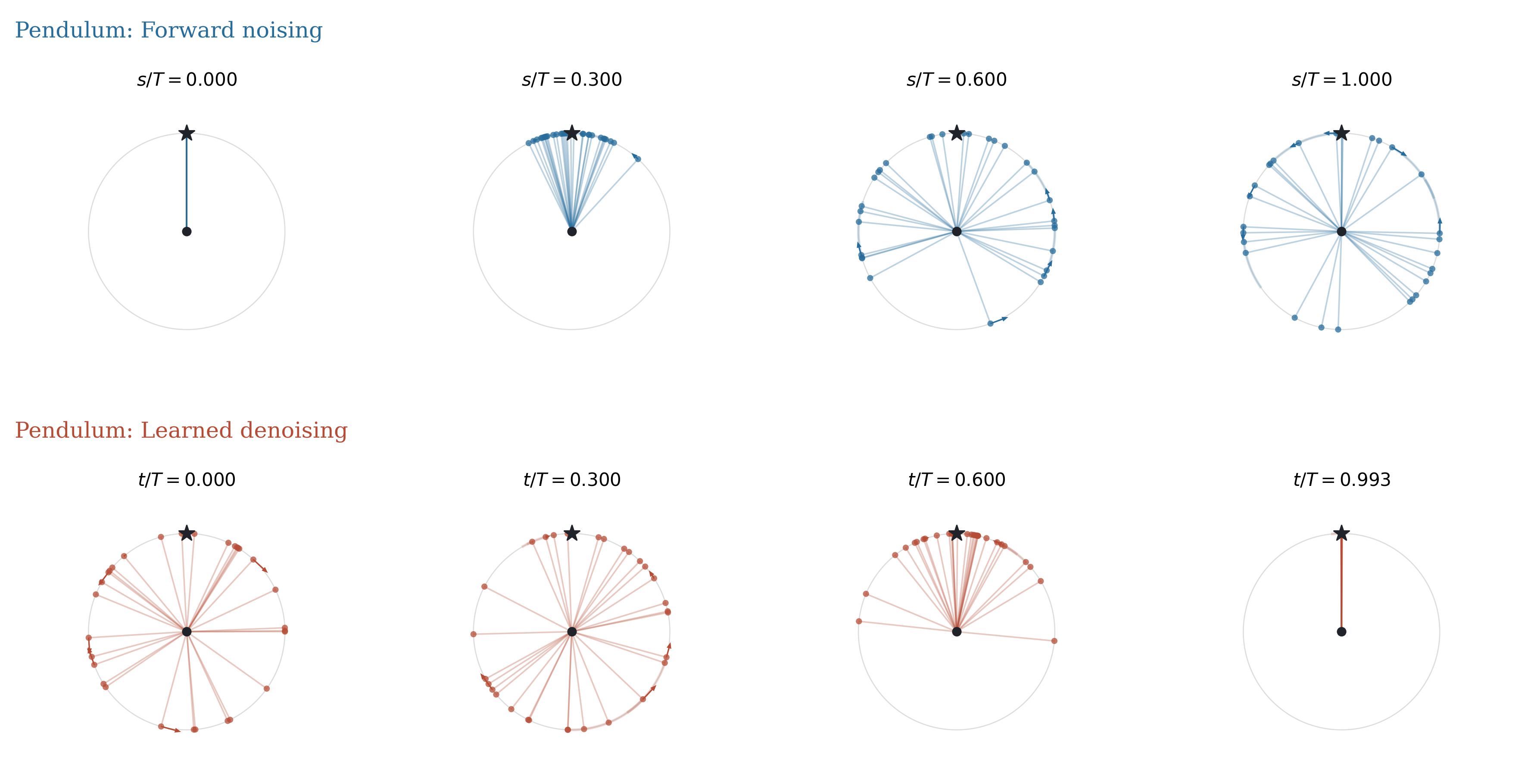}
	\caption{Noising and denoising for the inverted pendulum.
		Top: stochastic noising initialized at the upright
		equilibrium. Bottom: deterministic feedback using the
		learned horizontal score.
		Rods show orientations, arrows represent $\omega$. The star marks the upright target.}
	\label{fig:pendulum-noising-denoising}
\end{figure}

Figure~\ref{fig:pendulum-noising-denoising} illustrates the
construction for the inverted pendulum. In the noising phase,
particles start at the upright equilibrium
$(\theta,\omega)=(0,0)$ and spread under the time-reversed
drift and stochastic angular acceleration. A network learns
the horizontal score
$\mathcal Y_1\log\rho_t=\partial_\omega\log\rho_t$.
The deterministic denoising controller then uses this score
to steer independently sampled terminal noising states
toward the upright configuration.

Each rod represents a particle's orientation, while the
tangential arrows indicate its angular-velocity state $\omega$.
Thus, concentration of the rods near the upright position
shows angular alignment.

\hypertarget{score-matching-meets-flow-matching}{%
\section{Score Matching Meets Flow Matching}\label{score-matching-meets-flow-matching}}

In the preceding sections and chapters, we have seen
score matching and
flow matching as two
different but similar ideas that provide us with denoising controllers. Flow matching learns a velocity field as a
conditional expectation, whereas score matching learns the logarithmic
gradient of a density.

In the generative modeling community, it is known that the two methods
are closely related through \emph{Tweedie's formula}. One might even
find a very militantly pro-diffusion researcher scream this fact in
elation. Here, we will follow the exuberance of such a researcher.

Just as the flow-matched control is a conditional expectation, the
score-based controller can also be derived as a conditional expectation.
Instead of starting from Tweedie's formula, we will use an idea related
to \emph{Nelson's backward mean derivative}, though one will recognize
Tweedie's formula in some of the initial derivations.

We first consider Brownian motion to introduce the main idea and then
repeat the calculation for more general diffusions associated with
control systems.

Suppose the forward noising process is,

\[dX_t=\sqrt{2}\,dW_t,\]

and denote the density of \(X_t\) by \(\rho_t\) as usual. The score is

\[s(t,x)=\nabla\log\rho_t(x).\]

Let's define the \emph{backward conditional mean velocity}

\[b_-(t,x):=\lim_{h\downarrow0}\mathbb{E}\left[\left.\frac{X_t-X_{t-h}}{h}\right|X_t=x\right].\]

Let \(\phi\in C_c^\infty(\R^n)\) be a smooth test function. By the definition of conditional
expectation,

\[\int \phi(x)b_-(t,x)\rho_t(x)\,dx=\lim_{h\downarrow0}\frac{1}{h}\mathbb{E}\left[\phi(X_t)(X_t-X_{t-h})\right].\]

For Brownian motion,

\[X_t-X_{t-h}=\sqrt{2}(W_t-W_{t-h}).\]

{\it Stein's lemma} states that, if
$Z\sim\mathcal N(0,I_n)$ and $\phi:\R^n\to\R$ is sufficiently
regular, with the expectations, then
\[
\E[Z_i\phi(Z)]
=
\E[\partial_{z_i}\phi(Z)].
\]
Thus, multiplication by a Gaussian variable inside an
expectation can be exchanged for differentiation.
Applying this,
conditioning first on \(X_{t-h}\),

\[\mathbb{E}\left[\left.\phi(X_t)(X_t-X_{t-h})\right|X_{t-h}\right]=2h\,\mathbb{E}\left[\left.\nabla\phi(X_t)\right|X_{t-h}\right].\]

Taking expectation again and using the
tower property of iterated conditional expectations,

\[\mathbb{E}\left[\phi(X_t)(X_t-X_{t-h})\right]=2h\,\mathbb{E}\left[\nabla\phi(X_t)\right].\]

On the other hand, by conditioning on \(X_t\),

\[\mathbb{E}\left[\phi(X_t)\frac{X_t-X_{t-h}}{h}\right]=\mathbb{E}\left[\phi(X_t)\mathbb{E}\left[\left.\frac{X_t-X_{t-h}}{h}\right|X_t\right]\right].\]

Letting \(h\downarrow0\) and using the definition of \(b_-\),

\[\lim_{h\downarrow0}\mathbb{E}\left[\phi(X_t)\frac{X_t-X_{t-h}}{h}\right]=\mathbb{E}\left[\phi(X_t)b_-(t,X_t)\right].\]

Since \(X_t\) has density \(\rho_t\),

\[\mathbb{E}\left[\phi(X_t)b_-(t,X_t)\right]=\int \phi(x)b_-(t,x)\rho_t(x)\,dx.\]

Combining the two identities gives

\[\int \phi(x)b_-(t,x)\rho_t(x)\,dx=2\int \rho_t(x)\nabla\phi(x)\,dx.\]

Integrating by parts and assuming boundary terms vanish,

\[\int \phi(x)b_-(t,x)\rho_t(x)\,dx=-2\int \phi(x)\nabla\rho_t(x)\,dx.\]

Since this holds for every smooth test function \(\phi\),

\[b_-(t,x)\rho_t(x)=-2\nabla\rho_t(x).\]

Hence, wherever \(\rho_t>0\),

\[b_-(t,x)=-2\nabla\log\rho_t(x).\]

Thus the score can be written directly as a conditional expectation:

\[\nabla\log\rho_t(x)=-\frac{1}{2}b_-(t,x).\]

Equivalently,

\[\nabla\log\rho_t(x)=-\frac{1}{2}\lim_{h\downarrow0}\mathbb{E}\left[\left.\frac{X_t-X_{t-h}}{h}\right|X_t=x\right].\]

Since

\[X_t-X_{t-h}=\sqrt{2}(W_t-W_{t-h}),\]

we can also write

\[ \boxed{\nabla\log\rho_t(x)=-\frac{1}{\sqrt{2}}\lim_{h\downarrow0}\mathbb{E}\left[\left.\frac{W_t-W_{t-h}}{h}\right|X_t=x\right]}.\]

So the score measures the conditional mean noise that acted immediately
before the process arrived at \(x\). The right hand side of the formula is nearly identical  to the expressions for the control we derived in the flow matching chapter.

Our next goal is to generalize this heuristic derivation to more general
diffusions. Now consider the control-affine system with drift

\[\dot{x}=f_0(x)+\sum_{i=1}^m u_i g_i(x),\]

and the corresponding noising process

\[dX_t=-f_0(X_t)\,dt+\sqrt{2}\sum_{i=1}^m g_i(X_t)\circ dW_t^i.\]

Recall, 
\[d_i(x):=\nabla\cdot g_i(x),\qquad i=1,\ldots,m.\]

The generalized or horizontal score in the \(i\)-th control direction is

\[s_i(t,x):=\mathcal Y_i\log\rho_t(x).\]

We will also encounter the related adjoint quantity

\[u_i(t,x):=-\frac{\mathcal Y_i^\ast\rho_t(x)}{\rho_t(x)}.\]

As before, consider the backward conditional mean of the Brownian increment

\[\lim_{h\downarrow0}\mathbb{E}\left[\left.\frac{W_t^i-W_{t-h}^i}{h}\right|X_t=x\right].\]

Let \(\phi\) be a smooth test function. Conditioning on \(X_t\) gives

\[\mathbb{E}\left[\phi(X_t)\frac{\sqrt{2}(W_t^i-W_{t-h}^i)}{h}\right]=\mathbb{E}\left[\phi(X_t)\mathbb{E}\left[\left.\frac{\sqrt{2}(W_t^i-W_{t-h}^i)}{h}\right|X_t\right]\right].\]

We can compute the same quantity using Stein's lemma. Consider an Euler approximation of the noising
process in its {\it It\^o form}:
\[
X_t^h
=
X_{t-h}
+
h\,a(X_{t-h})
+
\sqrt{2}\sum_{j=1}^m
g_j(X_{t-h})\Delta W^j,
\]
where
\[
\Delta W^j=W_t^j-W_{t-h}^j,
\qquad
a=-f_0+\sum_{j=1}^m Dg_j\,g_j.
\]

Conditionally on $X_{t-h}$, $\Delta W^j$
are independent Gaussian variables with variance $h$.
Stein's lemma therefore gives
\begin{align*}
	\E\left[
	\phi(X_t^h)\frac{\sqrt{2}\Delta W^i}{h}
	\right]
	&=
	\sqrt{2}\,
	\E\left[
	\frac{\partial}{\partial\Delta W^i}\phi(X_t^h)
	\right]\\
	&=
	2\,\E\left[
	g_i(X_{t-h})\cdot\nabla\phi(X_t^h)
	\right].
\end{align*}
Formally passing to the limit $h\downarrow0$, we obtain
\[
\lim_{h\downarrow0}
\E\left[
\phi(X_t)\frac{\sqrt{2}(W_t^i-W_{t-h}^i)}{h}
\right]
=
2\,\E[\mathcal Y_i\phi(X_t)].
\]
Integration by parts yields
\[
2\,\E[\mathcal Y_i\phi(X_t)]
=
2\int_\Omega
\phi(x)\mathcal Y_i^*\rho_t(x)\,\dd x.
\]
Combining this with the conditional-expectation identity
identifies the limiting conditional mean, in the
test-function sense, as
\[
2\frac{\mathcal Y_i^*\rho_t}{\rho_t}
=
-2\left(
\nabla\cdot g_i+\mathcal Y_i\log\rho_t
\right).
\]

The Stein identity is exact for the Euler approximation. Passing to the continuous diffusion requires additional steps that would formalize the argument. Alternatively, the same
limiting relation can be derived directly for the diffusion
by applying It\^o's product rule to
$\phi(X_r)(W_r^i-W_{t-h}^i)$ over $r\in[t-h,t]$.

Therefore,

\[\lim_{h\downarrow0}\frac{\sqrt{2}}{h}\mathbb{E}\left[\phi(X_t)(W_t^i-W_{t-h}^i)\right]=2\int\rho_t(x)\mathcal Y_i\phi(x)\,dx.\]

Recall that the formal adjoint of \(\mathcal Y_i\) is

\[\mathcal Y_i^\ast f=-\mathcal Y_if-d_if.\]

By the definition of the adjoint,

\[\int\rho_t\mathcal Y_i\phi\,dx=\int\phi\,\mathcal Y_i^\ast\rho_t\,dx.\]

Combining this with the conditional-expectation calculation gives

\[\sqrt{2}\lim_{h\downarrow0}\mathbb{E}\left[\left.\frac{W_t^i-W_{t-h}^i}{h}\right|X_t=x\right]=2\frac{\mathcal Y_i^\ast\rho_t(x)}{\rho_t(x)}.\]

Therefore,

\[-\frac{\mathcal Y_i^\ast\rho_t(x)}{\rho_t(x)}=-\frac{1}{\sqrt{2}}\lim_{h\downarrow0}\mathbb{E}\left[\left.\frac{W_t^i-W_{t-h}^i}{h}\right|X_t=x\right].\]

The quantity on the left is exactly the feedback control that appears
in the time reversal derived in the
uncorrected noising process \eqref{eq:drift-first-noising-sde}.

Since

\[\mathcal Y_i^\ast\rho_t=-\mathcal Y_i\rho_t-d_i\rho_t,\]

we also have

\[-\frac{\mathcal Y_i^\ast\rho_t}{\rho_t}=\mathcal Y_i\log\rho_t+d_i=s_i+d_i.\]

Thus the conditional mean of the Brownian increment directly gives the
quantity

\[u_i^{\rm rev}(t,x)=-\frac{\mathcal Y_i^\ast\rho_{T-t}(x)}{\rho_{T-t}(x)}.\]

Equivalently,

\[u_i^{\rm rev}(t,x)=-\frac{1}{\sqrt{2}}\lim_{h\downarrow0}\mathbb{E}\left[\left.\frac{W_{T-t}^i-W_{T-t-h}^i}{h}\right|X_{T-t}=x\right].\]

So, for the original noising process, the conditional mean regression
recovers exactly the time-reversing controller we derived in the preceding sections.

Thus the time-reversing control itself is a conditional expectation of
the noise that acted immediately before the forward process arrived at
\(x\).

The preceding identity suggests that we do not need to learn the score
at all as we did in the
preceding sections. We can use the flow matching philosophy instead.

For a small finite \(h\), define 

\[U_{i,h}:=-\frac{W_t^i-W_{t-h}^i}{\sqrt{2}h}.\]

Take $t\in[h,T]$. If time is sampled during training, the conditioning below is on both $t$ and $X_t$; we display it at a fixed time for simplicity.

We can directly train a feedback controller \(u_\theta\) using the
ordinary regression loss

\[\mathcal{L}_{\mathrm{FM}}(\theta)=\mathbb{E}\left[\sum_{i=1}^m\left|u_{\theta,i}(T-t,X_t)-U_{i,h}\right|^2\right].\]

By the usual \(L^2\) regression property, the infinite sample minimizer
satisfies

\[u_{\theta,i}^\star(T-t,x)=\mathbb{E}\left[U_{i,h}\mid X_t=x\right].\]

Therefore, formally as \(h\downarrow0\),

\[u_{\theta,i}^\star(T-t,x)\longrightarrow-\frac{\mathcal Y_i^\ast\rho_t(x)}{\rho_t(x)}.\]

So the same feedback controller that we previously obtained through
score matching can instead be learned directly by flow matching.

\hypertarget{score-matching-vs-flow-matching}{%
\subsubsection{Score Matching vs Flow
Matching}\label{score-matching-vs-flow-matching}}

We therefore have two ways of learning the same denoising controller.

\begin{enumerate}
\def\labelenumi{\arabic{enumi}.}

\item
  \textbf{Score matching:} learn
\end{enumerate}

\[s_i(t,x)=\mathcal Y_i\log\rho_t(x),\]

and then construct the controller

\[u_i^{\rm rev}(t,x)=s_i(T-t,x)+d_i(x)=-\frac{\mathcal Y_i^\ast\rho_{T-t}(x)}{\rho_{T-t}(x)}.\]

\begin{enumerate}
\def\labelenumi{\arabic{enumi}.}
\setcounter{enumi}{1}

\item
  \textbf{Flow matching:} directly regress the random control
\end{enumerate}

\[U_{i,h}=-\frac{W_t^i-W_{t-h}^i}{\sqrt{2}h}\]

Unlike implicit score matching, this regression does not require
differentiating the neural network with respect to the state, and in our
numerical experiments this makes training considerably faster while
producing similar denoising behavior. Lastly, note that 
Brownian paths themselves do not have ordinary velocities, since they are almost nowhere differentiable. It is their conditional increment mean that is matched here. Though, informally it is tempting to think of the quantity we are trying to learn as,
$$-\frac{1}{\sqrt{2}}\mathbb{E}\left[\left.\dot{W_t}\right|X_{T-t}=x\right]. $$



\paragraph{Further reading.}
The denoising diffusion approach for steering to target distributions is from \cite{elamvazhuthi2025score}. While we have considered time-reversals which are deterministic, one can also consider time-reversals that are stochastic \cite{mei2025time}, which are more canonical in some sense due to the classical works in diffusion processes on time-reversals \cite{anderson1982reverse,haussmann1986time}.
The presentation on backward mean derivative, to relate score matching and flow matching, is inspired by the exposition of Yuri Gliklikh \cite{gliklikh2010global,grasse2017global}.

\chapter{Control Using Sub-Laplacians}
\label{chap:sub-laplacian-control}

Chapter~\ref{chap:denoising-diffusion} constructs transport by first generating a diffusion and then reversing its probability flow. This chapter gives a deterministic alternative. We prescribe the density path directly and solve an PDE for a feedback law. The same horizontal differential operators introduced in Chapter~\ref{chap:fokker-planck} reappear naturally and play an important role. 

Particularly, we consider the \textbf{sub-Laplacian}, also called the \textbf{sub-Riemannian} or \textbf{horizontal Laplacian}, a natural second-order operator associated with a family of admissible control vector fields. Just as the ordinary Laplacian encodes the geometry of a connected domain, the sub-Laplacian encodes geometry generated by the control directions. For our purposes, it can also be used to transport one probability density to another. This produces a route to control that is different from the stochastic noising--denoising constructions developed above.

\section{Moser Transport}

Let $\rho_0$ and $\rho_1$ be probability densities on a bounded connected domain $\Omega$. Our goal is to construct a vector-field transporting one probability density to another. As usual, we will consider the fully actuated situation to motivate the developments. Consider the Neumann Poisson problem
\begin{align}
    \Delta\phi
    &=
    \rho_1-\rho_0
    \qquad\text{in }\Omega,\\
    n_{\rm out}\cdot\nabla\phi
    &=0
    \qquad\text{on }\partial\Omega.
    \label{eq:moser-poisson}
\end{align}
where $\Delta$ is the usual Laplacian operator.
Since
\begin{equation}
    \int_\Omega(\rho_1-\rho_0)\dd x=0,
\end{equation}
on a connected domain, the solution is unique up to an additive constant under the usual regularity assumptions.

Now consider the {\it linear interpolation} between densities
\begin{equation}
    \rho_t
    =
    (1-t)\rho_0+t\rho_1,
    \qquad
    t\in[0,1].
    \label{eq:moser-linear-interpolation}
\end{equation}
Then
\begin{equation}
    \partial_t\rho_t
    =
    \rho_1-\rho_0.
\end{equation}
Assuming we have solved the Neumann problem and choosing
\begin{equation}
    v(t,x)
    =
    -\frac{\nabla\varphi(x)}{\rho_t(x)}.
    \label{eq:moser-velocity}
\end{equation}
we substitute into the continuity equation 
\begin{align}
    \partial_t\rho_t
    &=
    -\nabla\cdot(\rho_t v)\\
    &=
    \nabla\cdot(\nabla\phi)\\
    &=
    \Delta\phi\\
    &=
    \rho_1-\rho_0.
\end{align}
Thus the ODE 
\begin{equation}
	\dot x=v(t,x)= -\frac{\nabla\varphi(x)}{\rho_t(x)}
	\end{equation}
	transports $\rho_0$ to $\rho_1$.

This way of interpolating between densities is known as the {\it Moser transport}. 

\section{Nonholonomic Moser Transport}

The usual Moser transport is applicable to  fully actuated systems. As in Chapters~\ref{chap:fokker-planck} and~\ref{chap:denoising-diffusion}, we now replace ordinary derivatives by derivatives along the admissible control vector fields.

Consider the driftless control system
\begin{equation}
    \dot x
    =
    \sum_{i=1}^m u_i(t,x)g_i(x).
    \label{eq:moser-control-system}
\end{equation}
Recall the operators $\mathcal Y_i$ and their formal adjoints $\mathcal Y_i^*$ from \eqref{eq:global-directional-operators}. Define the sub-Laplacian
\begin{equation}
    \boxed{
    \Delta_H h
    :=
    -\sum_{i=1}^m\mathcal Y_i^*\mathcal Y_i h
    =
    \nabla\cdot
    \left(
        \sum_{i=1}^m
        g_i g_i^{\mathsf T}\nabla h
    \right).
    }
    \label{eq:horizontal-laplacian}
\end{equation}
We can see that for coordinate vector fields $g_i=e_i$, this reduces to the ordinary Euclidean Laplacian.

Use the same density interpolation \eqref{eq:moser-linear-interpolation} and solve the horizontal Poisson equation
\begin{equation}
    \Delta_H\phi
    =
    \rho_1-\rho_0
    \qquad\text{in }\Omega.
    \label{eq:horizontal-poisson}
\end{equation}
The natural no-flux boundary condition is
\begin{equation}
    n_{\rm out}\cdot
    \left(
        \sum_{i=1}^m
        g_i\mathcal Y_i\phi
    \right)
    =0
    \qquad\text{on }\partial\Omega.
    \label{eq:horizontal-neumann}
\end{equation}

Now choose the feedback law
\begin{equation}
    \boxed{
    u_i(t,x)
    =
    -\frac{\mathcal Y_i\phi(x)}{\rho_t(x)}.
    }
    \label{eq:nonholonomic-moser-control}
\end{equation}
The continuity equation for \eqref{eq:moser-control-system} is
\begin{equation}
    \partial_t\rho_t
    =
    \sum_{i=1}^m
    \mathcal Y_i^*\bigl(u_i\rho_t\bigr).
\end{equation}
Substituting \eqref{eq:nonholonomic-moser-control} gives
\begin{align}
    \partial_t\rho_t
    &=
    -\sum_{i=1}^m
    \mathcal Y_i^*\mathcal Y_i\phi\\
    &=
    \Delta_H\phi\\
    &=
    \rho_1-\rho_0.
\end{align}
Hence the control law realizes the desired density interpolation exactly, as long as the Poisson equation is solvable and $\rho_t$ stays positive where the control is evaluated.

\paragraph{Tracking Arbitrary Probability Densities}

While we started with linear interpolation, there is no reason to restrict ourselves to such trajectories. Suppose, more generally, we have a family of differentiable probability densities $(\rho_t)_{t \in [0,1]}$ such that $\rho_t >0$ are strictly positive densities for all $t$. Then mimicking the Moser transport we can consider exactly the same idea to solve the Poisson equation for every time $t>0$,

\begin{equation}
	\Delta_H\phi
	=
	 \partial_t \rho_t 
	\qquad\text{in }\Omega.
	\label{eq:horizontal-poisson2}
\end{equation}
and then once again we get a control law in \eqref{eq:nonholonomic-moser-control}, except we need to solve the Poisson equation for every $t \in [0,1]$. If $\rho_t$ integrate to $1$, $\frac{d}{dt} \int_{\Omega} \rho_t =  \int_{\Omega} \partial_t \rho_t = 0$, and the Poisson equation again has a unique solution up to additive constants.

Remarkably, if the control system is controllable, it is {\it trajectory controllable} along any family of sufficiently regular densities.

\paragraph{Numerical example.}

Consider the Brockett
integrator,
\begin{align}
	\dot x_1 &= u_1,\\
	\dot x_2 &= u_2,\\
	\dot x_3 &= x_1u_2-x_2u_1,
\end{align}
on the domain $\Omega=(-1,1)^3$. Its control vector fields are
\[
g_1=(1,0,-x_2)^{\mathsf T},
\qquad
g_2=(0,1,x_1)^{\mathsf T}.
\]
Together with their Lie bracket
$[g_1,g_2]=(0,0,2)^{\mathsf T}$, these fields span
$\R^3$ at every point.

We choose initial and terminal densities
\[
\rho_j(x)
=
\frac{1}{Z_j}
\exp\!\left(-\frac{|x|^2}{2\sigma_j^2}\right),
\qquad
Z_j
=
\int_\Omega
\exp\!\left(-\frac{|x|^2}{2\sigma_j^2}\right)\dd x,
\qquad j=0,1,
\]
with $\sigma_0=0.45$ and $\sigma_1=0.05$. Thus, the objective
is to transport a broad initial distribution to a density
concentrated near the origin. We prescribe the linear
interpolation
\[
\rho_t=(1-t)\rho_0+t\rho_1,
\qquad t\in[0,1].
\]

We solve the horizontal Poisson equation in a
\emph{physics-informed neural network} (PINN) fashion.
Specifically, we approximate the potential by a scalar
network $\phi_\theta(t,x)$ and minimize
\begin{align}
	\mathcal L(\theta)
	&=
	\mathbb E_{t,x}
	\left[
	\left|
	\Delta_H\phi_\theta(t,x)
	-(\rho_1(x)-\rho_0(x))
	\right|^2
	\right]\\
	&\quad+
	\lambda_{\mathrm{BC}}\,
	\mathbb E_{t,x_b}
	\left[
	\left|
	n_{\rm out}(x_b)\cdot
	\sum_{i=1}^{2}
	g_i(x_b)\mathcal Y_i\phi_\theta(t,x_b)
	\right|^2
	\right].
\end{align}
Here, $t$ is sampled uniformly from $[0,1]$, while $x$ and
$x_b$ represent interior and boundary collocation points.
The expectations correspond to uniform volume and surface
sampling, respectively, and $\lambda_{\mathrm{BC}}=10$.
Spatial derivatives are evaluated using automatic
differentiation. We fix the additive constant by enforcing
$\phi_\theta(t,0)=0$.

Although the network accepts time as an input, the exact
Poisson solution for this interpolation can be chosen
independent of time. The feedback remains time dependent
through the prescribed density:
\[
u_i^\theta(t,x)
=
-\frac{\mathcal Y_i\phi_\theta(t,x)}{\rho_t(x)},
\qquad i=1,2.
\]
After training, we select the parameters with the lowest
validation loss, sample initial states from $\rho_0$, and
integrate the controlled system.

\begin{figure}[htbp]
	\centering
	\includegraphics[width=\textwidth]
	{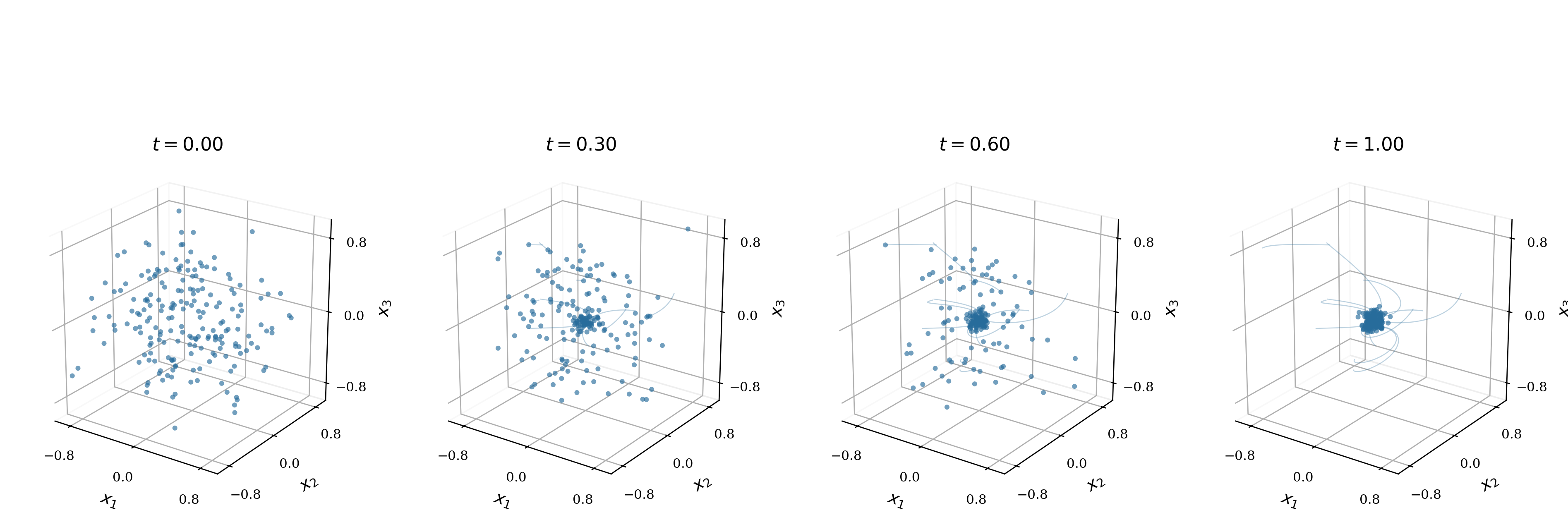}
	\caption{Nonholonomic Moser transport for the Brockett
		integrator using linear density interpolation.}
	\label{fig:nonholonomic-moser-snapshots}
\end{figure}

Figure~\ref{fig:nonholonomic-moser-snapshots} shows the
resulting evolution. The prescribed interpolation
concentrates the distribution near the origin, while the
particle trajectories must follow the admissible directions
of the nonholonomic system. For stabilization purposes, one should use this cautiously, as the inference generated $200$ samples there were $15$ samples that exited the domain (not seen in the plots). This indicates that maybe a better neural net or some additional fine-tuning is required is one is to use this for stabilizaiton objectives.

\paragraph{Computational limitations.}

The attractive feature of the method is that a nonlinear control problem has been converted into a linear PDE. The drawback is that one still has to solve a PDE on state space. Consequently, the method retains the familiar curse-of-dimensionality problem of grid-based HJB methods. Neural-network approximations can remove the grid, but then one introduces a sample complexity problem arising form the chosen samples at which the loss function is evaluated. Additionally, we now have a nonlinear, nonconvex approximation problem at the numerical level.

\section{Reachability to Transportability}

A theme that already appeared in
Chapters~\ref{chap:fokker-planck}
and~\ref{chap:denoising-diffusion} was that controllability
is closely related to whether derivatives along the
admissible directions are rich enough to make the
corresponding diffusions sufficiently non-degenerate,
and hence the construction of the associated controllers
well-posed. For driftless systems with controls of either
sign, reachability is symmetric: if one state can be
reached from another, the trajectory can also be traversed
in reverse.

The connection with controllability is analogous to the
connection between the ordinary Neumann Laplacian and
connectivity. Under suitable regularity and boundary
assumptions, bracket-generating controllability makes
the zero eigenvalue of the sub-Laplacian simple.
Together with the horizontal Poincar\'e inequality
discussed earlier, this ensures that, for a zero-mean
right-hand side such as $\rho_1-\rho_0$, the horizontal
Poisson equation has a solution unique up to an additive
constant.

Hence, with sufficient regularity of the resulting
feedback, the Moser transport of the previous section
is well-posed, and the system is {\it transportable}
from a smooth strictly positive initial density to
a smooth strictly positive target density. Therefore,
under these assumptions, we have the satisfying relation:
{\it controllability implies transportability}.

\section{Does Transportability Imply Reachability?}

Given that solvability of the Poisson equation can be
derived from controllability of the system, one can
ask if the converse is true. Does solvability of the
Poisson equation provide us with a way to conclude
controllability of the associated driftless control
system? The answer is: {\it Almost}. We can take two
{\it dual} views on this.

Throughout this discussion, trajectories are required
to remain in $\Omega$, and controls of either sign
are allowed. Reachable sets refer to states reachable
in some finite time.

\subsection{Perron--Frobenius View}

The Perron--Frobenius perspective is the idea that the
behavior of a dynamical system can be understood by
how a probability distribution evolves under its
action. Our perspective throughout these notes has
largely been the Perron--Frobenius perspective applied
to control problems.

In this context, we can ask: suppose we can achieve
Moser transport, and the Poisson equation
\eqref{eq:horizontal-poisson} is solvable for every
smooth strictly positive initial and target density
$\rho_0,\rho_1$. Is the associated system also
controllable? Here, we assume that the resulting
transport can be represented by admissible trajectories
of the control system. This problem would be trivial
if we could use Moser transport for Dirac initial and
final distributions, $\delta_{x_0}$ and $\delta_{x_1}$.
However, the Poission equation based construction requires the interpolating
densities to be strictly positive.

\begin{figure}[htbp]
	\centering
	\includegraphics[width=\textwidth]
	{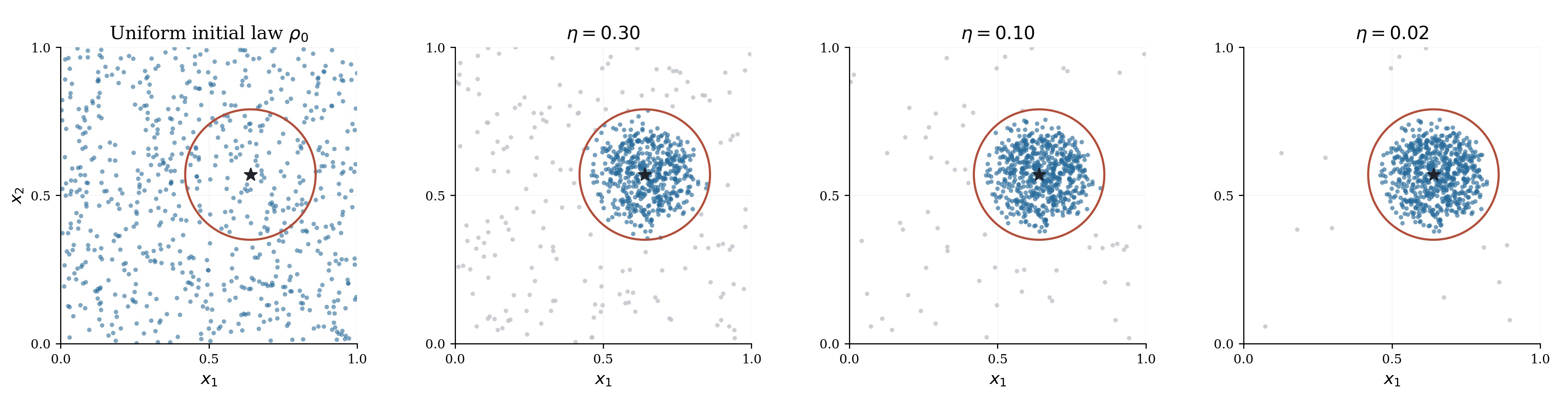}
	\caption{Schematic of Perron-Frobenius argument. The target ball in red is kept fixed. The parameter $\eta$ controls the volume of points transported to the ball. Dark blue shows the assigned endpoint inside ball. Light blue shows the assigned endpoint outside ball
		Target ball. }
	\label{fig:PFsnapshots}
\end{figure}

We can construct a {\it proof by pictures}. Take
$\rho_0=1/|\Omega|$ to be the uniform density. Fix a
ball $B_\varepsilon(x_1)\subset\Omega$, and choose a
smooth strictly positive target density
$\rho_1^{\varepsilon,\eta}$ such that
\[
\int_{B_\varepsilon(x_1)}
\rho_1^{\varepsilon,\eta}(x)\dd x
>
1-\eta,
\qquad \eta>0.
\]
For every $\varepsilon>0$ and every $\eta>0$, we can
transport the uniform distribution to this target
density, and hence more than a fraction $1-\eta$
of the mass into the $\varepsilon$-ball around $x_1$. See \ref{fig:PFsnapshots} for a schematic.

This means that a set of initial states occupying
more than a fraction $1-\eta$ of the domain can be
transported into that ball. Since $\eta$ can be made
arbitrarily small while keeping the ball fixed,
the set of states that can reach the ball has full
Lebesgue measure in $\Omega$.

By reversibility of the driftless system, this is
also the set of states reachable from the ball.
Writing $\mathcal R(B_\varepsilon(x_1))$ for this
reachable set, we obtain $
\left|
\Omega\setminus\mathcal R(B_\varepsilon(x_1))
\right|
=0.$
Technically speaking:

\begin{quote}
	{\it The set of reachable states from every
		$\varepsilon$-ball has full Lebesgue measure in
		$\Omega$.}
\end{quote}

One can think of this as a kind of
{\it approximate controllability}, since
$\varepsilon>0$ can be made arbitrarily small.
The final implication can be refined further. For a fixed target $x_1$, taking a countable sequence
$\varepsilon_k\downarrow0$ shows that almost every
initial state can be steered arbitrarily close to
$x_1$. This, of course, does not yet say that every initial state
can reach $x_1$ exactly.

\subsection{Koopman View}

One can also take the {\it Koopman perspective} on
this. The Koopman perspective looks at how the behavior
of {\it observables}, or functions evaluated along
trajectories, can be used to infer properties of
a dynamical system.

To see again why uniqueness for the Poisson equation implies
the same $\varepsilon$-ball kind of controllability,
let us suppose that $
A=\mathcal R(B_\varepsilon(x_1))$
does not have full measure. Then $A$ has positive
volume, since it contains the ball, and its complement
$B=\Omega\setminus A$ also has positive volume.

No state in $B$ can be reached from $A$. Otherwise,
concatenating trajectories would put that state in
$A$. By reversibility, no state in $A$ can be reached
from $B$ either. Thus, both sets are {\it invariant}
under admissible trajectories of the control system.

Consider their indicator functions,
\[
1_A(x)=
\begin{cases}
	1, & x\in A,\\
	0, & x\notin A,
\end{cases}
\qquad
1_B(x)=
\begin{cases}
	1, & x\in B,\\
	0, & x\notin B.
\end{cases}
\]
Since the sets are invariant under every admissible
control, they are also invariant under the individual
vector-field flows
\[
\dot x=g_i(x),
\]
for as long as these flows remain in $\Omega$.
For the control system
\eqref{eq:moser-control-system}, this corresponds
to setting $u_i\equiv1$ and all other controls
equal to zero.

Denoting the flow of $g_i$ by $\psi_t^i$, we have
\[
1_C(\psi_t^i(x))=1_C(x),
\qquad C\in\{A,B\}.
\]
Formally differentiating gives
\[
\frac{\dd}{\dd t}1_C(\psi_t^i(x))
=
(\mathcal Y_i1_C)(\psi_t^i(x))
=
0.
\]
Thus, the indicators have zero derivatives along
all the admissible directions. This also suggests
that
\[
\Delta_H1_C
=
-\sum_{i=1}^m
\mathcal Y_i^*\mathcal Y_i1_C
=
0.
\]

Here, we are ignoring the fact that $1_C$ need not
be differentiable in the ordinary sense. But the
idea can be made rigorous.

Now define
\[
f(x)=1_A(x)-1_B(x).
\]
Then
\[
\Delta_Hf=0.
\]
Since both $A$ and $B$ have positive volume, $f$
is not almost everywhere constant. This contradicts
uniqueness of the Poisson solution up to additive
constants. Equivalently, subtracting the spatial
mean gives
\[
\widetilde f
=
f-\frac{1}{|\Omega|}\int_\Omega f(x)\dd x,
\]
which is a nonzero, zero-mean solution of
$\Delta_H\widetilde f=0$. If uniqueness held after
fixing the mean, the only such solution would be
zero.

Thus, the domain cannot split into two invariant
pieces of positive volume. In particular, the
reachable set of every $\varepsilon$-ball must
have full measure, giving the same conclusion
as the Perron--Frobenius argument.
\paragraph{Further reading.}

The nonholonomic Moser construction  is due to \cite{khesin2009nonholonomic}. See also related geometric work on sub-Riemannian structures on groups of diffeomorphisms \cite{arguillere2017sub}. The relation between transport and controllability, as presented in the last section, has been derived in \cite{elamvazhuthi2025linear}. Classical connections between second-order operators and controllability include work of H\"ormander \cite{hormander1967hypoelliptic}, and support theorems of Stroock and Varadhan \cite{stroock1972support}.

\bibliography{ref}
\bibliographystyle{plain}

\end{document}